\documentclass[11pt]{article}

\usepackage[latin1]{inputenc} 
\usepackage{amsthm,amsmath,amssymb} 
\usepackage{graphicx}
\usepackage{color}
\usepackage{verbatim}

\usepackage{enumitem} 
\setlist[enumerate]{topsep=0pt,itemsep=-1ex,partopsep=1ex,parsep=1ex}
\renewcommand{\theenumi}{\roman{enumi}}

\usepackage{sidecap}

\usepackage{tikz}
\usetikzlibrary{shapes,positioning,intersections,quotes}

\theoremstyle{plain}
\newtheorem{theo}{Theorem}[section]

\newtheorem{lemma}[theo]{Lemma}
\newtheorem{cor}[theo]{Corollary}

\theoremstyle{definition}
\newtheorem{defn}[theo]{Definition}
\newtheorem{rem}[theo]{Remark}

\newtheorem{alg}[theo]{Algorithm}

\newcommand{\mc}[1]{\mathcal{#1}}
\newcommand{\mb}[1]{\mathbb{#1}}
\newcommand{\nib}[1]{\noindent {\bf #1}}
\newcommand{\nim}[1]{\noindent {\em #1}}
\newcommand{\brac}[1]{\left( #1 \right)}

\newcommand{\bsize}[1]{\left| #1 \right|}

\newcommand{\bgen}[1]{\left\langle #1 \right\rangle}
\newcommand{\sgen}[1]{\langle #1 \rangle}

\newcommand{\sub}{\subseteq}

\newcommand{\Lra}{\Leftrightarrow}

\newcommand{\sm}{\setminus}

\newcommand{\es}{\emptyset}
\newcommand{\pl}{\partial}

\newcommand{\aA}{\alpha}
\newcommand{\bB}{\beta}
\newcommand{\gG}{\gamma}

\newcommand{\lL}{\lambda}
\newcommand{\tT}{\theta}
\newcommand{\sS}{\sigma}
\newcommand{\oO}{\omega}

\newcommand{\GG}{\Gamma}

\newcommand{\OO}{\Omega}

\newcommand{\Ups}{\Upsilon}

\newcommand{\rk}{\operatorname{rank}}

\def\qed{\hfill $\Box$}

\title{The existence of designs}

\author{Peter Keevash\thanks{Mathematical Institute,
University of Oxford, Oxford, UK. Email: keevash@maths.ox.ac.uk.
\newline \hspace*{1.8em}Research supported
in part by ERC Consolidator Grant 647678 and ERC Advanced Grant 883810.}}

\begin{document}

\maketitle

\begin{abstract}
We prove the existence conjecture for combinatorial designs,
answering a question of Steiner from 1853.
More generally, we show that the natural divisibility conditions
are sufficient for clique decompositions of uniform hypergraphs
that satisfy a certain pseudorandomness condition. 
As a further generalisation, we obtain the same conclusion
only assuming an extendability property and the existence
of a robust fractional clique decomposition.
\end{abstract}

\section{Introduction}

A \emph{Steiner system} with parameters $(n,q,r)$ is a set $S$ 
of $q$-subsets of an  $n$-set\footnote{This means that $|X|=n$ and $S$ 
consists of subsets of $X$ each having size $q$.} $X$ such that 
every $r$-subset of $X$ is contained in exactly one element of $S$.
The question of whether there is a Steiner system with given 
parameters is one of the oldest problems in combinatorics, 
dating back to work of Pl\"ucker (1835), Kirkman (1846) 
and Steiner (1853); see \cite{RobinW} for a historical account. 

More generally, we say that a set $S$ of $q$-subsets of an 
$n$-set $X$ is a \emph{design} with parameters $(n,q,r,\lL)$ if 
every $r$-subset of $X$ is contained in exactly $\lL$ elements of $S$.
(This is often called an `$r$-design' in the literature.)
There are some obvious necessary `divisibility conditions' 
for the existence of such $S$, namely that $\tbinom{q-i}{r-i}$ 
divides $\lL \tbinom{n-i}{r-i}$ for every $0 \le i \le r$ 
(fix any $i$-subset $I$ of $X$ and consider the sets in $S$ 
that contain $I$). It is not known who first advanced the 
`Existence Conjecture' that the divisibility conditions are 
also sufficient, apart from a finite number of exceptional $n$ 
given fixed $q$, $r$ and $\lL$.

The case $r=2$ has received particular attention due to 
its connections to statistics, under the name of 
`balanced incomplete block designs'. 
We refer the reader to \cite{CD} for a summary
of the large literature and applications of this field. 
The Existence Conjecture for $r=2$ was a long-standing
open problem, eventually resolved by Wilson \cite{W1,W2,W3} 
in a series of papers that revolutionised Design Theory,
and had a major impact in Combinatorics. In this paper, 
we prove the Existence Conjecture in general, via a new method, 
which we will refer to as Randomised Algebraic Constructions.

\subsection{Results}

The Existence Conjecture will follow from a more general result 
on clique decompositions of hypergraphs that satisfy a certain 
pseudorandomness condition. 
To describe this we make the following definitions. 

\begin{defn}
A \emph{hypergraph} $G$ consists of a vertex set $V(G)$ and an 
edge set $E(G)$, where each $e \in E(G)$ is a subset of $V(G)$. 
We identify $G$ with $E(G)$.%
\footnote{So $|G|=|E(G)|$. We stress this point,
as most authors instead write $|G|=|V(G)|$.} 
If every edge has size $r$ we say that $G$ is an \emph{$r$-graph}. 
For $S \sub V(G)$, the \emph{neighbourhood} $G(S)$ 
is the $(r-|S|)$-graph $\{f \sub V(G) \sm S: f \cup S \in G\}$.
For an $r$-graph $H$, an \emph{$H$-decomposition} of $G$ 
is a partition of $E(G)$ into subgraphs isomorphic to $H$.
Let $K^r_q$ be the complete $r$-graph on $q$ vertices.
\end{defn}

Note that a Steiner system with parameters $(n,q,r)$
is equivalent to a $K^r_q$-decomposition of $K^r_n$.
It is also equivalent to a perfect matching 
(a set of edges covering every vertex exactly once)
in the auxiliary $\tbinom{q}{r}$-graph on $\tbinom{[n]}{r}$
(the $r$-subsets of $[n]:=\{1,\dots,n\}$) with edge set
$\{ \tbinom{Q}{r}: Q \in \tbinom{[n]}{q} \}$.
The next definition generalises the necessary
divisibility conditions described above.%
\footnote{Note that $|G(e)|$ denotes the number of edges in 
the neighbourhood of $e$, which is the degree of $e$ in $G$.}

\begin{defn} \label{Gdiv}
Suppose $G$ is an $r$-graph. We say that $G$ is 
\emph{$K^r_q$-divisible} if $\tbinom{q-i}{r-i}$ divides $|G(e)|$ 
for any $i$-set $e \sub V(G)$, for all $0 \le i \le r$. 
\end{defn}

Next we formulate our quasirandomness condition.
It is easy to see that it holds
with high probability (w.h.p.) if $G = G^r(n,p)$ 
is the standard binomial random $r$-graph
and $n$ is large given $p$, $c$ and $h$.

\begin{defn} \label{def:typ}
Suppose $G$ is an $r$-graph on $[n]$.
The density of $G$ is $d(G) = |G| \tbinom{n}{r}^{-1}$.
We say that $G$ is \emph{$(c,h)$-typical} if for any set $A$ 
of $(r-1)$-subsets of $V(G)$ with $|A| \le h$ we have 
$\bsize{\cap_{S \in A} G(S)} = (1 \pm |A|c) d(G)^{|A|} n$.
\end{defn}

Now we can state a simplified form of our main theorem.

\begin{theo} \label{main}
For any $q > r \ge 1$ there are $c_0, \aA >0$ 
and $h,n_0 \in \mb{N}$ such that if $G$ is a $K^r_q$-divisible 
$(c,h)$-typical $r$-graph on $n>n_0$ vertices,
where $d(G)>n^{-\aA}$ and $c < c_0 d(G)^{h^2}$,
then $G$ has a $K^r_q$-decomposition.
\end{theo}

Applying this with $G=K^r_n$, 
we deduce that for large $n$ the divisibility conditions
are sufficient for the existence of Steiner systems;
the existence of designs with any constant multiplicity $\lL$
follows from Theorem \ref{main+} below.
We have not tried to optimise our parameters,
although we do emphasise that the density of $G$
can decay polynomially in $n$, as this is used 
in \cite{K2} to estimate the number of designs.
Our method also gives a randomised algorithm 
for constructing designs.

Theorem \ref{main} gives new results 
even in the graph case ($r=2$); for example, 
it is easy to deduce that the standard random graph model $G(n,1/2)$
w.h.p.\ has a partial triangle decomposition that covers all 
but $(1+o(1))n/4$ edges:
deleting a perfect matching on the set of vertices of odd degree
and then at most two $4$-cycles
(to make the number of edges divisible by $3$)
gives a graph satisfying the hypotheses of the theorem.
This is the asymptotically best possible `leave', 
as w.h.p.\ there are $(1+o(1))n/2$ vertices of odd degree 
and any partial triangle decomposition must leave 
at least one edge uncovered at each vertex of odd degree.

We also note that if an $r$-graph $G$ on $n$ vertices 
satisfies $|G(S)| \ge (1-c)n$ for every $(r-1)$-subset $S$ 
of $V(G)$ then it is $(c,h)$-typical,
so we also deduce a minimum $(r-1)$-degree version of the theorem,
generalising Gustavsson's minimum degree version \cite{Gu} 
of Wilson's theorem.

To state our main theorem we introduce 
the following more general context of $r$-multigraphs.
Note that an $(n,q,r,\lL)$-design is equivalent
to a $K^r_q$-decomposition of the
$r$-multigraph $\lL \tbinom{[n]}{r}$.

\begin{defn}
An \emph{$r$-multigraph} $G$ on $[n]$ is a multiset
in which each element is an $r$-subset of $[n]$.
We identify $G$ with a vector%
\footnote{We identify $K^r_n$ with its 
set of edges $\tbinom{[n]}{r}$.} 
$G \in \mb{N}^{K^r_n}$,
where $G_e$ is the multiplicity of $e$ in $G$.
\end{defn}

We can also relax our pseudorandomness assumption,
with essentially the same proof,
obtaining a more general result in the spirit of \cite{KM},
namely that under certain conditions 
(`extendability' and `robust fractional decomposition'),
divisibility is the only obstruction to decomposition.
The next two definitions formulate our extendability assumption
(see Subsection \ref{sec:ext} for more discussion).

\begin{defn} \label{def:embed}
Suppose $H$ is an $r$-graph,
$G$ is an $r$-multigraph on $[n]$
and $\phi:V(H) \to [n]$ is injective.
We call $\phi$ an \emph{embedding} of $H$ in $G$
if $G_{\phi(f)}>0$ for all $f \in H$.
We write $K^r_q(G)$ for the set%
\footnote{We regard cliques as the same 
if they are identical as a subset of $K^r_n$:
we do not distinguish multiple edges.} 
of $\phi(Q)$ where $\phi$ 
is an embedding of $Q=K^r_q$ in $G$.
\end{defn}

\begin{defn} \label{def:ext}
Suppose $H$ is an $r$-graph with no isolated vertices,
$F \sub V(H)$ and $\phi: F \to [n]$ is injective.
We call $E=(\phi,F,H)$ an \emph{extension}.
We write $e_E=|H \sm H[F]|$, $v_E=|V(H)\sm F|$
and call $e_E$ the rank of $E$.
Now suppose $G$ is an $r$-multigraph on $[n]$.
We write $X_E(G)$ for the set or number of embeddings of $H$ 
in $G + \phi(H[F])$ that restrict to $\phi$ on $F$,
where sums of (multi)graphs are defined 
by viewing them as vectors over $\mb{N}$.
We say $E$ is \emph{$\oO$-dense} (in $G$) if $X_E(G) \ge \oO n^{v_E}$.
We say $G$ is \emph{$(\oO,h)$-extendable} if all 
extensions of rank $h$ are $\oO$-dense in $G$.
\end{defn}

Next we formulate our robust
fractional decomposition assumption
for an $r$-multigraph $G$,
saying that we can assign non-negative weights
to the $q$-cliques of $G$, maintaining an upper bound
on the ratio of any two weights, so that for any edge $e$
the total weight of cliques containing $e$ 
is roughly equal to the multiplicity of $e$ in $G$.

\begin{defn} \label{def:reg}
An $r$-multigraph $G$ on $[n]$ is 
\emph{$(K^r_q,c,\oO)$-regular} if there are
$w_{Q'} \in [\oO n^{r-q},\oO^{-1} n^{r-q}]$ for each $Q' \in K^r_q(G)$
with $\sum \{ w_{Q'} : e \in Q' \} = (1 \pm c) G_e$
for all $e \in \tbinom{[n]}{r}$.
\end{defn}

Note in particular that the upper bounds 
$w_{Q'} \le \oO^{-1} n^{r-q}$ in Definition \ref{def:reg}
imply $G_e \le \oO^{-1}$ for all $e \in K^r_n$.
We also reformulate our divisibility assumption
so that it applies to $r$-multigraphs $J \in \mb{N}^{K^r_n}$,
and more generally any $J \in \mb{Z}^{K^r_n}$.

\begin{defn} \label{Jdiv}
Suppose $J \in \mb{Z}^{K^r_n}$. 
We say that $J$ is \emph{$K^r_q$-divisible} 
if $\tbinom{q-i}{r-i}$ divides $\sum \{ J_e: f \sub e\}$ 
for any $0 \le i \le r$, $f \in \tbinom{[n]}{i}$. 
\end{defn}

Finally, we state our main theorem.

\begin{theo} \label{main+}
For any $q > r \ge 1$ there are $c_0>0$ and $n_0 \in \mb{N}$ 
such that if $h=2^{50q^3}$, $b = 2^{3^{r+q}}$,
$n>n_0$, $n^{-b^{-1}h^{-2}}<\oO<1$ and $c< c_0 \oO^h$,
then any $K^r_q$-divisible $(K^r_q,c,\oO)$-regular 
$(\oO,h)$-extendable $r$-multigraph on $n$ vertices
has a $K^r_q$-decomposition.
\end{theo}

\subsection{Previous work}

As a weaker version of the Existence Conjecture, 
Erd\H{o}s and Hanani \cite{EH} asked for approximate Steiner systems;
equivalently, finding $(1-o(1)) \tbinom{q}{r}^{-1} \tbinom{n}{r}$
edge-disjoint $K^r_q$'s in $K^r_n$. This was solved by 
R\"odl \cite{R}, who introduced a semi-random construction 
method known as the `nibble', which has since had a great impact 
on Combinatorics (see e.g.\ \cite{AKS,FR,Gr,Ka,Kim,KR,Kuz,PS,S,Vu} 
for related results and improved bounds).
It will also play an important role in this paper.
More recently, Ferber, Hod, Krivelevich and Sudakov \cite{FHKS}
gave a short probabilistic construction
in which every $r$-subset is covered by either one or two $q$-subsets.

Regarding exact results, we have already mentioned Wilson's theorem,
and Gustavsson's minimum degree generalisation thereof.
We should also note the seminal work of Hanani \cite{H1,H2},
which answers Steiner's problem 
for $(q,r) \in \{(4,2), (4,3), (5,2)\}$ and all $n$
(the case $(q,r)=(3,2)$ was solved by Kirkman, 
before Steiner posed the problem).
Besides these, we again refer to \cite{CD} as an introduction to
the huge literature on the construction of designs.
One should note that before the results of the current paper, 
there were only finitely many known Steiner systems with $r \ge 4$,
and it was not known if there were any Steiner systems with $r \ge 6$.

Even the existence of designs with $r \ge 7$ and any 
`non-trivial' $\lL$ was open before the breakthrough result 
of Teirlinck \cite{T} confirming this.
An improved bound on $\lL$ and a probabilistic method 
(a local limit theorem for certain random walks in high dimensions)
for constructing many other rigid combinatorial structures 
was given by Kuperberg, Lovett and Peled \cite{KLP}. 
Their result for designs is somewhat complementary to ours, 
in that they can allow the parameters $q$ and $r$ to grow with $n$,
whereas we require them to be (essentially) constant.
They also obtained much more precise estimates than we do
for the number of designs (within their range of parameters).

A different relaxation of the conjecture, 
which will play an important role in this paper,
is obtained by considering `integral designs', 
in which one assigns integers to the copies 
of $K^r_q$ in $K^r_n$ such that for every edge $e$ 
the sum of the integers assigned to the copies of $K^r_q$ 
containing $e$ is a constant independent of $e$.
Graver and Jurkat \cite{GJ} and Wilson \cite{W4} 
showed that the divisibility conditions 
suffice for the existence of integral designs 
(this is used in \cite{W4} to show the existence 
for large $\lL$ of integral designs with non-negative coefficients).
Wilson \cite{W5} also characterised the existence 
of integral $H$-decompositions for any $r$-graph $H$. 

\subsection{Subsequent work} \label{sec:subseq}

In the decade following the first arXiv version \cite{K} of this paper 
 there has been an explosion of solutions 
to long-standing open problems in Design Theory, 
some by developments of the ideas in this paper and some by 
developing new methods of `absorption' (see Section \ref{sec:strategy}).
For this `decennial edition' of the paper
we will just mention some highlights (a separate survey article 
would be needed to do justice to this task).

\begin{enumerate}
\item Besides various applications to Coding Theory, 
we mention the following applications in which not only our result 
but also the method of proof play an important role:
a conjectural analogue of the `expander mixing lemma'
for `high-dimensional permutations' proposed by Linial and Luria \cite{LL2};
the construction of $d$-dimensional coboundary expanders 
with bounded $(d-1)$-degrees by Lubotzky, Luria and Rosenthal \cite{LLR};
many w.h.p.~properties of random Latin Squares and Steiner Triple Systems,
starting with the existence of perfect matchings / transversals  by  Kwan \cite{Kw}
and including the solution of the McKay-Wanless conjecture
by Kwan, Sah and Sawhney \cite{KwSS} 
(see the latter paper for many further references).
\item Our proof strategy (see below) of applying a random construction and then absorption
leads to a fairly accurate estimate for the number of designs:
we gave the additional arguments needed for this in \cite{K2}.
Approximate enumeration of many other design-like structures,
such as high-dimensional permutations and sudoku squares,
follows from our more general results in \cite{K3}
(see also the expository articles \cite{K4,K5}).
\item Our counting method has also been combined 
with new absorption techniques for algebraically defined hypergraphs
to resolve the classical problem of approximately counting queens configurations,
in work with Bowtell \cite{BK}, independently by Simkin and Luria \cite{SL}, Simkin \cite{Sim}.
There are many design theoretic questions for algebraically defined hypergraphs,
so this seems to be a topic ripe for further developments. Another recent example 
with Sah and Sawhney \cite{KeSS} is the existence of subspace designs 
(which were conjectured not to exist!)
\item Most known counting results for design-like structures 
do not achieve the accuracy of an asymptotic formula --
this remains an open problem (even for Steiner Triple Systems).
In some cases, asymptotic formulae have been established by Fourier techniques,
notably for designs with large multiplicity $\lL$ and various other
 high multiplicity structures by Kuperberg, Lovett and Peled \cite{KLP}
and for semi-queens by Eberhard, Manners and Mrazovi\'c \cite{EMM}.
\item A second proof of the existence of designs, as well as some generalisations, 
such as the existence of $H$-decompositions for any hypergraph $H$
(a question from \cite{K}) was given by  Glock, K\"uhn, Lo and Osthus \cite{GKLO},
via the method of Iterative Absorption (discussed in more detail below), which has been a powerful tool
for many other problems, many discussed in the surveys by K\"uhn and Osthus \cite{KO} (on hamiltonicity)
and Kang,  Kelly, K\"uhn, Methuku and Osthus \cite{KKKMOsurvey} 
(on colouring, including an exposition of their proof  \cite{KKKMO2}  of the Erd\H{o}s-Faber-Lov\'asz conjecture).
\item In \cite{K3} we generalised the existence of designs to the setting of subset sums 
in lattices with coordinates indexed by labelled faces of simplicial complexes. 
As discussed in \cite{K5}, this general framework captures coloured and directed designs,
which lead to many further unexpected applications via combinatorial encoding techniques.
This approach was used in the solution of the Oberwolfach problem by Glock, Joos, Kim, K\"uhn and Osthus, 
and the generalised Oberwolfach problem in work with Staden \cite{KeSt}.
Staden and I \cite{KeSt2} also used it for Ringel's tree packing conjecture,
which was solved independently by Montgomery, Pokrovskiy and Sudakov \cite{MPS},
using a different set of ideas arising from recent progress on rainbow embedding problems
(see Pokrovskiy's survey \cite{Psurvey}). It was also used for progress on the (still open)
Gyarf\'as tree packing conjecture by Allen, B\"ottcher, Clemens, Hladk\'y, Piguet and Taraz \cite{ABCHPT}.
\item A conjecture of Erd\H{o}s \cite{Erd} on the existence of Steiner Triple Systems of high girth was solved by
Kwan, Sah, Sawhney and Simkin \cite{KwSSS}, building on the approximate solutions obtained
independently by Bohman and Warnke \cite{BW} and by Glock, K\"uhn, Lo, and Osthus \cite{GKLOsts}.
The more general High Girth Existence Conjecture posed in \cite{GKLOsts} was recently solved 
by Delcourt and Postle \cite{DPhi} using their new method (also giving a third proof of the existence of designs)
 of Refined Absorption \cite{DPref}, following the approximate solution via conflict-free matchings
independently by Delcourt and Postle \cite{DPalmost} and by
Glock, Joos, Kim, K\"uhn and Lichev \cite{GJKKL}.
\item The scope of the absorption method in sparse settings has been greatly extended by Montgomery,
via his method of Distributed Absorption introduced in \cite{M1} 
to resolve Kahn's Conjecture on the appearance threshold 
for bounded-degree spanning trees in random graphs. 
Among many other applications of Distributed Absorption (which needs its own survey),
we mention the solutions by M\"uyesser and Pokrovskiy \cite{MuPo}
to several old problems in Combinatorial Group Theory,
and Montgomery's spectacular work \cite{Mrbs} on the Ryser-Brualdi-Stein Conjecture
(following our recent progress with Pokrovskiy, Sudakov and Yepremyan \cite{KPSY}).
\item The theory of thresholds in random structures has been revolutionised by the solution
of the Kahn-Kalai Conjecture \cite{KK} by Park and Pham \cite{PP},
following the fractional solution (Talagrand's Conjecture \cite{Tal}) 
by Frankston, Kahn, Narayanan and Park \cite{FKNP}.
Despite this progress, it still seems challenging to determine
thresholds for designs or design-like structures.
One recent success was determining the threshold 
for Latin Squares / Steiner Triple Systems,
independently by me \cite{K6} and by Jain and Pham \cite{JP},
following recent progress by Sah, Sawhney and Simkin \cite{SSS} 
and by Kang, Kelly, K\"uhn, Methuku and Osthus \cite{KKKMO}.
The general problem remains open,
with some progress by Delcourt, Kelly and Postle \cite{DKP} via Refined Absorption.
\end{enumerate}

\subsection{Proof strategy} \label{sec:strategy}

Our main new idea is to use a Randomised Algebraic Construction:
the first step of our construction is to take a random subset 
of an algebraically defined `model' for designs.
This results in a partial decomposition 
that covers a constant fraction of the edge set,
and also carries a rich structure of possible local modifications.
We treat this partial decomposition as a template for the final decomposition.
By various applications of the nibble and greedy algorithms, 
we can choose another partial decomposition that covers 
all edges not in the template, which also spills over 
slightly into the template, so that every edge is covered
once or twice, and very few edges are covered twice
(we call the latter the `spill').
The crucial point is that the choice of the
template was such that the spill can be `absorbed',
converting the approximate decomposition 
into a (perfect) decomposition.

At this level of generality, our method sounds somewhat similar to the 
Absorbing Method of R\"odl, Ruci\'nski and Szemer\'edi \cite{RRS}
(see also the survey \cite{RR}). However, in the Absorbing Method
(in its basic form) as applied to the problem of designs,
the analogue of our template would be a random sparse partial 
decomposition (without any superimposed algebraic structure),
and it is not hard to see that local modifications
have a negligible probability of appearing in such a construction.

Another way to think about the failure of the naive Absorbing Method
is that there are too many possibilities for the `leave'
of the approximate decomposition. This viewpoint suggests the
more sophisticated approach of Iterative Absorption used
in \cite{GKLO}, in which the leave becomes gradually
more constrained, until there are so few options that 
each possible leave can have its own private `absorber'.
Moreover, it turns out that via an iterative approach
one can in fact handle many possibilities
for the leave in a single absorber:
this is the key idea behind Refined Absorption \cite{DPref}.

By contrast, our construction blends randomness with algebra,
in such a way that any approximate decomposition can be absorbed.
The rich rigid structures of Algebra make it a natural tool in
the construction of designs. For example,
it is not hard to see that orbits of $r$-transitive 
permutation groups are $r$-designs,
but apparently there are no $r$-transitive groups with $r>5$ 
other than the symmetric and alternating groups\footnote{
This claim is folklore among Group Theorists, assuming
the Classification of Finite Simple Groups (see \cite{GLS}).},
which points to the limitations of the purely algebraic approach.

Nevertheless, we will see that a suitable algebraically defined 
template has a dense well-distributed set of cliques 
that are `absorbable', in that they can be included
in the clique decomposition of the template
via a suitable local modification. Our template can be thought
of as a general absorber, which is almost as effective as that in  \cite{DPref},
except that it requires a somewhat stronger property of the edges to be absorbed
(linear boundedness, as opposed to just boundedness).

To make use of this structure, we first find an 
`integral decomposition' of the spill (mentioned above),
which can be thought of as a decomposition 
in which we can take each clique with any integer weight;
this is the point in the proof where 
the divisibility assumption is used.

Next we apply a `clique exchange algorithm'
that replaces the integral decomposition 
by a `signed decomposition', which can be
thought of as two partial decompositions,
called `positive' and `negative',
such that the underlying hypergraph 
of the negative decomposition is contained 
in that of the positive decomposition,
and the difference forms a `hole'
that is precisely equal to the spill.

We further ensure that each positive clique
can be absorbed into the template,
via a series of absorptions that we call a `cascade'
(this is the most technically challenging part of the paper).

Finally, deleting the positive cliques and replacing them by the
negative cliques eliminates one of the two uses of each edge in
the spill, so that we end up with a perfect decomposition.

The above remarks hopefully give the flavour of the proof strategy;
we will defer more detailed proof sketches of the various steps
to the relevant later sections of the paper, and also give a summary
of the proof in Section \ref{sec:conclude} for convenient reference.
In the next subsection
we make some further remarks on the implementation of the strategy
and how it differs from that in the first version \cite{K} of this paper.

\subsection{Implementation} \label{sub:imp}

While the overall proof strategy in this version
of the paper is the same as in the first version \cite{K},
the details of the implementation (introduced in 2018 for the second version)
are substantially different and considerably simpler.
The most important difference is that we now do not need
any inductive argument for reducing the vertex set.
There was an error in this part of our argument in \cite{K},
which was kindly pointed out by the authors of \cite{GKLO},
namely in the proof of \cite[Lemma 6.3]{K}. The lemma is true,
and the proof can be fixed with more sophisticated random greedy 
arguments, but this would make \cite{K} even more complicated,
whereas the issue is entirely avoided by our new approach.
Furthermore, we can
work entirely in the simpler setting of uniform hypergraphs,
rather than the more general setting of simplicial complexes
that was needed in \cite{K} for the purposes of induction.%
\footnote{The argument here does apply to the simplicial
complex setting, and so can be applied to the
results from \cite{K} that used simplicial complexes, namely
Theorems 6.6 and 6.7, but we omit this
for simplicity of exposition.}
The simpler method presented here may also be more amenable
to computer implementation with a view to constructing explicit designs.

\begin{SCfigure}
\caption{Two pictures of the associated octahedron for $xyz$. 
Red triangles belong to the template.
These may be swapped with the white triangles,
thus using $xyz$ while covering the same set of edges.}
\includegraphics[width=0.3\textwidth]{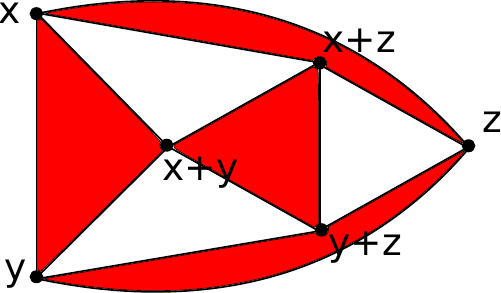}

\begin{tikzpicture}[line join=bevel,z=-5.5]
\coordinate (A1) at (0,0,-1);
\coordinate (A2) at (-1,0,0);
\coordinate (A3) at (0,0,1);
\coordinate (A4) at (1,0,0);
\coordinate (B1) at (0,1,0);
\coordinate (C1) at (0,-1,0);
\draw  (A1) -- (A2) -- (B1) -- cycle;
\draw [fill opacity=0.7,fill=red!70!black] (A4) -- (A1) -- (B1) -- cycle;
\draw [fill opacity=0.7,fill=red!70!black] (A1) -- (A2) -- (C1) -- cycle;
\draw (A4) -- (A1) -- (C1) -- cycle;
\draw [fill opacity=0.7,fill=red!70!black] (A2) -- (A3) -- (B1) -- cycle;
\draw (A3) -- (A4) -- (B1) -- cycle;
\draw  [fill opacity=0.7,fill=white]  (A2) -- (A3) -- (C1) -- cycle;
\draw [fill opacity=0.7,fill=red!70!black] (A3) -- (A4) -- (C1) -- cycle;
\end{tikzpicture}
\hspace{0.3cm}
\label{fig:octmove}
\end{SCfigure}

To develop some intuition for Randomised Algebraic Construction
it is helpful to first consider the special case of triangle
decompositions of typical graphs (see \cite{K2}).
Our algebraic model for triangle decompositions is the set 
of all triples $xyz$ with $x+y+z=0$ in some abelian group $\GG$.
Indeed, this is almost a triangle decomposition of the complete
graph on $\GG$, in that for any $xy$ there is a unique $z$
with $x+y+z=0$, but this ignores the possibility that $x,y,z$ 
may not be distinct, and also that our approach requires
decompositions of (hyper)graphs that are not complete.
Instead, to define the template of a graph $G$ in \cite{K2}, 
we randomly embed $V(G)$ in $\mb{F}_{2^a}$ for some $a$ 
such that $2^a$ is not much bigger than $|V(G)|$, 
and take all triangles $xyz$ satisfying $x+y+z=0$,
which gives a partial triangle decomposition of $G$.
In this construction, a triangle $xyz$ is absorbable
if $G$ contains the `associated octahedron' of $xyz$,
which is the complete $3$-partite graph with parts
$\{x,y+z\}$, $\{y,x+z\}$, $\{z,x+y\}$.
Indeed, this octahedron has two distinct triangle decompositions,
one of which contains $xyz$, and the other of which consists 
entirely of triangles with zero sum (see Figure~\ref{fig:octmove}).

In general, for our algebraic model of a $K^r_q$-decomposition,
we consider a vertex set that is a finite field, and
a set of $q$-cliques that correspond to the image of some
$q \times r$ matrix $M$ that is `generic'
(every square submatrix of $M$ is nonsingular).
The motivation for this model is that 
for every $r$-set $e$ of field elements
and injective map $\pi_e: e \to [q]$, 
we can reconstruct the unique vector $y$ such that 
$(My)_i = x$ for all $x \in e$, $i = \pi_e(x)$.
However, if we embed some $r$-graph $G$ in the field
and use this construction, 
then as in the triangle decomposition case, 
we can only use the subset of the model
that uses edges which are actually present in the given $r$-graph $G$.
Furthermore, as we must use each edge at most once, we make each edge $e$ 
randomly `decide' on some fixed injection $\pi_e$, and we only 
allow $q$-cliques that are compatible with these choices.

Similarly to the case of triangle decompositions,
we randomly embed $V(G)$ in $\mb{F}_{p^a}$,
for some prime $p$ which is large compared with $q$
but small compared with $n=|V(G)|$, and some $a$ 
such that $p^a$ is not much bigger than $n$.
Viewing $\mb{F}_{p^a}$ as a vector space over $\mb{F}_p$
we find a rich set of absorbable cliques via a construction 
somewhat analogous to the associated octahedra of triangles 
(this part of the argument is new to this version
and is much simpler than the approach used in the first version). 
In fact, rather than using just using one embedding 
of $V(G)$ in $\mb{F}_{p^a}$, we use $z$ such embeddings,
for some $z$ which is large compared with $q$
but small compared with $n$. The point is that
with positive probability every $r$-set in $V(G)$
has full dimension\footnote{
Regarding $\mb{F}_{p^a}$ as a vector space over $\mb{F}_p$,
we say that $e \sub V(G)$ has full dimension
in $f:V(G) \to \mb{F}_{p^a}$ 
if $f(e)$ spans a subspace of dimension $|e|$.}
in most of these embeddings,
which circumvents many technical difficulties from the first version
regarding the treatment of degenerate sets.

To illustrate the previous discussion 
in the case of triangle decompositions
we could take $M^t = \tbinom{1\ 1\ 1}{1\ 2\ 3}$.
Then to include a triangle $xyz$ of $G$ in the template
we require `compatibility' 
with respect to the random injections
$\pi_{xy}$, $\pi_{xz}$, $\pi_{xz}$,
say $\pi_{xy}(x)=\pi_{xz}(x)=1$,
$\pi_{xy}(y)=\pi_{yz}(y)=2$,
$\pi_{xz}(z)=\pi_{yz}(z)=3$,
and with respect to one of 
the random injections $f_j:V(G) \to \mb{F}_{p^a}$,
meaning that for some $u,v \in \mb{F}_{p^a}$ we have
$f_j(x)=u+v$, $f_j(y)=u+2v$, $f_j(z)=u+3v$.
(We also require certain `activation' events that account
for weights on triangles as in Definition \ref{def:reg}.)   
The lack of symmetry (resolved by the injections $\pi_e$)
significantly complicates the task of finding switching operations
analogous to the octahedral construction illustrated above;
the construction used for the general case
will be illustrated in Figure~\ref{fig:ab} below.

A further comment on the new implementation
is that we have found a considerably simpler approach
for constructing `bounded integral designs'.
As described above, Graver and Jurkat \cite{GJ} 
and Wilson \cite{W4} showed that the divisibility conditions
suffice for the existence of integral designs,
but our modification approach requires 
an additional local boundedness property.
Our new approach for bounded integral designs
relies on `robust local decodability' 
of the lattice of $K^r_q$-divisible vectors:
there is some constant $N=N(q)$ such that
for any $e \in K^r_n$ there are `many' integral 
combinations of $q$-cliques that equal
the vector in $\mb{Z}^{K^r_n}$ with $N$ 
in coordinate $e$ and $0$ otherwise.

It is interesting that local decodability 
was a key property in the general framework of \cite{KLP},
although we do not see any connection between this part
of our proof and their approach.
Furthermore, there are many natural related problems in design theory
that do not exhibit local decodability, such as 
`generalised partite hypergraph decompositions',
which encompass problems such as resolvable hypergraph designs,
large sets of hypergraph designs, decompositions of designs by designs,
high-dimensional permutations and Sudoku squares (see \cite{K3}).
Here the method from the first version of this paper can be applied:
the key idea is to solve the fractional relaxation of the integral design problem
(we allow rational weights of either sign),
and use this in an iterative rounding algorithm to obtain
finer approximations to an exact solution until the approximation
is so good that a trivial argument can be used to complete the solution.
However, the general integral relaxation has a much more complicated structure,
so there are many further difficulties to overcome (see \cite{K3}). 

Writing in 2024, we also remark on another feature of our 2018 strategy that 
in retrospect has independent interest, namely bounded generation (see subsection \ref{sub:bddgen}).
Here we generalised the results of \cite{GJ,W4} 
on generating nullspaces of shadow operators by octahedra,
by showing that it suffices to use a `thin' set of octahedra,
meaning that every edge is in only constantly many octahedra of the generating set.
This is analogous to the Refined Absorption recently introduced by Delcourt and Postle \cite{DPref}
which has several other applications already mentioned in subsection \ref{sec:subseq}.
Based on these ideas, we have recently discovered a new proof \cite{Knew} of the existence of designs
that is much shorter than any of the previous proofs and also gives reasonable bounds;
the methods in the current paper remain interesting 
for their other applications mentioned in subsection \ref{sec:subseq}.

\subsection{Organisation}

The organisation of this paper is as follows.
The next section contains various preliminary results
used throughout the paper, on concentration of probability,
almost perfect matchings in hypergraphs, and extensions.
In Section \ref{sec:template} we construct the template,
and establish its combinatorial extendability properties.
Section \ref{sec:approx} contains the nibble and cover arguments
that complete the template to an approximate decomposition,
namely a set of cliques such that every edge
is covered once or twice, and the set of edges covered twice
(the `spill') forms a suitably bounded subgraph of the template.
In Section \ref{sec:int} we find a suitably bounded 
integral decomposition of the spill.
In Section \ref{sec:ab} we analyse the algebraic properties 
of the template, showing that it has a rich structure
of absorbable and cascading cliques that
can be used for local modifications.
Section \ref{sec:cea} analyses the Clique Exchange Algorithm
that modifies the integral decomposition
so that the spill can be absorbed into the template.
In the final section we 
complete the proof of our main theorem
and make some concluding remarks.

\subsection{Notation and terminology}

Here we gather some notation and terminology 
that is used throughout the paper.
We write $[n] = \{1,\dots,n\}$.
For a set $S$, we write $\tbinom{S}{r}$ 
for the set of $r$-subsets of $S$.
We write $Q = \tbinom{[q]}{r}$ and also $Q = \tbinom{q}{r}$
(it will be clear from the context whether 
we are referring to the set or its size).
We identify $Q = \tbinom{[q]}{r}$ with the
edge set of $K^r_q$ (the complete $r$-graph on $[q]$).

For any set $S$ we write $K^r_q(S)$ for the complete
$q$-partite $r$-graph with parts of size $|S|$ 
where each part is identified with $S$.
If $S=[s]$ we write $K^r_q(S)=K^r_q(s)$.

We often use `concatenation notation' for sets,
for example $xyz$ may denote $\{x,y,z\}$,
and for function composition,
for example $fg$ may denote $f \circ g$.

We say that an event $E$ holds with high probability (w.h.p.) 
if $\mb{P}(E) = 1-e^{-\Omega(n^c)}$ for some $c>0$ as $n \to \infty$.
Whenever we make any such statement, 
we are implicitly assuming that
$n>n_0(q)$ is sufficiently large.
Then by union bounds we can assume that any specified 
polynomial number of such events all occur. 

Suppose $X$ and $Y$ are sets. 
We write $Y^X$ for the set of vectors with entries in $Y$
and coordinates indexed by $X$, which we also identify
with the set of functions $f:X \to Y$.
For example, we may consider $v \in \mb{F}_p^q$ 
as an element of a vector space over $\mb{F}_p$
or as a function from $[q]$ to $\mb{F}_p$.

We identify $v \in \{0,1\}^X$ with the set $\{x \in X: v_x=1\}$.
We identify $v \in \mb{N}^X$ with the multiset in $X$ 
in which $x$ has multiplicity $v_x$
(for our purposes $0 \in \mb{N}$). 
We also apply similar notation and terminology as for multisets 
to vectors $v \in \mb{Z}^X$. 
We often consider algorithms with input $v \in \mb{Z}^X$,
where each $x \in X$ is considered $|v_x|$ times,
with a sign attached to it (the same as that of $v_x$);
then we refer to $x$ as a `signed element' of $v$.

Arithmetic on vectors in $\mb{Z}^X$ is to be understood pointwise:
 $(v+v')_x = v_x + v'_x$ and $(vv')_x = v_x v'_x$ for $x \in X$.
For $v \in \mb{Z}^X$ we write $|v| = \sum_{x \in X} |v_x|$.
We also write $v=v^+-v^-$, where $v^+_x = \max\{v_x,0\}$ 
and $v^-_x = \max\{-v_x,0\}$ for $x \in X$.
For $X' \sub X$ we define $v[X'] \in \mb{Z}^{X'}$ 
by $v[X']_x = v_x$ for $x \in X'$.
We often identify any $v' \in \mb{Z}^{X'}$ 
with an element of $\mb{Z}^X$ 
by letting $v'_x=0$ for $x \in X \sm X'$.

If $G$ is a hypergraph, $v \in \mb{Z}^G$ and $e \sub V(G)$
we define $v(e) \in \mb{Z}^{G(e)}$ 
by $v(e)_f = v_{e \cup f}$ for $f \in G(e)$.

We say $J \in \mb{Z}^{K^r_n}$ is $\tT$-bounded 
if $\sum \{ |J_e|: f \sub e \in K^r_n \} < \tT n$ 
for all $f \in \tbinom{[n]}{r-1}$. 

We denote the standard basis vectors in $\mb{R}^d$
by $e_1,\dots,e_d$. Given $I \sub [d]$, 
we let $e_I$ denote the matrix 
with rows indexed by $I$ and columns indexed by $[d]$,
in which the row indexed by $i \in I$ is $e_i$.

We write $M \in \mb{F}_p^{q \times r}$ to mean that $M$ 
is a matrix with $q$ rows and $r$ columns
having entries in $\mb{F}_p$.
For $I \in Q = \tbinom{[q]}{r}$ we let $M_I$ be the 
square submatrix with rows indexed by $I$.
Note that $M_I = e_I M$.

We will regard $\mb{F}_{p^a}$ as a vector space 
over $\mb{F}_p$. For $e \sub \mb{F}_{p^a}$
we write $\dim(e)$ for the dimension of
the subspace spanned by the elements of $e$.
For $e \in \mb{F}_{p^a}^d$ we write $\dim(e)$
for the dimension of the set of coordinates of $e$.

When we use `big-O' notation, 
the implicit constant will depend only on $q$.

We write $a = b \pm c$ to mean $b - c \leq a \leq b + c$.

Throughout the paper we omit floor and ceiling symbols 
where they do not affect the argument. 

For convenient reference, we list here
several parameters used throughout the paper:
\begin{gather*}
Q = \tbinom{q}{r}, \quad z=h=2^{50q^3}, 
\quad b = 2^{3^{r+q}}, \quad
n^{-b^{-1}h^{-2}}<\oO<\oO_0(q), \quad 
p \text{ is a prime with } 2^{8q}<p<2^{9q}, \\
a \in \mb{N} \text{ with } p^{a-2} < n \le p^{a-1}, \quad
\gG = np^{-a}, \quad
\rho = \oO z^{-Q} q! (q)_r^{-Q} \gG^{q-r}, 
\text{ where } (q)_r=q!/(q-r)!,\\
N=(2q)^q !, \quad
c=\oO^h, \quad c_1 = (2Qc)^{1/2Q}, \quad
c_{i+1} = \oO^{-h/20Q} c_i \ \text{ for } \ i \in [4].
\end{gather*}
The multiplicative factor of $\oO^{-h/20Q}$ between
successive $c_i$'s is chosen so that there is plenty
of room to spare in the various inequalities below,
so we will omit detailed discussion of these
during the proof. We remark here that the tightest
inequality occurs during the cascade algorithm
in the proof of Theorem \ref{main+}, namely
$2^r p^{2q} r!\oO^{-p^{q^2}} c_4 < c_5/12$,
which holds easily as $p^{q^2} < 2^{9q^3} < h^{1/5}$ 
and $Q < 2^q < 2^{q^2} = h^{1/50q}$.
The assumption $\oO > n^{-b^{-1}h^{-2}}$ 
is much stronger than needed for the proof,
but we are only interested in establishing
some polynomial dependence, as in any case
the best bounds available from our proof are
presumably far from optimal.

\section{Preliminaries}

In this section we gather some results 
that will be used throughout the paper,
concerning concentration of probability,
almost perfect matchings in hypergraphs,
and extensions.

\subsection{Concentration of probability}

We make the following standard definitions.%
\footnote{In this paper all probability spaces are finite,
and will only be referred to implicitly via random variables.
We will only ever consider the natural filtration
$\mc{F}=(\mc{F}_i)_{i \ge 0}$ associated with a random process, 
where each $\mc{F}_i$ consists of all events 
determined by the history of the process up to step $i$.} 

\begin{defn}
Let $\OO$ be a (finite) probability space.
An \emph{algebra} (on $\OO$) is a set $\mc{F}$ of subsets of $\OO$ 
that includes $\OO$ and is closed under intersections 
and taking complements.
A \emph{filtration} (on $\OO$) is a sequence 
$\mc{F}=(\mc{F}_i)_{i \ge 0}$ 
of algebras such that $\mc{F}_i \sub \mc{F}_{i+1}$ for $i \ge 0$.
A sequence $A = (A_i)_{i \ge 0}$ of random variables on $\OO$
is a \emph{supermartingale} (w.r.t.\ $\mc{F}$)
if each $A_i$ is {\em $\mc{F}_i$-measurable}
(all $\{\oO: A_i(\oO)<t\} \in \mc{F}_i$) 
and $\mb{E}(A_{i+1}|\mc{F}_i) \le A_i$ for $i \ge 0$.
\end{defn}

Now we can state a general result of Freedman 
\cite[Proposition 2.1]{F} 
that essentially implies all of the bounds we will use 
(perhaps with slightly weaker constants).

\begin{lemma} \label{freed} 
Let $(A_i)_{i \ge 0}$ be a supermartingale 
w.r.t.\ a filtration $\mc{F}=(\mc{F}_i)_{i \ge 0}$.
Suppose that $A_{i+1}-A_i \le b$ for all $i \ge 0$, 
and let $E$ be the `bad' event that there 
exists $j \ge 0$ with $A_j \ge A_0 + a$ 
and $\sum_{i=1}^j Var[A_i \mid \mc{F}_{i-1}] \le v$.
Then $\mb{P}(E) \le \exp \brac{-\tfrac{a^2}{2(v+ab)}}$.
\end{lemma}

We proceed to give some useful consequences of Lemma \ref{freed}.
First we make another definition.

\begin{defn} \label{def:dom}
Suppose $Y$ is a random variable and 
$\mc{F} = (\mc{F}_0,\dots,\mc{F}_n)$ is a filtration.
We say that $Y$ is \emph{$(C,\mu)$-dominated (w.r.t.\ $\mc{F}$)} 
if we can write $Y = \sum_{i=1}^n Y_i$, 
where $Y_i$ is $\mc{F}_i$-measurable, $|Y_i| \le C$
and $\mb{E}[|Y_i| \mid \mc{F}_{i-1}] < \mu_i$ for $i \in [n]$, 
where $\sum_{i=1}^n \mu_i < \mu$.
\end{defn}

\begin{lemma} \label{dom}
If $Y$ is $(C,\mu)$-dominated then 
$\mb{P}(|Y|>(1+c)\mu) < 2e^{-\mu c^2/2(1+2c)C}$.
\end{lemma}

\nib{Proof.}
Let $A_i = \sum_{j<i} (Y_j-\mu_j)$ for $i \ge 0$; 
then $(A_i)_{i \ge 0}$ is a supermartingale and
\[ Var[A_i \mid \mc{F}_{i-1}] = Var[Y_i \mid \mc{F}_{i-1}] 
 \le \mb{E}[Y_i^2 \mid \mc{F}_{i-1}] 
\le C \mb{E}[|Y_i| \mid \mc{F}_{i-1}] \le C\mu_i.\]
By Lemma \ref{freed} applied with $a=c\mu$, $b=2C$ and $v=C\mu$
we obtain \[\mb{P}(Y>(1+c)\mu) < e^{-\mu c^2/2(1+2c)C}.\]
Similarly, considering $A_i = - \sum_{j<i} (Y_j+\mu_j)$
gives the same estimate for $\mb{P}(Y<-(1+c)\mu)$. \qed

\begin{rem} \label{rem:dom}
All of our applications of Lemma \ref{dom} 
will be such that we could also deduce concentration 
by coupling to a sum of bounded independent variables
and applying `Bernstein's inequality'
or the `Chernoff bound' for binomial variables
(see e.g.\ \cite[Remark 2.9]{JLR}).
In many cases, we will actually have 
a sum of bounded independent variables,
which for brevity we call `pseudobinomial',
where these inequalities can be directedly
applied with no need for a coupling.
\end{rem}

Next we record several consequences of 
the well-known inequality of Azuma \cite{Az}
(see e.g.\ \cite{McD}).

\begin{defn} \label{def:lip1}
Suppose $f:S \to \mb{R}$ where $S = \prod_{i=1}^n S_i$
and $b = (b_1,\dots,b_n)$ with $b_i \ge 0$ for $i \in [n]$.
We say that $f$ is \emph{$b$-Lipschitz} if for any 
$s,s' \in S$ that differ only in the $i$th coordinate
we have $|f(s)-f(s')| \le b_i$. 
We also say that $f$ is \emph{$B$-varying} where $B=\sum_{i=1}^n b_i^2$.
\end{defn}

\begin{lemma} \label{lip1}
Suppose $Z = (Z_1,\dots,Z_n)$ is a sequence 
of independent random variables,
and $X=f(Z)$, where $f$ is a $B$-varying function.
Then $\mb{P}(|X-\mb{E}X|>a) \le 2e^{-a^2/2B}$.
\end{lemma}

\begin{defn} \label{def:lip2}
Let $S_n$ be the symmetric group on $[n]$.
Suppose $f:S_n \to \mb{R}$ and $b \ge 0$.
We say that $f$ is \emph{$b$-Lipschitz} if 
whenever $\sS = \tau \circ \sS'$ 
for some transposition $\tau \in S_n$
we have $|f(\sS)-f(\sS')| \le b$. 
\end{defn}

\begin{lemma} \label{lip2} 
Suppose $f:S_n \to \mb{R}$ is $b$-Lipschitz,
$\sS \in S_n$ is uniformly random and $X=f(\sS)$.
Then $\mb{P}(|X-\mb{E}X|>a) \le 2e^{-a^2/2nb^2}$.
\end{lemma}

We will use a common generalisation 
of Lemmas \ref{lip1} and \ref{lip2},
which perhaps has not appeared before, but is proved in the same way.
It considers functions in which the input consists of $n$ independent 
random injections $\pi_i:[a'_i] \to [a_i]$:
if $a'_i=1$ this is a random element of $[a_i]$;
if $a'_i=a_i$ this is a random permutation of $[a_i]$.

\begin{defn} \label{def:genlip}
Let $a = (a_1,\dots,a_n)$ and $a' = (a'_1,\dots,a'_n)$, 
where $a_i \in \mb{N}$ and $a'_i \in [a_i]$ for $i \in [n]$,
and $\Pi(a,a')$ be the set of $\pi=(\pi_1,\dots,\pi_n)$
where $\pi_i:[a'_i] \to [a_i]$ is injective. 
Suppose $f:\Pi(a,a') \to \mb{R}$ and $b = (b_1,\dots,b_n)$ 
with $b_i \ge 0$ for $i \in [n]$.
We say that $f$ is \emph{$b$-Lipschitz} if for any $i \in [n]$
and $\pi,\pi' \in \Pi(a,a')$ such that $\pi_j = \pi'_j$ for $j \ne i$
and $\pi_i = \tau \circ \pi'_i$ for some transposition $\tau \in S_{a_i}$
we have $|f(s)-f(s')| \le b_i$. 
We also say that $f$ is \emph{$B$-varying} 
where $B=\sum_{i=1}^n a'_i b_i^2$.
\end{defn}

\begin{lemma} \label{lip3} 
Suppose $f:\Pi(a,a') \to \mb{R}$ is $B$-varying,
$\pi \in \Pi(a,a')$ is uniformly random and $X=f(\pi)$.
Then $\mb{P}(|X-\mb{E}X|>t) \le 2e^{-t^2/2B}$.
\end{lemma}

\begin{rem}
We require a slight generalisation
of Lemma \ref{lip3}, which is the same statement
under any distribution on $\pi \in \Pi(a,a')$
such that the $\pi_i$ are independent and uniform,
except for $i$ such that $a'_i=1$,
for which we allow any distribution
(thus generalising Lemma \ref{lip1}).
This follows from Azuma's inequality.
\end{rem}

\subsection{Almost perfect matchings}

The following theorem of Pippenger (unpublished, generalised in \cite{PS}) 
generalises the result of R\"odl mentioned in the introduction:
it gives a nearly perfect matching in any uniform hypergraph 
that is approximately regular and has small codegrees.

\begin{theo} \label{pip}
For any integer $k \ge 2$ and real $a>0$ there is $b > 0$
such that for any $D,n \in \mb{N}$, if $A$ is a $k$-graph on $[n]$
with $|A(x)| = (1\pm b) D$ for every vertex $x$
and $|A(xy)| < bD$ for every pair of vertices $x,y$,
then $A$ has a matching covering all but at most $an$ vertices. 
\end{theo}

For our purposes, 
$A$ will be a $Q$-graph, where $Q=\tbinom{q}{r}$,
where $V(A)$ is the set of edges (with multiplicity)
in some $r$-multigraph $G$ on $n$ vertices, and
$E(A) = \{ \tbinom{Q'}{r}: Q' \in \mc{Q} \}$
for some set $\mc{Q}$ of $q$-cliques.
The vertex degree assumption on $A$ translates
into saying that every edge of $G$ is in
roughly the same number of cliques in $\mc{Q}$.
The codegree assumption on $A$ will hold 
with plenty of room to spare, just using
the trivial bound that any pair of distinct $r$-sets
are contained in at most $n^{q-r-1}$ cliques.
The conclusion of Theorem \ref{pip} is that
we obtain a set of edge-disjoint cliques
covering almost all edges of $G$.
In fact, we will require the following stronger 
boundedness property of the `leave' 
(the submultigraph formed by the uncovered edges). 

\begin{defn} \label{bdd} 
Suppose $J$ is an $r$-multigraph on $[n]$ and $\tT>0$. 
We say that $J$ is \emph{$\tT$-bounded} if
$|J(e)| < \tT n$ for all $e \in \tbinom{[n]}{r-1}$. 
\end{defn}

Now we will add the required boundedness
property of the leave to the conclusion of Theorem \ref{pip},
and also quantify (to some extent) the dependency
of the size of the leave on the regularity of $A$.
There has been considerable effort in the literature
(see \cite{AKS, Gr, Kim, Vu}) regarding the latter point,
but for our purposes we only care that there is 
some polynomial dependence, as other arguments in our
paper only operate up to this level of accuracy.
The proof of the following lemma is an easy modification
of that given in \cite{Gr}, so we omit it.
(Addendum in 2024: these ideas have now been organised
into a generally applicable tool by Ehard, Glock and Joos \cite{EGJ}.)

\begin{lemma} \label{nibble+}
There are $b_0>0$ and $n_0 \in \mb{N}$ 
such that for $n>n_0$ and $n^{-1/Q} < b < b_0$, 
given any $r$-multigraph $G$ on $[n]$ vertices 
and a set\footnote{Note that we say `set',
not `multiset', so the auxiliary hypergraph
has codegrees $O(n^{q-r-1})$.} 
of $q$-cliques $\mc{Q}$ such that
every $r$-set $e$ is in $(1 \pm b) dn^{q-r} G_e$ 
elements of $\mc{Q}$, where $d>n^{-1/Q}$, there is 
$M^n \sub \mc{Q}$ such that\footnote{Here $\sum M^n$ denotes
the $r$-multigraph obtained by summing the cliques in $M^n$
and $L \ge 0$ (pointwise) means that its edge multiplicities 
are at most those in $G$, so if we regard copies 
of any given edge as distinct elements we can think 
of $M^n$ as a set of edge-disjoint cliques in $G$.} 
$L = G - \sum M^n \ge 0$ is $b^{1/2Q}$-bounded.
\end{lemma}

\subsection{Extensions} \label{sec:ext}

We conclude our preliminary section
with some basic properties of extensions
(see Definition \ref{def:ext})
that will be used throughout the paper.
First we make some comments on the definition.
It is important to note that edges of $H$ contained 
within $F$ have no effect on $X_E(G)$.
There is no loss of generality in assuming 
that $H$ has no isolated vertices,
which has the convenient consequence $v_E \le re_E$.
In the case $F=\es$ extendability gives 
a lower bound on the number of embeddings of $H$ in $G$.
In particular, if $H$ consists of a single edge
then we obtain the density bound $d(G) \ge \oO$.
We also note that
if $|V(H) \sm F|=1$ then $X_E(G)$ is an
intersection of neighbourhoods of the type
that appears in Definition \ref{def:typ}.
This explains the following result,
which gives an estimate
for the number of extensions in typical
$r$-graphs that is close to what would be
expected in a random $r$-graph of the same density.

\begin{lemma} \label{ext}
Let $G$ be a $(c,h)$-typical $r$-graph on $[n]$, where $c < h^{-2}$.
Suppose $E=(\phi,F,H)$ is an extension with $|H| \le h$. 
Then $X_E(G) = (1 \pm (e_E+1)c) d(G)^{e_E} n^{v_E}$.
\end{lemma}

\nib{Proof.}
Write $V(H) \sm F = \{v_1,\dots,v_{v_E}\}$
and suppose for $i \in [v_E]$ that there
are $e_i$ edges of $H$ using $v_i$
but not using any $v_j$ with $j>i$.
We can construct any embedding in $X_E(G)$
by choosing the images of the $v_i$'s
successively. By Definition \ref{def:typ},
the number of choices for $v_i$ 
given any previous choices is 
$(1 \pm e_i c) d(G)^{e_i} n$.
The lemma follows by multiplying these estimates,
using $e_E = \sum e_i$ and 
$(1+c)^h \le e^{hc} \le 1+hc+(hc)^2 \le 1+(h+1)c$. \qed

\medskip

We can now show that Theorem \ref{main} 
on $K^r_q$-decompositions of dense typical $r$-graphs
follows from Theorem \ref{main+}
on $K^r_q$-decompositions of dense $r$-graphs
satisfying the extendability condition
from Definition \ref{def:ext}
and the regularity condition
from Definition \ref{def:reg}.

\begin{cor}
Theorem \ref{main+} implies Theorem \ref{main}.
\end{cor}

\nib{Proof.}
It suffices to show that
the hypotheses of Theorem \ref{main} 
(we choose $\aA = (2b)^{-1} h^{-3}$)
imply those of Theorem \ref{main+}.
This follows from Lemma \ref{ext}.
Indeed, if $G$ is $(c,h)$-typical with
$d(G)>n^{-(2b)^{-1}h^{-3}}$ and $c < c_0 d(G)^{h^2}$
then $G$ is $(\tfrac{1}{2}d(G)^h,h)$-extendable 
and $(K^r_q,Qc,q!^{-1}d(G)^Q)$-regular,
so $(\oO,h)$-extendable and $(K^r_q,Qc,\oO)$-regular 
with $\oO = q!^{-1}d(G)^h > n^{-b^{-1}h^{-2}}$
and $Qc < Q c_0 d(G)^{h^2} < c'_0 \oO^h$,
for some $c'_0=c'_0(q)$. \qed

\medskip

We will also need the following estimate
on the number of extensions that use
an edge from some bounded $r$-graph $J$.

\begin{lemma} \label{extbdd}
Let $E=(\phi,F,H)$ be an extension.
Suppose $J \sub K^r_n$ is $c$-bounded.
Then \[ |\{\phi^* \in X_E(K^r_n): 
\phi^*(H\sm H[F]) \cap J \ne \es \}| 
< c|H\sm H[F]|n^{v_E}. \]
\end{lemma}

\nib{Proof.}
Fix any $e \in H\sm H[F]$. Let $r'=|e \sm F|$.
As $J$ is $c$-bounded, there are
at most $cn^{r'}$ choices of the restriction 
$\phi^*\!\mid_e$ of $\phi^*$ to $e$
such that $\phi^*(e) \in J$. 
Each such choice has fewer than $n^{v_E-r'}$ 
extensions to $\phi^* \in X_E(K^r_n)$.
Summing over $e$ proves the lemma. \qed
 
\medskip

Next we turn to typicality properties of random $r$-graphs. 
We write $L \sim K^r_n(\nu)$ to mean that 
$L$ is a binomial random hypergraph
where each $e \in K^r_n$ is independently
included in $L$ with probability $\nu$.
The following lemma shows that
random $r$-graphs are w.h.p.\ typical.

\begin{lemma} \label{extrandom}
Suppose $L \sim K^r_n(\nu)$, where $\nu > n^{-1/3s}$.
Then w.h.p.\ $L$ is $(n^{-1/9},s)$-typical.
\end{lemma}

\nib{Proof.}
By a Chernoff bound w.h.p.\ $d(L) = \nu + O(n^{-0.4})$.
Let $E=(\phi,F,H)$ be any extension with $|H| \le s$. 
Note that $\mb{E}X_E(L) = (1+O(n^{-1})) \nu^{e_E} n^{v_E}$.
Also, for any $k \in [r]$ there are $O(n^k)$ edges $e \in K^r_n$
with $|e \sm \phi(F)|=k$, and for each such $e$, changing whether 
$e \in L$ affects $X_E(L)$ by $O(n^{v_E-k})$.
Thus $X_E(L)$ is $O(n^{2v_E-1})$-varying,
so by Lemma \ref{lip1} w.h.p.\
$X_E(L) = (1 \pm n^{-1/9}) d(L)^{e_E} n^{v_E}$. \qed

\medskip

We conclude this section by defining a refined
notion of boundedness that operates with respect
to all small extensions in some $r$-graph $L$.
The lemma following the definition shows
that if $J$ is bounded and has no `heavy' edges
and $L$ is random then w.h.p.\ $J$ is bounded w.r.t.\ $L$.

\begin{defn}
Let $E=(\phi,F,H\sm\{e\})$ be an extension,
$L \sub K^r_n$ and $J \in \mb{Z}^{K^r_n}$.
Define $X^e_E(L,J) = \sum_{\phi^* \in X_E(L)} |J_{\phi^*(e)}|$.
We say that $J$ is $(\tT,s)$-bounded w.r.t.\ $L$
if $X^e_E(L,J) < \tT d(L)^{e_E} n^{v_E}$
for any extension $E=(\phi,F,H\sm\{e\})$ with $|H| \le s$
and $e \in H \sm H[F]$.
\end{defn}

\begin{lemma} \label{Jrandom}
Suppose $J \in \mb{Z}^{K^r_n}$ 
is $\tT$-bounded with $\tT > n^{-1/20}$
and $|J_e| < n^{0.1}$ for all $e \in K^r_n$.
Let $L \sim K^r_n(\nu)$, where $\nu > n^{-1/3s}$.
Then w.h.p.\ $J$ is $(1.1\tT, s)$-bounded w.r.t.\ $L$.
\end{lemma}

\nib{Proof.}
By a Chernoff bound w.h.p.\ $d(L) = \nu + O(n^{-0.4})$.
Let $E=(\phi,F,H\sm\{e\})$ be an extension 
with $|H| \le s$ and $e \in H \sm H[F]$.
Write $X^e_E(L,J) = \sum_{\phi^* \in X_E(K^r_n)} 
1_{\phi^* \in X_E(L)} |J_{\phi^*(e)}|$.
As $J$ is $\tT$-bounded,
$\sum_{\phi^* \in X_E(K^r_n)} |J_{\phi^*(e)}| < \tT n^{v_E}$.
For each $\phi^* \in X_E(K^r_n)$ we have
$\mb{P}(\phi^* \in X_E(L)) = \nu^{e_E}$,
so $\mb{E}X^e_E(L,J) < \tT \nu^{e_E} n^{v_E}$.
For any $k \in [r]$, there are $O(n^k)$ choices 
of $f \in K^r_n$ with $|f \sm \phi(F)|=k$.
For each such $f$, changing whether $f \in L$
affects $X^e_E(L,J)$ by $O(n^{v_E-k+0.1})$.
Thus $X^e_E(L,J)$ is $O(n^{2v_E-0.8})$-varying, 
so by Lemma \ref{lip1} w.h.p.\ 
$X^e_E(L,J) < 1.1\tT d(L)^{e_E} n^{v_E}$. \qed

\section{Template} \label{sec:template}

In this section we construct the template,
and establish its combinatorial extendability properties.
(We defer the analysis of its algebraic extendability 
properties to Section \ref{sec:ab}.)
Henceforth, we fix $G$ as in the statement of Theorem \ref{main+}, 
and assume without loss of generality that 
$\oO < \oO_0(q)$ is sufficiently small,
so $G$ is a $K^r_q$-divisible $(K^r_q,c,\oO)$-regular 
$(\oO,h)$-extendable $r$-multigraph on $[n]$,
where $n^{-b^{-1}h^{-2}}<\oO<\oO_0(q)$, 
without loss of generality $c=\oO^h$,
and $n>n_0$ is sufficiently large.
(For convenient notation here we have slightly altered the statement
of Theorem \ref{main+}, assuming $c=\oO^h$ and $\oO<\oO_0(q)$
instead of $c< c_0 \oO^h$. This is valid 
as increasing $c$ or $\oO$ only makes our hypotheses weaker.
We will not compute an explicit bound for $\oO_0$ or $n_0$.)

\subsection{Overview}

As discussed in subsection \ref{sub:imp},
our algebraic model for designs will be the image 
$M(\mb{F}_{p^a}^r) \sub \mb{F}_{p^a}^q$ 
of a suitable matrix $M \in \mb{F}_p^{q \times r}$, 
as in Definition \ref{defgen} below.
The template will be a set of edge-disjoint $q$-cliques
determined by this algebraic model and 
various random choices described below.

Informally, these random choices are
\begin{enumerate}
\item activation events for each clique $Q'$, with probabilities proportional
to the weights $w_{Q'}$ satisfying the regularity condition for $G$,
\item random injections $\pi_e:e \to [q]$ for each edge%
\footnote{The clique exchange algorithm in Section \ref{sec:cea}
will operate with edges in the complete $r$-graph, not just in $G$.}
$e \in K^r_n$, describing the required position of $e$ in a certain vector in $\mb{F}_{p^a}^q$ 
corresponding to a template clique containing $e$, 
thus ensuring that any edge is covered by at most one template clique,
\item random injections $f_j: [n] \to \mb{F}_{p^a}$ for each $j \in [z]$,
via which we identify our basic combinatorial objects ($q$-cliques)
with basic algebraic objects (vectors in $\mb{F}_{p^a}^q$),
\item random choices $T_e \in [z]$ for each edge $e \in K^r_n$;
 we think of $e$ as \emph{following} the injection $f_{T_e}$.
\end{enumerate}
We will define the template $M^*$ as $\bigcup_{j \in [z]} M^*_j$,
where each $M^*_j$ consists of all activated $q$-cliques $Q'$
in which each edge $e$ follows $f_j$ and the image $q$-set 
$f_j(V(Q')) \sub \mb{F}_{p^a}$ can be ordered as a vector in  $\mb{F}_{p^a}^q$ 
where each image $r$-set $f_j(e)$ of an edge 
appears in the coordinates prescribed by $\pi_e$
and is \emph{non-degenerate}, that is, linearly independent over $\mb{F}_p$.

The point of using many independent injections $f_j$ 
is that it is unlikely that the template `aborts', meaning that 
for some edge $e$ there are few $j$ such that $f_j(e)$ is non-degenerate.
We assume without further comment throughout the paper
that the template does not abort.
Strictly speaking, we include the event `template aborts' 
in our union bound of all bad events for the template,
so all statements concerning the template of
the form `w.h.p.\ P' should be understood as
`w.h.p.\ P or the template aborts';
henceforth we will suppress such qualifications.

We will see that $M^*$ is a disjoint union of cliques,
so we can define its underlying $r$-graph $G^* = \bigcup_{j \in [z]} G^*_j$,
where each $G^*_j$ (the `colour $j$' edges of $G$)
is the union of the cliques in $M^*_j$.

We present the formal details of the construction in the next subsection.
In subsection \ref{sub:loc} we start the analysis by estimating the probability 
that any given edge belongs to the template. In itself this is a fairly routine calculation,
but in fact the precise setup and statement will be somewhat technical, 
as we will require these probabilities conditional on certain `local events'
containing information about the random choices made for other edges.
These estimates will be used throughout the paper, and in particular
in subsection \ref{sub:ext} to show extendability of (the underlying $r$-graph) of the template,
where we may consider `colour blind' or `rainbow' extensions.

\subsection{Construction} \label{sub:con}
 
We now proceed to the formal details of the construction outlined above.
We start by defining the algebraic model.

\begin{defn} \label{defgen}
Let $p$ be a prime\footnote{This exists by Bertrand's postulate.}
with $2^{8q}<p<2^{9q}$.
Let $M \in \mb{F}_p^{q \times r}$ be a $q \times r$
matrix with entries in $\mb{F}_p$.
We call $M$ {\em generic} if
every square submatrix of $M$ is nonsingular.
\end{defn}

To see that $M$ as in Definition \ref{defgen} exists,
consider a uniformly random choice of $M$.
For any fixed $j$ by $j$ submatrix, 
revealing its rows in sequence, 
the $i$th row is in the span of the previous rows
with probability at most $p^{i-1-j}$, so the 
matrix is singular with probability at most $2p^{-1}$.
Thus the required property fails with probability 
at most $2^{q+r+1}p^{-1} < 1$, so $M$ exists.

\medskip

Next we choose the dimension $a$ for which we will embed $[n]$ in $\mb{F}_{p^a}$,
so that the density $\gG$ of occupied vertices is of order $1/p$.
Let $a \in \mb{N}$ be such that $p^{a-2} < n \le p^{a-1}$.
We write \[\gG = np^{-a}, \ \ \text{ noting that } p^{-2} < \gG \le p^{-1}.\] 

We also recall that $G$ is $(K^r_q,c,\oO)$-regular, so there are weights
$w_{Q'} \in [\oO n^{r-q},\oO^{-1} n^{r-q}]$ for each $Q' \in K^r_q(G)$
with $\sum \{ w_{Q'} : e \in Q' \} = (1 \pm c) G_e$ for all $e \in K^r_n$. 

Now we list the sequence of independent random choices
used to construct the template.

\begin{enumerate} 
\item Let $A_{Q'}$ be independent Bernoulli variables 
for each $q$-clique $Q' \in K^r_q(G)$
with each $\mb{P}(A_{Q'}=1)=w_{Q'} \oO n^{q-r}$.
We say that $Q'$ is \emph{activated} if $A_{Q'}=1$.
\item Let $\pi_e:e \to [q]$ be independent random injections
for each edge $e \in K^r_n$ of the complete $r$-graph on $[n]$;
these will correspond to `positions' of edges in vectors.
\item Let $f = (f_j: j \in [z])$, with%
\footnote{We use a different letter here for clarity.}
$z=h$, where we choose independent uniformly random
injections $f_j: [n] \to \mb{F}_{p^a}$;
these identify combinatorial objects with algebraic objects.
\item We choose $T_e \in [z]$ for all $e \in K^r_n$ 
independently and uniformly at random;
we think of $e$ as `following' the embedding $f_{T_e}$.
\end{enumerate}

Given $f$, for each $e \in K^r_n$ we let\footnote{
Recall that we regard $\mb{F}_{p^a}$ as a vector space 
over $\mb{F}_p$ and for $X \sub \mb{F}_{p^a}$
we write $\dim(X)$ for the dimension of
the subspace spanned by the elements of $X$.} 
\[\mc{T}_e = \{j \in [z]: \dim(f_j(e))=r\}.\]
We \emph{abort} if any $|\mc{T}_e| \le z-2r$, which occurs with probability
at most $\tbinom{n}{r} \tbinom{z}{2r} (p^r/n)^{2r} = O(n^{-r})$.

We say $Q' \in K^r_q(G)$ is \emph{compatible} with $j$
if we can write $Q'=\phi(Q)$ for some
injection $\phi:[q] \to [n]$ such that
$T_e=j \in \mc{T}_e$ for all $e \in Q'$,
and for some $y \in \mb{F}_{p^a}^r$ we have
$f_j(\phi(i)) = (My)_i$ for all $i \in [q]$.

For future reference, we note the following lemma, which implies that
for any non-degenerate embedding $f_j\vert_e$ of an edge $e$ in some clique $Q'=\phi(Q)$
there is a choice of $f_j\vert_{V(Q')}$ such that if $T_{e'}=j$ for all $e' \in Q$
then  $Q'$ is compatible with $j$.

\begin{lemma} \label{lem:existcompat}
Let $e \in Q'=\phi(Q) \in K^r_q(G)$ and $y \in \mb{F}_{p^a}^r$.
Suppose $j \in \mc{T}_e$ and $f_j(\phi(i)) = (My)_i$ for all $\phi(i) \in e$.
Then $j \in \mc{T}_{e'}$ for all $e' \in Q'$.
In particular, $My$ has distinct coordinates.
\end{lemma}

\begin{proof}
As $j \in \mc{T}_e$ we have $\dim(f_j(e))=r$.
Write $f_j(e) = ((My)_i: i \in I)$, where $I=\phi^{-1}(e)$.
Then $y = M_I^{-1} f_j(e)$, where $M_I$ is the square submatrix of $M$
with rows indexed by $I$, which is nonsingular by Definition \ref{defgen}.
For any $e' = \phi(I') \in Q'$ we have $f_j(e') = M_{I'} y = M_{I'} M_I^{-1} f_j(e)$.
We claim that $j \in \mc{T}_{e'}$, that is, $\dim(f_j(e'))=r$.
Indeed, if this fails then $a \cdot f_j(e') = 0$ for some nonzero $a \in \mb{F}_p^r$.
However, then $a' := a M_{I'} M_I^{-1} \ne 0$ by nonsingularity of $M_{I'} M_I^{-1}$,
so $a' \cdot f_j(e) = 0$ contradicts $\dim(f_j(e))=r$. The claim follows.
In particular, $My$ has distinct coordinates.
\end{proof}

Let $\pi = (\pi_e: e \in K^r_n)$ 
where we choose independent uniformly random
injections $\pi_e:e \to [q]$.
We say $Q' \in K^r_q(G)$ is \emph{compatible} with $\pi$
if there is some bijection $\phi:V(Q) \to V(Q')$ 
such that $\pi_e\phi(i)=i$ whenever $\phi(i) \in e \in Q'$
(for brevity we write this as $\pi_e\phi=id$).
Note that if such $\phi$ exists then it is unique
(we have $\phi(i)=\pi_e^{-1}(i)$ for all $e \ni i$),
so we can write $Q'=\phi(Q)$ unambiguously,
and we will often identify template cliques with such embeddings $\phi$.

Now we define the template;
the lemma following the definition shows that
it is an edge-disjoint union of compatible cliques.

\begin{defn} \label{def:template} $ $
\begin{enumerate} 
\item
Let $M^*_j$ for $j \in [z]$ be the set of all activated 
$q$-cliques compatible with $\pi$ and $j$.
\item
The template is $M^* = \cup_{j \in [z]} M^*_j$.
\item
The underlying $r$-graph of the template is
$G^* = \cup_{j \in [z]} G^*_j$, where $G^*_j := \bigcup M^*_j$.
\end{enumerate} 
\end{defn}

In informal discussions below we will often abuse language 
and refer to the underying $r$-graph $G^*$ simply as `the template'.
By the following lemma, each $G^*_j$ 
is an $r$-graph (with no multiple edges), consisting
of all $e \in K^r_n$ such that $e \in Q'$ for some $Q' \in M^*_j$.
As $T_e=j$ for all $e \in G^*_j$,
we have $G^*_1,\dots,G^*_z$ edge-disjoint.

\begin{lemma}
$M^*$ is a clique decomposition of $G^*$.
\end{lemma}

\nib{Proof.}
It suffices to show for fixed $j \in [z]$ that
any $e \in G^*_j$ belongs to a unique clique $Q' \in M^*_j$.
To see this, note that
as each square submatrix of $M$ is nonsingular,
there is a unique $y \in \mb{F}_{p^a}^r$ such that 
$(My)_i = f_j(x)$ for all $x \in e$, $i=\pi_e(x)$,
which determines $V(Q') = f_j^{-1}(My)$. \qed

\medskip

We conclude with some further notation 
that will be used in the analysis of the template.

\begin{defn}
For $e \in G^*$ let $M^*(e) \in M^*$ be 
the $q$-clique such that $e \in M^*(e)$.

For $J \sub G^*$ let $M^*(J) = \sum_{e \in J} M^*(e) \in \mb{N}^{G^*}$.
\end{defn}

\subsection{Local events}  \label{sub:loc}

Here we give estimates for edge probabilities, 
conditional on certain `local events' $\mc{E}^e$ for each $e \in K^r_n$
that determine whether $e$ is in the template.

Formally, writing $\OO$ for the probability space of the template,
for each $e \in K^r_n$ and $\oO \in \OO$ the local event $\mc{E}^e = \mc{E}^e(\oO)$ 
will be a subset of $\OO$ containing $\oO$
defined by specifying the values of certain random variables,
such that  $1_{e \in G^*}$ is constant ($0$ or $1$) on $\mc{E}^e$.
Together, the set of possible local events $\mc{E}^e$ for fixed $e$
will be the atoms of an algebra $\mc{F}^e$ on $\OO$.
We can construct $\mc{F}^e$ recursively in stages,
where in each stage we have some partition $P$ of $\OO$,
and if there is some cell $C$ of $P$ on which $1_{e \in G^*}$
is not constant then we refine $P$ by partitioning $C$
according to the values of some further random variables,
which we informally think of as being `revealed'.

An equivalent formulation, which we adopt in the following definition,
is to recursively construct $\mc{E}^e(\oO) \sub \OO$ for given $\oO$ 
by specifying the values of some sets $\mc{X}_0, \mc{X}_1, \dots$
of random variables until we obtain $\oO \in \mc{E}^e(\oO) \sub \OO$
on which $1_{e \in G^*}$ is constant such that the construction
satisfies $\mc{E}^e(\oO') = \mc{E}^e(\oO)$ for any $\oO' \in \mc{E}^e(\oO)$.
All random variables in the following definition (e.g.~$T_e=T_e(\oO)$)
are functions of $\oO \in \OO$, which we usually suppress in our notation.\footnote{
So it is unlikely that $\oO \in \OO$ will be confused with $\oO>0$
appearing in the statement of Theorem \ref{main+}.
 }
For expository purposes we will say everything twice,
first intuitively and then formally.

\begin{defn} \label{local} (Local events)

Suppose $e \in K^r_n$. 
If $G_e=0$ we note that $e \notin G^*$ and let 
$\mc{E}^e$ be the trivial event that always holds.

Suppose $G_e>0$, reveal $T_e=j$ and $f_j\vert_e=\aA$.
If $\dim(\aA)<r$ then $\mc{E}^e$ is the event
that $T_e=j$ and $f_j\vert_e=\aA$, which witnesses $e \notin G^*$.

(Formally, for any $\oO \in \OO$ we start with $\mc{X}_0 = \{T_e\}$, 
write $j:=T_e(\oO)$, let $\mc{X}_1 = \{T_e,f_j\vert_e\}$,
and write $\aA:=f_j\vert_e(\oO)$. For any $\oO$ such that
$\dim(\aA)<r$ then we have already defined the local event, 
namely $\mc{E}^e(\oO) = \{ \oO': T_e(\oO')=j, f_j\vert_e(\oO')=\aA \}$.
For any other $\oO$ we continue.)

Now suppose $\dim(\aA)=r$, reveal $\pi_e$,
and let $y \in \mb{F}_{p^a}^r$ 
with $f_j(x) = (My)_i$ for all $x \in e$, $\pi_e(x)=i$;
note that $y$ is unique as $M$ is generic.
We reveal $f_j^{-1}((My)_i)$ for all $i \in [q] \sm \pi_e(e)$.
If all $f_j^{-1}((My)_i)$ are defined%
\footnote{If $(My)_i \notin f_j([n])$ then $f_j^{-1}((My)_i)$ is undefined.}
we let $\phi:[q] \to [n]$ be such that%
\footnote{Recall our concatenation notation 
and that we identify vectors with functions.}
$f_j\phi=My$. 
If any $f_j^{-1}((My)_i)$ is undefined or $\phi(Q) \notin K^r_q(G)$
then this defines the local event $\mc{E}^e$ witnessing $e \notin G^*$.

(Formally, given $j,\aA$ as above with $\dim(\aA)=r$,
we enlarge $\mc{X}_1$ to $\mc{X}_2 = \{T_e,f_j\vert_e,\pi_e\}$,
use $\pi_e(\oO)$ to define $y=y(\oO)$ as above,
then enlarge $\mc{X}_2$ to $\mc{X}_3$ 
by including the random variables $f_j^{-1}((My)_i)$ 
for all $i \in [q] \sm \pi_e(e)$. If any $f_j^{-1}((My)_i)$ is undefined 
or $\phi(Q) \notin K^r_q(G)$ then we have already specified the local event, 
namely $\mc{E}^e(\oO) = \{ \oO': X(\oO')=X(\oO) \ \forall X \in \mc{X}_3\}$.
For any other $\oO$ we continue.)

Finally, suppose $\phi(Q) \in K^r_q(G)$,
reveal whether $\phi(Q)$ is activated,
and reveal $(T_{e'},\pi_{e'})$ 
for all $e' \in \phi(Q) \sm \{e\}$.
Then $\mc{E}^e$ is defined by all 
the random variables revealed so far,
which determine whether $e \in G^*$:
given $T_e=j$, $f_j\vert_e=\aA$ with $\dim(\aA)=r$,
$y$ obtained from $\pi_e$,
and $\phi$ obtained from all $f_j^{-1}((My)_i)$,
we have\footnote{This statement follows from Lemma \ref{lem:existcompat}.}
 $e \in G^*$ iff $\phi(Q)$ is activated and
$T_{e'}=j$ and $\pi_{e'}\phi=id$
for all $e' \in \phi(Q)$.

(Formally, given $X(\oO)$ for all $X \in \mc{X}_3$
such that $\phi$ as above is defined and $\phi(Q) \in K^r_q(G)$,
we enlarge $\mc{X}_3$ to $\mc{X}_4$ by including
the random variables $A_{\phi(Q)}$
and $(T_{e'},\pi_{e'})$ for all $e' \in \phi(Q)$.
The local event is then
$\mc{E}^e(\oO) = \{ \oO': X(\oO')=X(\oO) \ \forall X \in \mc{X}_4\}$.)

We say that a vertex $x$ is \emph{touched} by $\mc{E}^e$ if $f_{T_e}(x)$ is revealed%
\footnote{Formally, this means $G_e>0$ and $f_{T_e}(x)$ is the same for all $\oO \in \mc{E}^e$.} 
by $\mc{E}^e$.

We say that an edge $e'$ is \emph{touched} by $\mc{E}^e$
if $T_{e'}$ is revealed by $\mc{E}^e$.
\end{defn}

Note that if an edge $e'$ is touched by $\mc{E}^e$
then $e'=e$ or $\phi$ is defined and $e' \in \phi(Q)$,
so in particular all vertices of $e'$ are touched.
On the other hand, $\mc{E}^e$ can touch vertices 
of an edge without touching the edge.
The next lemma gives estimates for edge probabilities
in the template conditional on certain combinations of local events
(or with no conditioning if $S=\es$). 

For $S \sub K^r_n$, we let $\mc{F}^S$ be the algebra on $\OO$ 
generated by the algebras $\{ \mc{F}^f: f \in S\}$.
Each atom $\mc{E}^S$ of $\mc{F}^S$ is
a non-empty intersection of atoms $\mc{E}^f$ of $\mc{F}^f$ with $f \in S$.
Without loss of generality $\mc{E}^S = \bigcap_{f \in S} \mc{E}^f$,
as if for some $f \in S$ no atom of $\mc{F}^f$ is included 
then we can delete $f$ from $S$.
We say that a vertex or edge is touched 
by $\mc{E}^S = \cap_{f \in S} \mc{E}^f$
if it is touched by any $\mc{E}^f$ with $f \in S$.
Let 
\[ \rho := \oO z^{-Q} q! (q)_r^{-Q} \gG^{q-r}. \]

\begin{lemma} \label{e|E}
Let $S \sub K^r_n$ with $|S|<h=z$
and $\mc{E}^S = \cap_{f \in S} \mc{E}^f$ be an atom of $\mc{F}^S$.
Suppose $e \in K^r_n$ is not touched by $\mc{E}^S$
and $j \in [z] \sm \{T_f: f \in S\}$.
Then $\mb{P}(e \in G^*_j \mid \mc{E}^S)
= (1 \pm 1.1c) \rho G_e$.
\end{lemma}

\nib{Proof.}
We fix any $e \in Q' \in K^r_q(G)$
and estimate the probability
that $e \in G^*_j$ with $M^*(e) = Q'$.
Throughout we exclude any clique $Q'$ with 
some $e' \in Q'$ touched by $\mc{E}^S$;
there are $O(n^{q-r-1})$ such cliques,
as there are $O(1)$ choices of $e'$,
and $|e \cup e'| \ge r+1$.
We activate any clique $Q'$ 
with probability $w_{Q'} \oO n^{q-r}$. 
The probability that $T_{e'}=j$
for all $e' \in Q'$ is $z^{-Q}$.
We fix one of the $q!$ labellings $Q'=\phi(Q)$
and condition on $\pi_{e'}$ for all $e' \in \phi(Q)$
such that $\pi_{e'}\phi=id$; this occurs with
probability $(q)_r^{-Q}$.
We condition on $f_j\!\mid_e$ such that $\dim(f_j(e))=r$;
as $j \in [z] \sm \{T_f: f \in S\}$
this occurs with probability $1-O(n^{-1})$.
As $M$ is generic, there is a unique $y \in \mb{F}_{p^a}^r$ 
such that $(My)_i = f_j(x)$ for all $x \in e$, $i=\pi_e(x)$.
By Lemma \ref{lem:existcompat},
for any $I \in Q$ we have $\dim((My)_i: i \in I)=r$
and $My$ has distinct coordinates.
With probability $(1+O(n^{-1}))(p^{-a})^{q-r}$
we have $f_j(\phi(i))=(My)_i$ 
for all $i \in [q] \sm \pi_e(e)$.
Multiplying the probabilities, 
recalling $\gG = np^{-a}$, we obtain
\[ \mb{P}(M^*(e) = Q' \mid \mc{E}^S)
= (1+O(n^{-1})) w_{Q'} \oO n^{q-r} z^{-Q} 
q! (q)_r^{-Q} (p^{-a})^{q-r} 
= (1+O(n^{-1})) w_{Q'} \rho.\]
Summing over $Q'$, recalling
$\sum \{ w_{Q'} : e \in Q' \} = (1 \pm c) G_e$,
gives $\mb{P}(e \in G^*_j \mid \mc{E}^S)
= (1 \pm 1.1c) \rho G_e$. \qed

\begin{rem} \label{rem:e|E}
The proof of Lemma \ref{e|E} also shows 
for any $j \in [z] \sm \{T_f: f \in S\}$,
$v: e \to \mb{F}_{p^a}$ with $\dim(v)=r$,
and injection $\pi: e \to [q]$ that
\begin{enumerate}
\item  $\mb{P}(e \in G^*_j \mid \mc{E}^S \cap \{f_j\!\mid_e=v\})
= (1 \pm 1.1c) \rho G_e$,
\item $\mb{P}(e \in G^*_j \mid \mc{E}^S \cap \{T_e=j\})
= (1 \pm 1.1c) z \rho G_e$,
\item $\mb{P}( \{e \in G^*_j\} \cap \{\pi_e=\pi\} \mid \mc{E}^S)
= (1 \pm 1.1c) (q)_r^{-1} \rho G_e$.
\end{enumerate}
\end{rem}

\subsection{Extendability} \label{sub:ext}

We conclude this section by showing that
(the underlying $r$-graph $G^*$ of) the template is w.h.p.\ extendable. 
This will be deduced from extendability of $G$ 
and the estimates for edge probabilities conditional on local events
obtained above in Lemma \ref{e|E}. The idea of the proof is that
for any extension $\phi^+$ in some $X_E(G)$ we can bound  
the probability of $\phi^+ \in X_E(G^*)$ by applying Lemma \ref{e|E} to each edge in turn,  
using a new colour for each edge to avoid dependencies on previous edges.
This gives a bound for the expectation, from which we obtain a w.h.p.~statement
by standard concentration inequalities. Note that using a new colour for each edge
is useful both as a proof device and for later applications where we need `rainbow' extensions.

\begin{lemma} \label{ext*}
Suppose $E=(\phi,F,H)$ is an extension with $|H| \le z/3$.
Then w.h.p.\ $X_E(G^*) > \oO n^{v_E} (z\rho/2)^{e_E}$.
\end{lemma}

\nib{Proof.} 
We consider the modified variable $X'$ 
defined as the set or number of $\phi^+$ in $X_E(G^*)$
such that $(V(M^*(e)) \sm e) \cap \phi^+(V(H)) = \es$ for all $e \in \phi^+(H \sm H[F])$.
We claim that $X_E(G^*) \ge X' \ge X_E(G^*)  - O(n^{v_E-1})$,
deterministically given any outcome of the random choices in the template.
Here the first inequality is trivial. To see the second inequality, 
note that for any $\phi^+$ counted by $X_E(G^*)$ but not $X'$
we can list the vertices of $V(H)$ in some order 
such that some $e \in \phi^+(H \sm H[F])$ determines
(via $f_{T_e}$ and $\pi_e$) some $q$-clique $M^*(e)$
with vertices $f_{T_e}^{-1}(My)$ for some $y \in \mb{F}_{p^a}^r$
that contains some vertex $\phi^+(u)$ with $u$ earlier 
in the order than the last $v$ with $\phi^+(v) \in e$.
There are $O(1)$ choices for the roles of $u,v,e$ in $H$
and the position for $v$ to appear in $f_{T_e}^{-1}(My)$.
Such a choice determines a linear equation with coefficients in $\mb{F}_p$
for $f_{T_e}(v)$ in terms of $f_{T_e}(v')$ with $v'$ earlier than $v$,
thus determining $v$ uniquely given the earlier choices.
There are  $O(n^{v_E-1})$ choices for the remaining vertices
of the extension, so the claim follows.

Thus it suffices to prove the stated w.h.p.~lower bound for $X'$ instead of $X_E(G^*)$.
We will do so conditional on any event $\mc{E}_0$ of the form
$\{ f_j(x) = a_{j,x}\  \forall j \in [z], x \in \phi(F) \}$
on which the template does not abort.
As $G$ is $(\oO,h)$-extendable there are at least $\oO n^{v_E}$
choices of $\phi^+ \in X_E(G)$.  We fix any such $\phi^+$
and estimate $\mb{P}(\phi^+ \in X')$
by repeated application of Lemma \ref{e|E}.
Consider any $e=\phi^+(f)$ with $f \in H \sm H[F]$
and let $\mc{E}$ be the intersection of $\mc{E}_0$ 
and the local events of all previously considered edges,
that is, $\mc{E}=\mc{E}_0 \cap \mc{E}^S$ where $S$ is the set 
of previously considered edges and $\mc{E}^S = \mc{E}^S(\oO)$ 
is the atom of $\mc{F}^S$ containing the random element 
$\oO$ of the probability space of the template.

We may assume $(V(M^*(e')) \sm e') \cap \phi^+(V(H)) = \es$ 
for all $e' \in S$, and so $e$ is not touched by $\mc{E}^S$.
As the template does not abort on $\mc{E}_0$,
there are at least $2z/3$ choices of $j \in [z]$
not used by any previous edge such that Lemma \ref{e|E} 
with Remark \ref{rem:e|E}.i applies 
to give $\mb{P}(e \in G^*_j \mid \mc{E}) > 0.9 \rho$.
Furthermore, similarly to above, with probability $1-O(1/n)$
under the choice of the local event $\mc{E}^e$
we have $(V(M^*(e)) \sm e) \cap \phi^+(V(H)) = \es$.
Multiplying all conditional probabilities 
$2z/3 \cdot 0.9\rho - O(1/n) > 0.6z\rho$
and summing over $\phi^+$
gives $\mb{E}X_E(G^*) > \oO n^{v_E} (0.6 z\rho)^{e_E}$.

Next we show concentration.
First we show concentration of $\mb{E}[X' \mid f]$, 
which is the conditional expectation where we reveal 
the embeddings $f =(f_j: j \in [z])$ 
(consistently with $\mc{E}_0$)
but not the other random choices 
in the construction of the template.
Changing any $f_j(x)$ with $x \notin \phi(F)$
affects $X'$ by $O(n^{v_E-1})$, 
so $\mb{E}[X' \mid f]$ is $O(n^{2v_E-1})$-varying
in the sense of Definition \ref{def:genlip},
applied with $a'=(a'_1,\dots,a'_z)=(n,\dots,n)$,
$a=(a_1,\dots,a_z)=(p^a,\dots,p^a)$,
$b=(O(n^{v_E-1}),\dots,O(n^{v_E-1}))$,
$B = \sum a'_i b_i^2 = O(n^{2v_E-1})$. 
Applying Lemma \ref{lip3} with $t=.01\mb{E}X'$,
we have $\mb{P}(|\mb{E}[X' \mid f]-\mb{E}X'|>t)
 \le e^{-\OO(n^{v_E})^2/O(n^{2v_E-1})}$, so w.h.p.\
$\mb{E}[X' \mid f] > \oO n^{v_E} (0.59z\rho)^{e_E}$.

Now we fix $f$ consistent with $\mc{E}_0$ such that
$\mb{E}[X' \mid f] > \oO n^{v_E} (0.59z\rho)^{e_E}$
and show concentration of $X'$ 
under the remaining independent random choices 
during the construction of the template,
namely the activation of cliques and choices
of $T_e$ and $\pi_e$ for each $e \in K^r_n$.
We classify $e \in G$ according to the possible values
of $|Im(\phi') \cap \phi(F)|$ where $e \in \phi'(Q) \in K^r_q(G)$
and there is some $y \in \mb{F}_{p^a}^r$ and $j \in [z]$
with $f_j(\phi'(i)) = (My)_i$ for all $i \in [q]$.
Given $s \in [r]$, there are 
$O(n^s)$ such $\phi'$ with $|Im(\phi') \cap \phi(F)|=r-s$,
changing whether $\phi'(Q)$ is activated or
any $T_e$ or $\pi_e$ for $e \in \phi'(Q)$
affects $X'$ by $O(n^{v_E-s})$,
and all other choices do not affect $X'$.
Thus\footnote{
The applications of Definition \ref{def:genlip}
and Lemma \ref{lip3} here and henceforth 
are similar to those in the previous paragraph,
so we will not spell them out in full detail.
One should note that we use Definition \ref{def:genlip}
in more generality here than above, as the vectors $a,a',b$
are indexed by all the remaining independent random choices
and we cannot take $b$ to be a constant vector.
} 
$X'$ is $B$-varying, where
$B = \sum_{s \in [r]} O(n^s) O(n^{v_E-s})^2 = O(n^{2v_E-1})$, 
 so by Lemma \ref{lip3} w.h.p.\ 
$X_E(G^*) \ge X' > \oO n^{v_E} (z\rho/2)^{e_E}$, as required. \qed

\begin{rem} \label{rem:X'}
Let $X'_E(G^*)$ be the set or number of 
$\phi^+ \in X_E(G^*)$ that are `rainbow',
meaning that $j \ne j'$ 
whenever $\{f,f'\} \sub H \sm H[F]$,
$\phi^+(f) \in G^*_j$, $\phi^+(f') \in G^*_{j'}$.
The proof of Lemma \ref{ext*} shows
w.h.p.\ $X'_E(G^*) > \oO n^{v_E} (z\rho/2)^{e_E}$.
Furthermore, by Remark \ref{rem:e|E}.iii the proof also shows
that given any fixed injections $\pi_f: f \to [q]$
for all $f \in H \sm H[F]$ we can find
at least $\oO n^{v_E} (z\rho/2(q)_r)^{e_E}$
choices of $\phi^+ \in X'_E(G^*)$ with
$\pi_e \phi^+\!\mid_f = \pi_f$ for all
$e=\phi^+(f)$ with $f \in H \sm H[F]$.
\end{rem}

\section{Approximate decomposition} \label{sec:approx}

In this section we complete the template $G^*$
 to an approximate decomposition,
namely a set of cliques such that every edge
is covered once or twice, and the edges covered twice
form a suitably bounded subgraph of the template.
Throughout the paper, we think of the template
as a deterministic object that satisfies
all w.h.p.\ statements that we make about it,
except that for convenience of exposition
we have deferred some of these w.h.p.\ statements
to the places where they are used.
In particular, we assume henceforth that $G^*$ satisfies
the extendability conditions established in Lemma \ref{ext*},
so e.g.~in Lemma \ref{cover} below we omit `w.h.p.'
when referring to these properties. By contrast,
in Lemma \ref{nibble} below we retain `w.h.p.'
as we are establishing a new w.h.p.~property of $G^*$.

\subsection{Nibble} \label{sec:nibble}

Here we show how to partition almost all 
of the multigraph $G-G^*$ into $q$-cliques.
This will follow from the approximation decomposition result Lemma \ref{nibble+},
although to make this work we need to pass to 
a carefully chosen subcollection $\mc{Q}'$ of the $q$-cliques $K^r_q(G)$ in $G$,
where we only consider rainbow cliques to eliminate dependencies,
and we also randomly subsample cliques according to a `rejection sampling' distribution
that corrects for biases introduced by the template construction
depending on the edge multiplicities $G_e$.

\begin{lemma} \label{nibble}
w.h.p.~the template $G^*$ is such that there is a set of $q$-cliques $M^n$ 
such that the {\em leave} $L := G - G^* - \sum M^n$ 
is $c_1$-bounded and $L_e \ge 0$ for all $e \in G$.
\end{lemma}

\nib{Proof.}
We will apply Lemma \ref{nibble+} with $G-G^*$
in place of $G$ and some $\mc{Q}'$ in place of $\mc{Q}$.
We construct $\mc{Q}' \sub K^r_q(G)$ as follows,
using some randomness from the template probability space $\OO$
and an additional independent random choice.
Consider any $Q' \in K^r_q(G)$ and reveal the local events 
$\mc{E}^e=\mc{E}^e(\oO)$ for each $e \in Q'$,
where $\oO \in \OO$ is the random element determining the template.
If $Q'$ is not activated or $T_e=T_{e'}$ for any $e \ne e'$ in $Q'$
then we do not include $Q'$ in $\mc{Q}'$.
For $v \in \{0,1\}^{Q'}$ let $\mc{E}^{Q'}_v$ be the 
event that $v_e = 1_{e \in G^*}$ for all $e \in Q'$.
If $Q'$ is activated,
all $T_e$ for $e \in Q'$ are distinct and $\mc{E}^{Q'}_v$ holds 
then we include $Q'$ in $\mc{Q}'$ independently with probability 
$p^v_{Q'} = \prod_{e \in Q'} (1-v_e/G_e)$.
Note that if $G_e=G^*_e=1$ for some $e \in Q'$
then $p^v_{Q'}=0$, so $\mc{Q}' \sub K^r_q(G-G^*)$. 

Now we fix any $e' \in G$ and estimate the number $|\mc{Q}'(e')|$
of cliques in $\mc{Q}'$ containing $e'$.
We consider any activated $Q'$ with
$e' \in Q' \in K^r_q(G)$ and condition on 
the local event $\mc{E}^{e'}$ and any event 
$\mc{C} = \cap_{e \in Q'} \{T_e=j_e\}$ 
such that all $j_e$ are distinct
(the latter occurs with probability $(z)_Q z^{-Q}$).

For any $v  \in \{0,1\}^{Q'}$ with $v_{e'}=G^*_{e'}=1_{e' \in G^*}$, 
by repeated application of Lemma \ref{e|E} 
(with Remark \ref{rem:e|E}.ii), we have
$\mb{P}(\mc{E}^{Q'}_v \mid \mc{E}^{e'} \cap \mc{C} )
=  \prod_{e \in Q' \sm \{e'\}} ((1 \pm 1.1c)p^{v_e}_e)
= ( 1 \pm Qc) \prod_{e \in Q' \sm \{e'\}} p^{v_e}_e$, where 
$p^1_e = z\rho G_e$ and $p^0_e = 1-z\rho G_e$. Then 
\begin{align*}
\mb{P}( Q' \in \mc{Q}' \mid \mc{E}^{e'} \cap \mc{C} )
& = \sum_{v: v_{e'}=G^*_{e'}} 
\mb{P}( Q' \in \mc{Q}' \mid \mc{E}^{e'} \cap \mc{C} \cap \mc{E}^{Q'}_v )
\mb{P}(\mc{E}^{Q'}_v \mid \mc{E}^{e'} \cap \mc{C} ) \\
& = \sum_{v: v_{e'}=G^*_{e'}} 
( 1 \pm Qc)  p^v_{Q'} \prod_{e \in Q' \sm \{e'\}} p^{v_e}_e \\
& =  ( 1 \pm Qc) (1-G^*_{e'}/G_{e'})
 \prod_{e \in Q' \sm \{e'\}} \sum_{v_e \in \{0,1\}}  (1-v_e/G_e) p^{v_e}_e \\
& = ( 1 \pm Qc) (1-G^*_{e'}/G_{e'}) (1-z\rho)^{Q-1},
\end{align*}
as $(z\rho G_e )(1-1/G_e) + (1-z\rho G_e)(1-0/G_e)=1-z\rho$
for any $G_e>0$.

Recalling that we activate $Q'$ independently with
probability $w_{Q'} \oO n^{q-r}$ and
$\sum \{ w_{Q'} : e' \in Q' \} = (1 \pm c) G_{e'}$ we have
$\mb{E}[|\mc{Q}'(e')| \mid \mc{E}^{e'}] = (z)_Q z^{-Q} \cdot
( 1 \pm Qc) (1-G^*_{e'}/G_{e'}) (1-z\rho)^{Q-1} \cdot
(1 \pm c) \oO n^{q-r} G_{e'}
= (1 \pm 1.1Qc) d' n^{q-r} (G_{e'}-G^*_{e'})$,
where $d' = (z)_Q z^{-Q} \oO (1-z\rho)^{Q-1}$.

To show concentration of 
$X := |\mc{Q}'(e')| \mid \mc{E}^{e'}$ we apply 
Lemma \ref{lip3}, similarly to the proof of Lemma \ref{ext*},
with appropriate modifications for the conditioning.
Let $U$ be the set of vertices touched by $\mc{E}^{e'}$.
To see concentration of $\mb{E}[X \mid f]$, 
we note that changing any $f_j(x)$ with $x \notin U$
affects $X$ by $O(n^{q-r-1})$, 
so $\mb{E}[X \mid f]$ is $O(n^{2(q-r)-1})$-varying,
and by Lemma \ref{lip3} w.h.p.\ 
$\mb{E}[X \mid f] = (1 \pm c) \mb{E}X$.

Now we fix such $f$ and show concentration
under the remaining random choices.
We classify $e \in G$ according to the possible values
of $|Im(\phi') \cap U|$ where $e \in \phi'(Q) \in K^r_q(G)$
and there is some $y \in \mb{F}_{p^a}^r$ and $j \in [z]$
with $f_j(\phi'(i)) = (My)_i$ for all $i \in [q]$.
Given $s \in [r]$, there are 
$O(n^s)$ such $\phi'$ with $|Im(\phi') \cap U|=r-s$,
and changing whether $\phi'(Q)$ is activated or
any $T_e$ or $\pi_e$ for $e \in \phi'(Q)$ untouched by $\mc{E}^{e'}$
affects $|\mc{Q}'(e')|$ by $O(n^{q-r-s})$.
Thus $X$ is $O(n^{2(q-r)-1})$-varying,
so by Lemma \ref{lip3} w.h.p.\ on any local event $\mc{E}^{e'}$
we have $|\mc{Q}'(e')| = (1 \pm 2Qc) d' n^{q-r} (G_{e'}-G^*_{e'})$.

As $c_1=(2Qc)^{1/2Q}$, Lemma \ref{nibble} now follows
from Lemma \ref{nibble+}. \qed

\subsection{Cover}

To complete the approximate decomposition,
we will cover the leave $L$ by a set $M^c$ of $q$-cliques,
each of which has one edge in $L$ 
and all remaining edges in $G^*$.
These remaining edges
$S := (\sum M^c) - L \sub G^*$
constitute the {\em spill} referred
to in Subsection \ref{sec:strategy}.
We require that $M^*(S) \sub G^*$ is a set,
that is, $S$ uses at most one edge 
of any given clique in the template $M^*$.
We also require that $M^*(S)$ is bounded 
in the sense of the following algebraic condition
that implies boundedness in the sense 
of Definition \ref{bdd} (consider lines
where all but one coordinate is fixed,
e.g.\ $v=e_i$ for some $i$).

\begin{defn} \label{def:linbdd}
For $I \in Q$ and $u,v$ in $\mb{F}_{p^a}^I \sm \{0\}$ 
we call $L = \{ u + cv : c \in \mb{F}_{p^a}\}$ 
an \emph{$I$-line}.
 
Suppose $J \sub G^*$. 
We say $J$ is linearly \emph{$\tT$-bounded} 
if for each $j \in [z]$, $I \in Q$ and $I$-line $L$
at most $\tT p^a$ edges $e \in J \cap G^*_j$
have $\pi_e(e)=I$ and $f_j(e) \in L$,
regarding $f_j(e)$ in $\mb{F}_{p^a}^I$ via 
$f_j(e)_i = f_j \pi_e^{-1} (i)$.
\end{defn}

\begin{rem} \label{rem:linbdd}
When trying to simplify the proof
in the first version of this paper, 
in the second version
we replaced all linear boundedness
assumptions by boundedness assumptions.
However, as pointed out by Lisa Sauermann,
this made the proof incorrect: 
specifically, Lemma 6.10 of the second version is wrong,
as the parameter $r'$ in its statement should be
replaced by the parameter $r_{a\phi'}$ 
in Definition \ref{def:ee'} of the current version.
Although this correction has several knock-on effects
in Sections \ref{sec:ab} and \ref{sec:cea},
the proof here is otherwise 
quite similar to that in the second version.
\end{rem}

The following lemma is immediate from the observation
that any affine linear space of dimension at least one
can be partitioned into lines.

\begin{lemma} \label{linbdd+dim}
Suppose $J \sub G^*$ is linearly $\tT$-bounded.
Let  $j \in [z]$, $I \in Q$ and $A$ be an
affine linear subspace of $\mb{F}_{p^a}^I$
with $\dim A \ge 1$.
Then at most $\tT |A|$ edges $e \in J \cap G^*_j$
have $\pi_e(e)=I$ and $f_j(e) \in A$.
\end{lemma}


We conclude this section with the following
lemma that implements the cover step.

\begin{lemma} \label{cover}
Suppose $L$ is a $c_1$-bounded 
submultigraph of $G-G^*$.
Then there is a set $M^c$ of $q$-cliques,
each of which contains exactly one edge of $L$, with 
{\em spill} $S := (\sum M^c) - L \sub G^*$, such that 
$M^*(S)$ is a set,%
\footnote{As $M^*(S)$ is a multiset {\em a priori},
we are asserting here that no edge
has multiplicity greater than $1$.}
$c_2$-bounded and
linearly $c_2$-bounded 
(recall $c_2 = \oO^{-h/20Q} c_1$).
\end{lemma}

\nib{Proof.} 
We order $L$ as $(e_i: i \in [|L|])$, and apply a random
greedy algorithm to select $q$-cliques $(K_i: i \in [|L|])$.
Write $S_i = (\cup_{i'<i} K_{i'} ) \sm L$.
At step $i$, we let $K_i = \phi_i(Q)$ be a uniformly
random $q$-clique containing $e_i$ such that
$K_i \sm \{e_i\} \sub G^*$ and
$M^*(K_i)$ is a set disjoint from $M^*(S_i)$.
(If no such $\phi_i$ exists then we abort.)
Note that the disjointness condition is equivalent 
to $K_i \cap M^*(S_i) = \es$. 

To develop some intuition for this algorithm, 
it is helpful to first consider the simpler process of choosing 
$K'_i = \phi'_i(Q)$ ignoring all disjointness conditions
bar requiring each $M^*(K'_i)$ to be a set,
so that $K'_1,\dots,K'_{|L|}$ are independent.
We denote the number of choices for $\phi'_i$ by $X_q(e_i)$, 
and claim that 
\[ X_q(e_i) > \oO (z\rho/2)^{Q-1} n^{q-r} - O(n^{q-r-1}), \]
where we apply Lemma \ref{ext*} and claim that
$O(n^{q-r-1})$ choices are excluded
due to $M^*(K'_i)$ having some repeated edge.
To see this, we fix any $e \ne e' \in Q \sm \{e_i\}$
and bound the number of choices for which
$M^*(\phi'_i(e)) \cap M^*(\phi'_i(e')) \ne \es$.
We consider the choices of $\phi'_i(v)$
for $v \in Q \sm e_i$ sequentially in some order
ending with some $v^*$ in exactly one of $e,e'$.
Say $v^*$ is the last vertex of $e'$
and let $v'$ be the last vertex of $e$.
Choosing $\phi'_i(v')$ determines $M^*(\phi'_i(e))$,
and there are $O(1)$ such choices with $e' \in M^*(\phi'_i(e))$.
Similarly, choosing $\phi'_i(v^*)$ determines $M^*(\phi'_i(e'))$,
and there are $O(1)$ such choices with 
$M^*(\phi'_i(e)) \cap M^*(\phi'_i(e')) \ne \es$.
The claim follows.

For each $e \in G^*$ let 
$E_e = \sum_{i \in [|L|]} \mb{P}(e \in M^*(K'_i))$.
We claim that $E_e < (2q)^{2q} \oO^{-1} (z\rho/2)^{1-Q} c_1$.
To see this, we write $E_e = \sum_{e' \in M^*(e)}
\sum_{i \in [|L|]} \mb{P}(e' \in K'_i)$.
For any $i$ and $e' \in M^*(e)$ 
there are at most $q! n^{q-r-|e' \sm e_i|}$
choices of $\phi'_i$ such that $e' \in K'_i$, 
so $\mb{P}(e' \in K'_i) < q!n^{q-r-|e' \sm e_i|}  X_q(e_i)^{-1}
< 1.1 q!\oO^{-1} (z\rho/2)^{1-Q} n^{-|e' \sm e_i|}$.
Also, as $L$ is $c_1$-bounded, for any $r' \in [r]$
there are at most $\tbinom{r}{r'} c_1 n^{r'}$ choices of $i$ 
with $|e_i \sm e'| = r'$. Summing over $e'$ and $r'$
we deduce $E_e < Q \sum_{r' \in [r]} 1.1 q!\oO^{-1} 
(z\rho/2)^{1-Q} n^{-r'} \cdot \tbinom{r}{r'} c_1 n^{r'} 
< (2q)^{2q} \oO^{-1} (z\rho/2)^{1-Q} c_1$, as claimed.

Now for any $f \in \tbinom{[n]}{r-1}$ we have 
$|M^*(S)(f)| = \sum \{ 1_{e \in M^*(K'_i)}:
 i \in [|L|], f \sub e \}$ pseudobinomial
(a sum of bounded independent variables)
with mean at most 
$(2q)^{2q} \oO^{-1} (z\rho/2)^{1-Q}  c_1 n < c_2 n/2$,
so w.h.p.\ $M^*(S)$ is $c_2$-bounded by Bernstein's inequality
(see Remark \ref{rem:dom}).

We now turn to the analysis of the algorithm.
The idea is to show that w.h.p.\ in each step the 
disjointness condition $K_i \cap M^*(S_{i-1}) = \es$
forbids at most half of the possible choices,
so the estimates from the independent process hold
in the actual process up to a factor of two. 

For $i \in [|L|]$ we let $\mc{B}_i$ be the bad event 
that $M^*(S_i)$ is not $c_2$-bounded.
We define a stopping time%
\footnote{This means that each $\{ \tau \le i\}$ is an event determined 
by the history of the process up to step $i$.} 
$\tau$ as the smallest $i$ for which $\mc{B}_i$ 
holds or the algorithm aborts, 
or $\infty$ if there is no such $i$. 
It suffices to show w.h.p.\ $\tau=\infty$. 

We fix $i_0 \in [|L|]$ and bound $\mb{P}(\tau=i_0)$ as follows.
For any $i<i_0$, since $\mc{B}_i$ does not hold, 
$M^*(S_i)$ is $c_2$-bounded. Then by Lemma \ref{extbdd}
the condition $K_i \cap M^*(S_i) = \es$ forbids at most 
$Qc_2 n^{q-r} < \tfrac{1}{2}X_q(e_i)$ choices of $\phi_i$.

For each $e \in G^*$ let 
$r_e = \sum_{i<i_0} \mb{P}'(e \in M^*(K_i))$,
where $\mb{P}'$ denotes conditional probability
given the choices made before step $i$.
By the bound on excluded choices, 
$\mb{P}'(e \in M^*(K_i)) < 2 \mb{P}(e \in M^*(K'_i))$, 
so $r_e < 2E_e$.

Now consider any $f \in \tbinom{[n]}{r-1}$
and let $X = \sum_{i<i_0} X_i$, where $X_i = |M^*(K_i)(f)|$.
Then $\sum_{i<i_0} \mb{E}'X_i = \sum \{ r_e: f \sub e \}
 \le 2 (2q)^{2q} \oO^{-1} (z\rho/2)^{1-Q} c_1 n < c_2 n/2$. 
We deduce that $X$ is $(Q,c_2 n/2)$-dominated
(see Definition \ref{def:dom}) with respect
to the natural filtration $\mc{F}$ of the process
(each $X_i$ is $\mc{F}_i$-measurable with $0 \le X_i \le Q$),
so w.h.p.\ $X < c_2 n$ by Lemma \ref{dom}.

Finally, consider any $j \in [z]$, $I \in Q$ and $I$-line $L$,
and let $X = \sum_{i<i_0} X_i$, where 
$X_i = |\{e \in M^*(K_i): f_j(e) \in L\}|$.
Then $\sum_{i<i_0} \mb{E}'X_i = \sum \{ r_e: f_j(e) \in L \}
\le 2 (2q)^{2q} \oO^{-1} (z\rho/2)^{1-Q} c_1 p^a$, 
so $X$ is $(Q,c_2 p^a/2)$-dominated,
so w.h.p.\ $X < c_2 p^a$ by Lemma \ref{dom}.

Thus w.h.p.\ $M^*(S_i)$ is $c_2$-bounded 
and linearly $c_2$-bounded for all $i<i_0$,
so $\tau>i_0$. Taking a union bound over $i_0$, 
w.h.p.\ $\tau=\infty$, as required. \qed

\section{Integral decomposition} \label{sec:int}

In this section we prove an analogue of the
results of Graver and Jurkat \cite{GJ} and Wilson \cite{W4}
on integral decompositions in which we can also impose a boundedness requirement. 
Their results, stated as Lemma \ref{GJ/W} below, show that the $K^r_q$-divisibility conditions
which are \emph{necessary} for decomposition are \emph{sufficient} for integral decomposition
in $K^r_n$ with $n \ge q+r$, that is, any $K^r_q$-divisible vector in $K^r_n$ can be expressed
as an integral linear combination of $q$-cliques in $K^r_q(K^r_n)$. In the next subsection
we will introduce their method of `octahedral decomposition' and use it to prove Lemma \ref{GJ/W} 
(we include the proof for expository purposes, as it illustrates some ideas needed later).

The main result of the section is Lemma \ref{bddint} below, which adds a boundedness property:
if the $K^r_q$-divisible target vector $J \in \mb{Z}^{K^r_n}$ is $\tT$-bounded
then we can express $J$ using an integral combination of cliques $\Phi \in \mb{Z}^{K^r_q(K^r_n)}$ 
so that both the positive and negative contributions are $O_q(\tT)$-bounded.
It is convenient to only  consider boundedness rather than linear boundedness at this stage,
as we will be able to enforce linear boundedness later via the Clique Exchange Algorithm in Section \ref{sec:cea}.

We will prove Lemma \ref{bddint} in subsection \ref{sub:bddint},
using a result on bounded generation in subsection \ref{sub:bddgenII},
proved via a result in subsection \ref{sub:bddgen},
which extends the results of \cite{GJ, W4} 
on generating nullspaces of shadow operators by octahedra,
by showing that it suffices to use a `thin' set of octahedra,
meaning that every edge is in only constantly many octahedra of the generating set.
As mentioned in subsection \ref{sub:imp}, we introduced this method in 2018 as a simplification
of our 2014 argument, and in retrospect from 2024 it has independent interest
as an analogy of Refined Absorption.

\subsection{Octahedral decomposition} \label{sub:oct}

A key idea in \cite{GJ, W4}, which we will also use,
is `octahedral decomposition', which we will discuss 
in this subsection. We will make some definitions
and then state the main result of \cite{GJ, W4}.
The following definition concerns inclusion matrices
(standard objects in combinatorics) and their
associated $\mb{Z}$-linear maps. One can think
of the latter as `shadow operators'
(which is reflected in our notation $\pl_r$)
by analogy with the standard combinatorial $r$-shadow,
which maps any hypergraph $G$ 
to the $r$-graph consisting of all $r$-sets
that are contained in some edge of $G$.
Our operator $\pl_r$ is defined similarly,
except that we keep track of integer multiplicities:
$(\pl_r J)_e$ is the sum of $J_f$ 
over all $f$ containing $e$.
For example, if $J \in \mb{Z}^{K^s_n}$ is $\{0,1\}$-valued
(so $J$ can be identified with an $s$-graph) then $\pl_r J$
is the $r$-multigraph where the multiplicity
of each $e$ in the combinatorial $r$-shadow of $J$
counts edges of $J$ containing $e$.

\begin{defn} \label{def:inclusion}
Suppose $s \ge r \ge 0$ and $J \in \mb{Z}^{K^s_n}$.
We define $\pl_r J \in \mb{Z}^{K^r_n}$ by
$\pl_r J_e = \sum \{ J_f: e \sub f \in K^s_n \}$.
Equivalently, $\pl_r J = M^r_s(n) J$,
where $M^r_s(n)$ is the inclusion matrix
with rows indexed by $K^r_n$,
columns indexed by $K^s_n$,
and $ef$-entry $M^r_s(n)_{ef} = 1_{e \sub f}$.

We write $\pl_r J = \pl J$ 
if $r$ is clear from the context.
We apply the same notation to vectors
of $q$-cliques identifying $Q'$ with $V(Q')$:
for $\Phi \in \mb{Z}^{K^r_q(K^r_n)}$ we define
$\pl \Phi \in \mb{Z}^{K^r_n}$ by
$\pl \Phi_e = \sum \{ \Phi_{Q'}: e \in Q' \}$.

If $\pl \Phi = J$ we call $\Phi$ 
an integral decomposition of $J$.
\end{defn}

The following result of
Graver and Jurkat \cite{GJ} and Wilson \cite{W4}
shows that the necessary divisibility conditions on $J$
are sufficient for an integral decomposition $\Phi$,
that is, an assignment of integer weights to the $q$-cliques
in $K^r_n$ such that the total weight of cliques
on any edge $e$ is $J_e$. 

\begin{lemma} \label{GJ/W} \cite{GJ,W4}
Suppose $n \ge q+r$ and 
$J \in \mb{Z}^{K^r_n}$ is $K^r_q$-divisible.

Then there is $\Phi \in \mb{Z}^{K^r_q(K^r_n)}$ 
such that $\pl \Phi = J$. 
\end{lemma}

Now we will introduce the tools
of the proof of Lemma \ref{GJ/W}.

\begin{defn} \label{oct}
The $j$-octahedron $O_j$ is 
the complete $j$-partite $j$-graph with 
parts $\{i_0, i_1\}$ for $i \in [j]$. 
We denote its edges by $\{e_x: x \in \{0,1\}^j\}$,
where $e_x = \{ i_{x_i} : i \in [j]\}$. 
We define the sign of $e_x$ and $x$ 
by $s(e_x)=s(x)=(-1)^{\sum x}$.
\end{defn}

We view a copy $\phi(O_j)$ of $O_j$ in $K^j_n$
as a set of signed edges, 
or as a vector in $\mb{Z}^{K^j_n}$, where
each $\phi(O_j)_{\phi(e_x)} = s(e_x)$
and $\phi(O_j)_e$ is $0$ otherwise.

Given $S \sub \mb{Z}^d$, the integer span of $S$ is 
$\bgen{S} = \{ \sum_{x \in S} \Phi_x x : \Phi \in \mb{Z}^S \}$.
The next definition and lemma characterise
the integer span of octahedra.

\begin{defn} \label{null}
We say $J \in \mb{Z}^{K^j_n}$ is null
if $\sum \{ J_e : f \sub e \} = 0$ 
for all $f \in \tbinom{[n]}{j-1}$.
Note that any $j$-octahedron is null.
Let $\mc{N}_j$ be the set of null $J \in \mb{Z}^{K^j_n}$.
Let $\mc{O}_j$ be the set of all $j$-octahedra in $K^j_n$.
\end{defn}

\begin{lemma} \label{GJ/W:oct} \cite{GJ,W4}
$\bgen{\mc{O}_j} = \mc{N}_j$.
\end{lemma}

\nib{Remarks.} $ $
\begin{enumerate}
\item If $n<2j$ then $\mc{O}_j=\es$ and 
there are no non-trivial null $J \in \mb{Z}^{K^j_n}$.
\item In \cite{GJ} it is shown that one can even select a subset
of the octahedra that forms an integer basis of $\mc{N}_j$
(we mention this for the sake of interest,
but we do not use it in this paper).
\end{enumerate}

\medskip

Next we give a construction that
implements octahedra using $q$-cliques.
Suppose $\phi(O_j) \in \mc{O}_j$ and 
$Y \in \tbinom{[n] \sm V(\phi(O_j))}{q-j}$.
Define $\phi(O_j)*Y = \sum_{e \in \phi(O_j)}
s(e) Q^e \in \mb{Z}^{K^r_q(K^r_n)}$
where each $V(Q^e)=e \cup Y$. 
To put this in words, $\phi(O_j)*Y$ is obtained
from the octahedron $\phi(O_j)$ by extending
each of its signed edges to a $q$-clique with the same sign,
where each extension adds the same $(q-r)$-set $Y$.

\begin{lemma} \label{lem:makeoct}
$\pl_j (\phi(O_j)*Y) = \phi(O_j)$.
\end{lemma}

\nib{Proof.}
Every $e \in \phi(O_j)$ appears in a unique 
$q$-clique of $\phi(O_j)*Y$ with sign $s(e)$.
Any other $e \in K^j_n$ appears in $q$-cliques of 
$\phi(O_j)*Y$ the same number of times with each sign,
so does not contribute to $\pl_j (\phi(O_j)*Y)$. \qed

\begin{rem} \label{rem:div}
The $K^r_q$-divisibility constants will appear
when we use the above construction for $K^r_q$-decompositions:
if $\Phi \in \mb{Z}^{K^r_q(K^r_n)}$ then
$\pl_j \pl_r \Phi = \tbinom{q-j}{r-j} \pl_j \Phi$
for any $0 \le j \le r$.
\end{rem}

We conclude this subsection with a proof of Lemma \ref{GJ/W}
(which we do not use, but we include it for expository
purposes, as it illustrates some ideas of 
the proof of Lemma \ref{bddint}).
The idea of the proof is to modify $J$ by repeatedly
subtracting $q$-cliques so that it becomes `more null',
until it becomes zero. Here, and throughout the section,
we note that if $J$ is $K^r_q$-divisible and
$\Phi \in \mb{Z}^{K^r_q(K^r_n)}$
then $J - \pl \Phi$ is $K^r_q$-divisible. 
In particular, in the proof below each $J_j$ is $K^r_q$-divisible. 

We say $J \in \mb{Z}^{K^r_n}$ 
is $j$-null if $\pl_j J = 0$;
thus $J$ is $(r-1)$-null $\Lra$ $J$ is null,
$J$ is $0$-null $\Lra$ $\sum_e J_e = 0$,
and $J$ is $r$-null $\Lra$ $J=0$.
For clarity of notation in the proof,
we reserve $\pl$ for the shadow operators $\pl_j$ on $\mb{Z}^{K^j_q(K^j_n)}$
and include the subscript for the shadow operators $\pl_j$ on $\mb{Z}^{K^r_n}$.

\medskip

\nib{Proof of Lemma \ref{GJ/W}.}
Suppose $n \ge q+r$ and 
$J \in \mb{Z}^{K^r_n}$ is $K^r_q$-divisible.
We will define $\Phi_0,\dots,\Phi_r \in \mb{Z}^{K^r_q(K^r_n)}$
and $J_0 = J-\pl \Phi_0$, $J_j = J_{j-1}-\pl \Phi_j$ for $j \in [r]$,
proving by induction on $j \ge 0$ that each $J_j$ is $j$-null.
This will prove the lemma, as then $J_r=0$, so 
$\Phi = \sum_{j=0}^r \Phi_j$ satisfies $\pl \Phi = J$.

We start with $\Phi_0 = Q^{-1} \sum_e J_e \{Q^0\}$
for any fixed $Q^0 \in K^r_q(K^r_n)$, which is the vector in
$\mb{Z}^{K^r_q(K^r_n)}$ with $Q^{-1} \sum_e J_e$
in coordinate $Q^0$ and zero otherwise,
noting that $\sum_e J_e$ is divisible by $Q$, 
as $J$ is $K^r_q$-divisible.
Thus $\pl_0 J_0 = \sum_e J_e - Q \sum_\phi (\Phi_0)_\phi = 0$
gives the base case $j=0$ of the induction.

Now suppose inductively we have constructed $J_{j-1}$ for some $j \in [r]$.
Let $J^* = \tbinom{q-j}{r-j}^{-1} \pl_j J_{j-1}$. 
Then $J^* \in \mb{Z}^{K^j_n}$ as $J_{j-1}$ is $K^r_q$-divisible,
and $J^*$ is null (that is, $(j-1)$-null) as $J_{j-1}$ is $(j-1)$-null;
both deductions in this sentence use Remark \ref{rem:div}.

By Lemma \ref{GJ/W:oct} we have
$J^* \in \bgen{\mc{O}_j}$,
so there is $\Psi^j \in \mb{Z}^{\mc{O}_j}$
with $J^* = \sum_{X \in \mc{O}_j} \Psi^j_X X$.
Let $\Phi_j = \sum_{X \in \mc{O}_j} \Psi^j_X (X*Y_X) \in \mb{Z}^{K^j_q(K^j_n)}$,
choosing each $Y_X \in \tbinom{[n] \sm V(X)}{q-j}$ arbitrarily.
Then $\pl \Phi_j = J^*$ by Lemma \ref{lem:makeoct},
so $\pl_j \pl \Phi_j = \tbinom{q-j}{r-j} \pl \Phi_j =  \tbinom{q-j}{r-j} J^*= \pl_j J_{j-1}$
by Remark \ref{rem:div}.

Thus $J_j = J_{j-1} - \pl \Phi_j$ is $j$-null,
proving the induction step. The lemma follows. \qed

\subsection{Bounded generation I: octahedra}  \label{sub:bddgen}

As discussed in the introduction of this section,
we will show in this subsection that
the octahedral sets $\mc{O}_j$ can be replaced
by `thin' subsets $\mc{O}'_j$ that still generate the null spaces $\mc{N}_j$.

First we require some more notation.
We define a partial product $*$ on 
$\mb{Z}^{\{0,1\}^n}$ as follows.
If $v,v' \in \mb{Z}^{\{0,1\}^n}$ with $e \cap e' = \es$
whenever $v_e v_{e'} \ne 0$ then
$v*v' = \sum_{e,e'} v_e v_{e'} \{e \cup e'\}$;
otherwise $v*v'$ is undefined.
Note that 
\begin{enumerate}
\item $(u+v)*v'=u*v'+v*v'$
if both sides are defined,
\item any octahedron can be 
expressed as a product of $1$-octahedra: 
\[\phi(O_j) = *_{i=1}^j 
( \{\phi(i_0)\} - \{\phi(i_1)\} ).\]
\end{enumerate}

Next we introduce some more notation
for specifying octahedra.

\begin{defn} $ $
\begin{enumerate}
\item We define addition cyclically on $[n]$:
 $x+y$ is $x+y$ or $x+y-n$,
whichever is in $[n]$.
\item Suppose $f, \aA \in [n]^j$.
We define a copy $\phi^\aA_f(O_j)$ of $O_j$
by $\phi^\aA_f(i_0) = f_i$ and 
$\phi^\aA_f(i_1) = f_i+\aA_i$,
if all such vertices are distinct,
otherwise $\phi^\aA_f$ is undefined.
\item We say that $\phi^\aA_f(O_j)$
is thin if $\aA \in [2j]^j$.
\item Let $\mc{O}'_j$ be the set 
of all thin $j$-octahedra.
\end{enumerate}
\end{defn}

Note that any $j$-octahedron in $K^j_n$
can be written (in several ways)
in the form $\phi^\aA_f(O_j)$.

The key lemma of this subsection is that thin octahedra span all octahedra.

\begin{lemma} \label{octj}
$\sgen{\mc{O}'_j} = \bgen{\mc{O}_j}$.
\end{lemma}

\nib{Proof.}
We need to show that any $\phi^\aA_f(O_j)$
is in the integer span of $\mc{O}'_j$.
Say that $\phi^\aA_f(O_j)$ is $j'$-thin
if $\aA_i \in [2j]$ for $i \le j'$.
We show by induction on $j'=j,j-1,\dots,0$
that any $j'$-thin $j$-octahedron 
is in $\sgen{\mc{O}'_j}$.
This will prove the lemma, 
as any octahedron is $0$-thin.

For $j'=j$ note that any $j$-thin
octahedron is thin, so in $\mc{O}'_j$.
For the induction step, suppose $j' \in [j]$,
that $\phi^\aA_f(O_j)$ is $(j'-1)$-thin
and any $j'$-thin $j$-octahedron 
is in $\sgen{\mc{O}'_j}$.
Consider $\phi^{\aA'}_f(O_j) \in 
\phi^\aA_f(O_j) + \sgen{\mc{O}'_j}$
with $\aA'_i=\aA_i$ for all $i \ne j'$
and minimal $\aA'_{j'}>0$.
We claim that $\aA'_{j'} \in [2j]$,
so $\phi^{\aA'}_f(O_j)$ is $j'$-thin.
The induction step clearly follows,
so it remains to prove the claim.

Suppose for contradiction that $\aA'_{j'}>2j$.
Fix $\bB \in [2j]$ such that
$f_{j'}+\aA'_{j'}-\bB 
\notin V(\phi^{\aA'}_f(O_j))$.

Write $\phi^{\aA'}_f(O_j) 
= \phi_1(O_1)*\phi_2(O_{j-1})$,
where $\phi_1(O_1) 
= \{ f_{j'} \} - \{ f_{j'}+\aA'_{j'} \}$.

Let $\psi(O_1) = \{ f_{j'}+\aA'_{j'}-\bB \} 
- \{ f_{j'}+\aA'_{j'} \}$.

Note that $\psi(O_1)*\phi_2(O_{j-1})$
is $j'$-thin, so in $\sgen{\mc{O}'_j}$
by induction hypothesis. But then
$\phi^{\aA'}_f(O_j) - \psi(O_1)*\phi_2(O_{j-1})
= (\phi_1(O_1) - \psi(O_1))*\phi_2(O_{j-1})
= \phi^{\aA' - \bB e_{j'}}_f(O_j)$
contradicts minimality of $\aA'_{j'}$.
This proves the lemma. \qed

\subsection{Bounded generation II: cliques}  \label{sub:bddgenII}

In this subsection we use the results from the previous subsection
on bounded generation via octahedra to obtain analogous results
for $q$-cliques in $K^r_n$. It will be convenient and sufficient for our purposes 
to relax the pointwise bounds of the previous subsection to boundedness;
thus, instead of requiring that every edge of $K^r_n$ 
is covered by $O(1)$ elements of the generating set,
we will only require boundedness with respect to $(r-1)$-tuples.
A straightforward application of Lemma \ref{octj} above
gives $O(1)$-boundedness in Lemma \ref{bddgencomplete} below; 
however, our application will require $o(1)$-boundedness,
for which we will need a more complicated inductive construction,
to be explained below.

We start by constructing an $O(1)$-bounded generating 
subset of the $q$-cliques in $K^r_n$;
we also include a weak pointwise bound  for covering edges $e \in K^r_n$,
as this will be helpful when we apply the lemma for our more general construction below.

\begin{lemma} \label{bddgencomplete} 
For $n>n_0(q)$ sufficiently large
there is $S \sub K^r_q(K^r_n)$
with\footnote{Here we consider integer spans
of subsets of $\mb{Z}^{K^r_n}$,
identifying each clique $Q' \in K^r_q(K^r_n)$
with the vector $Q' \in \mb{Z}^{K^r_n}$
defined by $Q'_e = 1_{e \in Q'}$.} 
 $\bgen{S}=\bgen{K^r_q(K^r_n)}$ 
such that $\pl S$ is $(8r)^r q!$-bounded
and $|\pl S_e| < n^{0.01}$ for all $e \in K^r_n$.
\end{lemma}

\nib{Proof.}
Recall that $\mc{O}'_j$ denotes 
the set of all thin $j$-octahedra. 
For each $j \in [r]$ and $X \in \mc{O}'_j$ 
we choose independent uniformly random
$Y_X \in \tbinom{[n] \sm V(X)}{q-j}$
and add $\{e \cup Y_X: e \in X\}$ to $S$.
The proof that $\bgen{S}=\bgen{K^r_q(K^r_n)}$
is the same as that of Lemma \ref{GJ/W},
replacing $\mc{O}_j$ by $\mc{O}'_j$.

To show boundedness, we claim that 
$\mb{E} \pl S_e < 0.9 (8r)^r q!$ for all $e \in K^r_n$. 
To see this, note that for each $0 \le r' \le j \le r$
there are fewer than 
$2^{j-r'} (2j)^j \tbinom{r}{r'} n^{r'}$ choices of $(X,e')$
with $e' \in X = \phi^\aA_f(O_j) \in \mc{O}'_j$ such that
$|e' \sm e|=r'$ and $V(X) \cap e = e' \cap e$.
For each such $(X,e')$ we have a contribution 
of $2^{r'}$ to $\pl S_e$ with probability
$\mb{P}(e \sm e' \sub Y_X) 
\le (1+O(n^{-1})) (q-j)! n^{-|e \sm e'|}$,
where $|e \sm e'|=r-j+r'$.
Thus $\mb{E} \pl S_e < 
\sum_{r',j} 2^{j-r'} (2j)^j \tbinom{r}{r'} n^{r'}
\cdot 2^{r'} \cdot  (q-j)! n^{j-r-r'} 
= (8r)^r (q-r)! + O(n^{-1}) < 0.9 (8r)^r q!$,
as claimed. By Bernstein's inequality, we deduce
w.h.p.\ $|\pl S_e| < n^{0.01}$ for all $e \in K^r_n$.
We also deduce w.h.p.\ $|\pl S(f)| < (8r)^r q! n$
for all $f \in \tbinom{[n]}{r-1}$,
so $\pl S$ is $(8r)^r q!$-bounded. \qed

\medskip

The main result of this subsection is the following version of the previous
lemma relative to a bounded subgraph $L$ of $K^r_n$. 
It shows that if $L$ is $o(1)$-bounded then there is an $o(1)$-bounded generating
set of cliques for the $K^r_q$-divisible vectors supported in $L$.

\begin{lemma} \label{bddgen} 
Let\footnote{The notational use of $L$ and $S$ in this section
is unrelated to our global notation for the leave and the spill.
} 
$L \sub K^r_n$ be $\nu$-bounded,
where $n^{-1} \le \nu \le \nu_0(q,r)$,
with $\nu_0(q,r)$ sufficiently small,
and $n \ge n_0(q,r)$ sufficiently large.
Then there is $S \sub K^r_q(K^r_n)$ with%
\footnote{Here we identify $\mb{Z}^L$ 
with the set of $v \in \mb{Z}^{K^r_n}$
supported in $L$.} 
$\bgen{K^r_q(K^r_n)} \cap \mb{Z}^L \sub \bgen{S}$
such that $\pl S$ is $2q! (8r)^r \nu^{1/b}$-bounded,
where $b = 2^{3^{r+q}}$.
\end{lemma}

Lemma \ref{bddgen} is immediate from the following lemma, 
in which we strengthen the conclusion so that 
it is amenable to proof by induction.
The lemma shows that there is a probability distribution on generating sets,
which w.h.p.~satisfies the conclusion of Lemma \ref{bddgen},
and has further properties useful for the inductive proof,
namely a weak pointwise bound for covering edges $e \in K^r_n$
and bounds on the probabilities of using specific cliques $Q' \in K^r_q(K^r_n)$.

\begin{lemma} \label{bddgen:IH} 
Let $L \sub K^r_n$ be $\nu$-bounded,
where $n^{-1} \le \nu \le \nu_0(q,r)$,
with $\nu_0(q,r)$ sufficiently small,
and $n \ge n_0(q,r)$ sufficiently large.
Let $b_r = 2^{3^{r+q}}$ and $a_r=0.9 b_r$.
Then there is a probability distribution
on subsets $S$ of $K^r_q(K^r_n)$ such that,
\begin{enumerate}
\item $\mb{P}(Q' \in S) < 2q! n^{r-q}$
for all $Q' \in K^r_q(K^r_n)$, 
\item $\mb{P}(Q' \in S) < \nu^{1/a_r} n^{r-q}$
for all $Q' \in K^r_q(K^r_n)$ with $Q' \cap L = \es$, and 
\item w.h.p.\ $\pl S$ is $2q! (8r)^r \nu^{1/b_r}$-bounded,
$|\pl S_e| < n^{0.1(1-1/2r)}$ for all $e \in K^r_n$,
and $\bgen{K^r_q(K^r_n)} \cap \mb{Z}^L \sub \bgen{S}$.
\end{enumerate}
\end{lemma}

Before giving the formal proof, we illustrate the idea of the proof
by sketching its application to the case of triangles, that is $(q,r)=(3,2)$;
this case can be handled more simply by the methods in \cite{K2},
but here we use it to illustrate the general argument for the sake of exposition.

Suppose then that $L \sub K_n$ is $\nu$-bounded, that is, has maximum degree $\le \nu n$.
We will construct $S = S^0 \cup S^1 \cup S^2$ as a random set of triangles in $K_n$.
The plan is that any triangle-divisible $J$ supported in $L$ will be expressed as integral
combination of triangles in $S$ via the following process:
\begin{enumerate}[label=(\arabic*)] \setcounter{enumi}{-1}
\item fix a random set $M$ of $m = \nu^c n$ vertices 
for some small $c$ (say $M=[m]$ by relabelling),
\item use triangles in $S^0$ to eliminate the support of $J^0=J$ 
on edges disjoint from $[m]$, thus reducing to some triangle-divisible $J^1$
supported on edges that intersect $[m]$,
\item use triangles in $S^1$ to eliminate the support of $J^1$ 
on edges with exactly one vertex in $[m]$,
thus reducing to some triangle-divisible $J^2$ supported in $K_m$,
\item use triangles in $S^2$ to express $J^2$.
\end{enumerate}

A suitable construction of $S^2$ is already provided by Lemma \ref{bddgencomplete} with $m$ in place of $n$,
which gives the required boundedness properties of $\pl S^2$ (as $m/n$ is small).
Applying a random permutation of $[m]$ ensures that any fixed triangle of $K_n$
appears in $S^2$ with probability $|S^2| \tbinom{n}{3}^{-1} = O(m^2/n^3)$.

We construct $S^0$ as a set of triangles $Q^e$, where for each edge $e \in L$
independently we choose a triangle $Q^e$ of $K_n$ uniformly at random
subject to containing $e$ and some third vertex in $[m]$.
Clearly $S^0$ can be used to eliminate the support of $J$
on edges disjoint from $[m]$, as required for (1) above.
Any fixed triangle $Q$ is chosen with probability at most $3 \cdot m/n \cdot 1/m = 3/n$,
as for each edge $e$ of $Q$, if $e \in L$ then we may choose $Q^e=Q$ if the third vertex
of $Q$ is chosen for $M$ and then chosen as the third vertex of $Q^e$.
This will be the only place in the construction where some fixed triangle has such a large
appearance probability (which is allowed for triangles containing an edge of $L$).

To see the required boundedness conditions on $\pl S^0$, we note that
any edge $xy \in K_n$ with $x \in [m]$ and $y \in e$ is covered by $Q^e$ with probability $1/m$.
Any vertex $v$ is thus covered $|L(v)| \le \nu n$ times in the role of $y$,
and in the role of $x$ a pseudobinomial number
with mean at most $|L|/m \le 2\nu n^2/m = 2\nu^{1-c} n$.
Similarly, the number of times any given edge is covered is pseudobinomial
with mean at most $\nu n/m = \nu^{1-c}$,  so the required boundedness conditions 
on $\pl S^0$ hold w.h.p~by Chernoff bounds. 

We let $L^1$ be the random graph of edges in $Q^e \sm \{e\}$ for some $e \in L$,
which satisfies $\mb{P}(e' \in L^1) \le \nu^{1-c}$ for any edge $e'$ 
and is w.h.p.~$2\nu^{1-c}$-bounded. 
For $v \in [n] \sm [m]$ we write $L^v := L^1(v) \sub [m]$.

Now we come to the construction of $S^1$, using the inductive hypothesis 
for $(q',r')=(2,1)$ in the following form:~for any vertex $v \in [n] \sm [m]$, 
there is a probability distribution on subgraphs $S_v$ of $K_m$ such that 
\begin{enumerate}
\item any edge appears in $S_v$ with probability at most $4/m$,
\item any edge disjoint from $L^v$ appears with probability at most $\nu^{c'}/m$ for some fixed $c'>0$,
and \item w.h.p. $|S_v| < \nu^{c'} m$, the maximum degree of $S_v$ is at most $m^{0.1/2}$,
and any $J_v \in \mb{Z}^{L^v}$ with $\sum_x (J_v)_x$ even
can be expressed as an integer combination of edges in $S_v$.
\end{enumerate}
The construction of $S_v$ is similar to that of $S^0 \cup S^2$ above,
but simpler as we are choosing random edges rather than triangles:
we fix a random small set $M'$, include for each $x \in L^v$ a random edge 
from $x$ to $M'$, and a suitably sparse random set of pairs in $M'$
spanning all vectors supported in $M'$ with even sum, 
with the latter provided by Lemma \ref{bddgencomplete} (or a simple explicit construction).

We construct $S^1$ by sampling $S_v$ as above independently for each $v \in [n] \sm [m]$
and taking all triangles of the form $vxy$ with $xy \in S_v$.
Thus $S^1$ can be used to reduce the support to $K_m$, as required for (2) above:\
given some triangle-divisible $J^1$ supported on edges that intersect $[m]$,
for each $v \in [n] \sm [m]$, by property (iii) of $S_v$ we can write 
$J^1(v) = \sum_{xy \in S_v} \Phi_{xy} (1_x + 1_y)$ for some $\Phi \in \mb{Z}^{S_v}$,
then eliminate the support of $J^1$ on edges incident to $v$
using $\Phi^v \in \mb{Z}^{S^1}$ defined by $\Phi^v_{vxy} = \Phi_{xy}$.

We verify the other required properties of $S^1$ 
using properties (i) and (ii) of the distribution on $S_v$.
For any vertex $v \in [n] \sm [m]$ and edge $xy \in K_n$ we estimate $\mb{P}(xy \in S_v)$
by considering whether $xy$ intersects the random set $L^v:=L^1(v) \sub [m]$.

We have $\mb{P}(xy \cap L^v \ne \es) \le \mb{P}(xv \in L^1) +  \mb{P}(yv \in L^1) \le 2\nu^{1-c}$
and $\mb{P}(xy \in S_v \mid xy \cap L^v \ne \es) \le 4/m = 4\nu^{-c}/n$ by (i),
$\mb{P}(xy \in S_v \mid xy \cap L^v = \es) \le \nu^{c'}/m =  \nu^{c'-c}/n$ by (ii),
so $\mb{P}(xy \in S_v) \le (8\nu^{1-2c}+ \nu^{c'-c})/n$.  
Taking $c=\min\{c'/2,1/4\}$, we thus have a suitable bound on the probability
that any fixed triangle is chosen for $S^1$.
We also deduce for any fixed edge $e \in K_m$ that the expected number 
of triangles in $S^1$ covering $e$ is at most $2(8\nu^{1-2c}+ \nu^{c'-c})$,
so the required boundedness conditions on $\pl S^1$ hold w.h.p~by Chernoff bounds. 
(The contribution from edges of $\pl S^1$ containing some $v \in [n] \sm [m]$
is non-random given $S_v$ and bounded by property (iii) of $S_v$:~at most
$|S_v| < \nu^{c'} m$ edges cover $v$ and any $uv$ with $u \in [m]$
has multiplicity at most $m^{0.1/2} < n^{0.1 \cdot 3/4}$.)

The final argument above was the motivation for strengthening the inductive hypothesis 
to include properties of a random $S_v$, rather than just the existence of $S_v$.
Now we give the general proof.

\medskip

\nib{Proof of Lemma \ref{bddgen:IH}.}
We use induction on $r \ge 1$;
the induction hypothesis is the statement of the lemma,
which we assume for any $(q',r')$ with $1 \le r'<r$.
We will take $S = \cup_{i=0}^r S^i$
for some $S^i \sub K^r_q(K^r_n)$.
We start by giving the constructions 
of $S^r$ and $S^0$. These do not use
the induction hypothesis, and in the base
case $r=1$ we will take $S = S^0 \cup S^r = S^0 \cup S^1$.

\medskip

\nim{Step 1:\ defining $S^r$, $S^0$, $L^1$; the base case.} 
(See $S^2$ and $S^0$ in the sketch for triangles.)

\medskip

Let $m = \nu^{3/b_r} n$ and $S^r \sub K^r_q(K^r_m)$
be given by Lemma \ref{bddgencomplete}:
we have $\bgen{S^r}=\bgen{K^r_q(K^r_m)}$ 
and $\pl S^r$ is $(8r)^r q!$-bounded in $[m]$,
and $|\pl S^r_e| < m^{0.01}$ for all $e \in K^r_m$.
Let $\pi:[m] \to [n]$ be a uniformly random injection
and let $\pi(S^r) = \{ \pi(Q'): Q' \in S^r\}$.
Then $\bgen{\pi(S^r)} = \bgen{K^r_q(\tbinom{\pi([m])}{r})}$
and $\mb{P}(Q' \in \pi(S^r)) = |S^r|\tbinom{n}{q}^{-1} 
< (8r)^r q! \tbinom{m}{r} \tbinom{n}{q}^{-1} 
< (8r)^r q!^2 \nu^{3r/b_r} n^{r-q}
< \nu^{2/a_r} n^{r-q}$ for all $Q' \in K^r_q(K^r_n)$. 
For convenient notation we relabel so that 
$\pi$ is the identity embedding of $[m]$ in $[n]$.

We let $S^0 = \{ Q^e : e \in L \}$,
where for each $e \in L$ independently 
we choose $Q^e \in K^r_q(K^r_n)$ 
uniformly at random subject to $e \in Q^e$
and $V(Q^e) \sm e \sub [m]$.
We also let $S^{0-} = \{ Q^e \sm \{e\} : e \in L \}$
and $L^1 = \bigcup_{e \in L} Q^e \sm \{e\}$.

We claim for any $e' \in K^r_n$ that
$\mb{E} \pl S^{0-}_{e'} < q! \nu (n/m)^r$.
To see this, we fix any $e \in L$
with $e \ne e'$ and estimate $\mb{P}(e' \in Q^e)$.
We can assume $e' \sm e \sub [m]$
and $r' := |e' \sm e| \le q-r$,
otherwise the probability is $0$.
There are at least $\tbinom{m-r}{q-r}$ choices for $Q^e$, 
of which at most $m^{q-r-r'}$ contain $e'$, 
so $\mb{P}(e' \in Q^e) < (1+O(m^{-1})) (q-r)! m^{-r'}$.
As $L$ is $\nu$-bounded, there are at most
$\nu \tbinom{r}{r'} n^{r'} + O(n^{r'-1})$ such choices of $e$,
so summing over $r'$ gives the claim.

We deduce that 
$\mb{P}(e' \in L^1) < q! \nu (n/m)^r < \sqrt{\nu}$ (say), 
w.h.p.\ $|\pl S^0_{e'}| < m^{0.01}$ for all $e' \in K^r_n$
and w.h.p.\ $\pl S^0$ is $\sqrt{\nu}$-bounded,
and so $L^1$ is $\sqrt{\nu}$-bounded. 

Furthermore, for any $Q' \in K^r_q(K^r_n)$,
there are at most $\tbinom{q}{r}$ choices of $e \in Q' \cap L$,
and for each, the probability of choosing $\pi$ 
such that $V(Q') \sm e \sub [m]$ is 
$\tbinom{n-q+r}{m-q+r} \tbinom{n}{m}^{-1}$
and then $\mb{P}(Q^e=Q' \mid \pi) \le \tbinom{m-r}{q-r}^{-1}$,
so $\mb{P}(Q' \in S^0) \le 
\tbinom{q}{r} \tbinom{n-r}{m-q+r} \tbinom{n}{m}^{-1}
\tbinom{m-r}{q-r}^{-1} < q! n^{r-q}$.
 
In the base case $r=1$ of the lemma,
we now claim that taking $S = S^0 \cup S^1$
completes the proof. It remains to show that
$\bgen{K^1_q(K^1_n)} \cap \mb{Z}^L \sub \bgen{S}$.
To see this, we consider any
$J \in \bgen{K^1_q(K^1_n)} \cap \mb{Z}^L$.
We define $J' = J - \pl_1 \Phi'$, where
for each $e \in L$ we add $J^0_e \{Q^e\}$ to $\Phi'$,
that is, $\Phi'_Q = \sum \{ J^0_e: Q^e=Q\}$ for each $Q$;
this cancels the coefficients of all such $e$,
and all new signed elements $e'$ of $J'$
are contained in $[m]$.
Thus we obtain $J' \in \bgen{K^1_q(K^1_m)}
= \bgen{S^1} \sub \bgen{S}$, as required.

\medskip

\nim{Step 2:\ defining $S^i$ and $L^{i+1}$ for $i \in [r-1]$,
with $S^i$ obtained from the induction hypothesis for
$L^f := L^i(f)[[m]]$ for each $f \in \tbinom{[n] \sm [m]}{r-i}$.} 
(See definition of $S^1$ in the sketch for triangles.)

\medskip

Now suppose $r>1$. We construct $S^i$ sequentially
for $1 \le i \le r-1$ using the induction hypothesis.
Let $\nu_0=\nu$ and $\nu_{i+1} = \nu_i^{1/b_i}$
for $0 \le i \le r-1$. 
At the start of round $i$ we will have some random
$L^i \sub K^r_n$ that is $\nu_i$-bounded,
such that all $\mb{P}(e \in L^i) < 0.1 \nu_i$;
this holds for $i=1$ as $\sqrt{\nu}<0.1\nu_1$. 
Note that each
$\nu_j = \nu^{1/\prod_{i=0}^{j-1} b_i} < \nu^{1/\sqrt{b_j}}$, 
as $\sum_{i=0}^{j-1} 3^{i+q} < \tfrac{1}{2} 3^{j+q}$.

For each $f \sub [n] \sm [m]$ with $|f|=r-i$
we let $L^f := L^i(f)[[m]]$ be the restriction
of the neighbourhood $L^i(f)$ to $[m]$, 
and note that $L^f \sub K^i_m$ is $\nu'_i$-bounded,
where $\nu'_i = \nu_i n/m \ge m^{-1}$,
as $\nu_i \ge \nu \ge n^{-1}$.
By the induction hypothesis we can choose
(independently for each $f$)
a random $R^f \sub K^i_{q-r+i}(K^i_m)$ such that
\begin{enumerate}
\item $\mb{P}(X \in R^f) < 2(q-r+i)! m^{r-q}$
for all $X \in K^i_{q-r+i}(K^i_m)$, 
\item $\mb{P}(X \in R^f) < (\nu'_i)^{1/a_i} m^{r-q}$ 
for all $X \in K^i_{q-r+i}(K^i_m)$
with $X \cap L^f = \es$, and 
\item w.h.p.\ $\pl R^f$ is 
$2(q-r+i)! (8i)^i (\nu'_i)^{1/b_i}$-bounded in $[m]$,
$|\pl R^f_e| < m^{0.1(1-1/2i)}$ for all $e \in K^i_m$, and
$\bgen{K^i_{q-r+i}(K^i_m)} \cap \mb{Z}^{L^f} \sub \bgen{R^f}$.
\end{enumerate}
We obtain $S^f = R^f*f = \{ X*f: X \in R^f \}
\sub K^r_q(K^r_n)$ by adding $f$ 
to the vertex-set of each clique in $R^f$.
We let $S^i$ be the union of all such $S^f$ and let 
$L^{i+1} = L^i \cup \{e: \pl S^i_e > 0 \}$.

This defines round $i$ of the construction;
the boundedness of $L^{i+1}$ required at the start of round $i+1$ above 
will be established as part of the next step.

\medskip

\nim{Step 3:\ boundedness of $\pl S^i$,  $L^{i+1}$ 
and the probability of choosing any fixed $q$-clique. }
(See analysis of $S^1$ in the sketch for triangles.)

\medskip

We require the following claim, which we prove by induction on $i$.
(This secondary induction argument is internal to each step 
of the primary induction on $r$ used to prove the lemma;
the boundedness of $L^i$ is also proved by the secondary induction,
but to lighten the exposition we will not state it formally as part of the claim.)
 
 \medskip

\nim{Claim.} $\mb{E} \pl S^i_e \le 0.1 \nu_{i+1}$
for any $e \in K^r_n$ and $i \in [r-1]$.

\medskip

To prove this, note first that
$\mb{P}(e' \in L^i) \le \mb{E} \pl S^{i-1}_{e'}$
and $\mb{E} \pl S^{i-1}_{e'} \le 0.1 \nu_i$,
by the claim induction hypothesis for $i>1$, 
or  in the claim base case $i=1$
by the  bound $\mb{P}(e' \in L^1) < \sqrt{\nu}$
obtained in Step 1.
Thus for any $f \sub [n] \sm [m]$ with $|f|=r-i$
and $X \in K^i_{q-r+i}(K^i_m)$ we have 
\begin{align*}
\mb{P}(X \in R^f)
& = \mb{P}(X \cap L^f = \es )
\mb{P}(X \in R^f \mid X \cap L^f = \es)
+ \mb{P}(X \cap L^f \ne \es )
\mb{P}(X \in R^f \mid X \cap L^f \ne \es) \\
& \le 1 \cdot  (\nu'_i)^{1/a_i} m^{r-q}
+ Q \cdot 0.1 \nu_i \cdot 2(q-r+i)! m^{r-q}
\le 2 (\nu'_i)^{1/a_i} m^{r-q},
\end{align*}
using $\mb{P}(e' \in L^f)=\mb{P}(e' \cup f \in L^i) \le 0.1\nu_i$
in a union bound over $e' \in X$.

Here we digress from the proof of the claim to deduce
the required bounds on the probability of choosing
any fixed $q$-clique $Q' \in K^r_q(K^r_n)$:\ we have
$\mb{P}(Q' \in S^i) 
\le 2(\nu'_i)^{1/a_i} m^{r-q} 
< 2 \nu_i^{1/a_i} (n/m)^q n^{r-q}
< 2 \nu_{i+1} \nu^{-3q/b_r} n^{r-q}
< \nu^{1/2\sqrt{b_r}} n^{r-q}$ (say). 
Summing over $i$ we deduce
$\mb{P}(Q' \in S) < 2q! n^{r-q}$
for all $Q' \in K^r_q(K^r_n)$, and
$\mb{P}(Q' \in S) < \nu^{1/a_r} n^{r-q}$
for all $Q' \in K^r_q(K^r_n)$ with $Q' \cap L = \es$.

Returning to the proof of the claim,
we now consider any $e \in K^r_n$ 
with $r' = |e \cap [m]| \ge i$.
There are at most 
$n^{r'-i}$ choices for an $(r-i)$-set $f$
with $e \sm [m] \sub f \sub [n] \sm [m]$,
and fewer than $m^{q-r+i-r'}$ 
choices for $X \in K^i_{q-r+i}(K^i_m)$ 
with $e \cap [m] \sub V(X)$.
Then for each such $f$ we have $\mb{E} \pl S^f_e 
\le m^{q-r+i-r'} \cdot 2 (\nu'_i)^{1/a_i} m^{r-q}
= 2 (\nu'_i)^{1/a_i} m^{i-r'}$,
so summing over $f$ gives
$\mb{E} \pl S^i_e
\le 2 (\nu'_i)^{1/a_i} (n/m)^{r'-i}
= 2 (\nu'_i)^{1/a_i} \nu^{3(i-r')/b_r} 
< 2 \nu_{i+1}^{1.1} \nu^{-3r/b_r} 
< 0.1 \nu_{i+1}$, where for $i=r-1$ 
we recall $\nu_r < \nu^{1/\sqrt{b_r}}$.
This proves the claim. 

\medskip

We deduce that
 $|\pl S^i_e|$ is $(m^{0.1(1-1/2i)},1)$-dominated,
as the maximum contribution from
each $f$ is at most $m^{0.1(1-1/2i)}$.
Thus w.h.p.\ $|\pl S^i_e| < m^{0.1(1-1/2r)}$
by Lemma \ref{dom}. 

We also claim that w.h.p.\ $\pl S^i$ is $\nu_{i+1}/2$-bounded.
To see this, we fix any 
$f' \in \tbinom{[n]}{r-1}$ with 
$|f' \sm [m]| \le r-i$, so
$r' := |f' \cap [m]| \ge i-1$,
and estimate $|(\pl S^i)(f')|$.
If $r' \ge i$ then by the above estimates
$|(\pl S^i)(f')|$ is $(m^{0.1},0.1 \nu_{i+1} n)$-dominated,
so w.h.p.\ $|(\pl S^i)(f')| < \nu_{i+1}n/2$.
On the other hand, if $r'=i-1$ then 
\[ |(\pl S^i)(f')| = |(\pl R^{f'})(f' \cap [m])| 
< 2q! (8i)^i (\nu'_i)^{1/b_i} m
= 2q! (8i)^i (\nu^{1/b_r})^{3-1/b_i} \nu_{i+1} n
< \nu_{i+1} n/2,\]
 as $\nu < \nu_0(q,r)$, as claimed.

We deduce that $L^{i+1}$ is $\nu_{i+1}$-bounded,
as required at the start of round $i+1$ in Step 3,
and also that $\pl S$ is $2q! (8r)^r \nu^{1/b_r}$-bounded,
as $\pl S^r$ is $(8r)^r q! m/n$-bounded
and $\sum_{i=1}^r \nu_i < 2\nu_r < \nu^{1/b_r}$.

\medskip

\nim{Step 4:\ verifying that $S$ generates
all $K^r_q$-divisible vectors supported in $L$.}
(See `the plan' in the sketch for triangles.)

\medskip

It remains to show that
$\bgen{K^r_q(K^r_n)} \cap \mb{Z}^L \sub \bgen{S}$.
To see this, we consider any
$J \in \bgen{K^r_q(K^r_n)} \cap \mb{Z}^L$.
We let $J^0=J$ and construct 
$J^i = J^{i-1} - \pl \Phi^i \in \mb{Z}^{L^i}$ 
for $i \in [r]$ where $\Phi^i \in \mb{Z}^S$ 
such that $J^i_e=0$ whenever $|e \cap [m]|<i$.
To define $J^1 = J^0 - \pl \Phi^1$, for each $e \in L$ 
we add $J^0_e \{Q^e\}$ to $\Phi^1$;
this cancels the coefficients of all such $e$,
and all new signed elements $e'$ of $J^1$ 
have $e' \cap [m] \ne \es$ and $e' \in L^1$.

Given $J^i$ with $0<i<r$, for each 
$f \sub [n] \sm [m]$ with $|f|=r-i$ we note that 
$J^i(f) \in \bgen{K^i_{q-r+i}(K^i_m)} \cap \mb{Z}^{L^f}
\sub \bgen{R^f}$, so $J^i(f) = \pl_i \Phi^f$
for some $\Phi^f \in \mb{Z}^{R^f}$.
We define $\Phi^{i+1} \in \mb{Z}^{S^i}$ 
as the sum over all such $f$ of
$\sum_{X \in R^f} \Phi^f_X \{X*f\}$.
Then $\pl_r \Phi^{i+1}_e = J^i_e$ for all
$e \in L^i$ with $|e \cap [m]|=i$,
so all such coefficients are cancelled
in $J^{i+1} = J^i - \pl \Phi^i$,
and all new signed elements $e'$ of $J^{i+1}$ 
have $|e' \cap [m]|>i$ and $e' \in L^{i+1}$.
Thus we obtain $J^r \in \bgen{K^r_q(K^r_m)}
= \bgen{S^r} \sub \bgen{S}$. \qed

\subsection{Bounded integral decomposition}  \label{sub:bddint}

The main result of this section is an analogue
of Lemma \ref{GJ/W} on integral decomposition
in which we also impose a boundedness condition on $\Phi$.
For $\Phi \in \mb{Z}^{K^r_q(K^r_n)}$ we define $\pl^\pm \Phi \in \mb{N}^{K^r_n}$ by 
\[ \pl^+ \Phi_e =  \sum 
\{ \Phi_{Q'}: e \in Q', \Phi_{Q'}>0 \}
\ \text{ and } \
\pl^- \Phi_e =  \sum 
\{ -\Phi_{Q'}: e \in Q', \Phi_{Q'}<0 \}.\]

\begin{lemma} \label{bddint} 
For any $K^r_q$-divisible $\tT$-bounded $J \in \mb{Z}^{K^r_n}$
where $n>n_0(q)$ is sufficiently large 
and $\tT>n^{-1/4Qb}$, where $b = 2^{3^{r+q}}$,
there is $\Phi \in \mb{Z}^{K^r_q(K^r_n)}$ 
such that $\pl \Phi = J$ and $\pl^\pm \Phi$ 
are $N^2 \tT$-bounded, where $N:=(2q)^q !$.
\end{lemma}

We will require several other lemmas
for the proof of Lemma \ref{bddint}.
Our first two lemmas will prove it in 
the `highly divisible' case of $J \in N\mb{Z}^{K^r_n}$,
using `robust local decodability' 
of the lattice of $K^r_q$-divisible vectors:
for any $e \in K^r_n$ there are many ways to write
$N\{e\} = \pl \Psi$ where $\Psi \in \mb{Z}^{K^r_q(K^r_n)}$
is `small'. We will use the bounded local generating set 
for a sparse random subgraph of $K^r_n$ 
obtained in the previous subsection to reduce the
general case of Lemma \ref{bddint} to the highly divisible case.

\begin{lemma} \label{decode} 
There is $\Psi^* \in \mb{Z}^{K^r_q(K^r_{r+q})}$
with $\pl \Psi^* =  N \{[r]\}$ and%
\footnote{Recall that for $v \in \mb{Z}^X$ 
we write $|v| = \sum_{x \in X} |v_x|$.} 
$|\Psi^*| < N^2$.
\end{lemma}

\nib{Proof.}
By Gottlieb's Theorem \cite{Go},
the inclusion matrix $M=M^r_q(q+r)$
(see Definition \ref{def:inclusion}) has full rank.
By Cramer's rule, every entry of $M^{-1}$ is rational
with absolute value and denominator both at most $N$;
indeed, as $M$ has entries in $\{0,1\}$, any submatrix
has determinant at most $\tbinom{q+r}{r}! \le N$.
Let $\Psi^* = N M^{-1} v$, where 
$v \in \mb{Z}^{K^r_{r+q}}$ 
with all $v_e=1_{e=[r]}$. \qed

\begin{lemma} \label{bddint:N}
Suppose $n>n_0(q)$ is large, $\tT>n^{-1/2}$ and
$J \in N\mb{Z}^{K^r_n}$ is $\tT$-bounded.
Then there is $\Phi \in \mb{Z}^{K^r_q(K^r_n)}$ such that 
$\pl \Phi = J$ and $\pl^\pm \Phi$ are $(2q)^q N \tT$-bounded.
\end{lemma}

\nib{Proof.}
For each signed element $e$ of $N^{-1}J$
we choose independent uniformly random
$\psi^e(K^r_{r+q}) \sub K^r_n$ with $\psi^e([r])=e$
and add $s(e) \psi^e(\Psi^*)$ to $\Phi$.
Then $\pl \Phi = J$. For the boundedness condition,
for any $e' \in K^r_n$ we estimate
$\GG_{e'} := \sum_{e \ne e'} \mb{P}(e' \in \psi^e(K^r_{r+q}) \sm \{e\})$.
As $J$ is $\tT$-bounded, for each $r' \in [r]$ there are fewer 
than $\tbinom{r}{r'} N^{-1}\tT n^{r'}$ signed elements $e$ of $N^{-1}J$
with $|e \sm e'|=r'$. For each such $e$
there are $\tbinom{n-r}{q}$ choices for $\psi^e(K^r_{r+q})$,
of which at most $n^{q-r'}$ contain $e'$, so
$\mb{P}(e' \in \psi^e(K^r_{r+q}) \sm \{e\}) 
< (1+O(n^{-1})) q! n^{-r'}$.
Summing over $r'$ we obtain $\GG_{e'} < 2^r q! N^{-1}\tT$.
Then by Chernoff bounds and Lemma \ref{decode}  
w.h.p.\ $\pl^\pm \Phi$ are $(2q)^q N \tT$-bounded. \qed

\medskip

The next lemma allows us to `flatten' 
any $J \in \mb{Z}^{K^r_n}$
without incurring any 
significant loss in boundedness.

\begin{lemma} \label{bddint:flat} 
For any $\tT$-bounded $J \in \mb{Z}^{K^r_n}$,
where $n>n_0(q)$ is large and $\tT>n^{-1/2}$,
there are $J' \in \mb{Z}^{K^r_n}$
and $\Phi \in \mb{Z}^{K^r_q(K^r_n)}$ 
such that $\pl \Phi = J-J'$,
all $|J'_e| < n^{0.1}$ and $J'$ and
$\pl^\pm \Phi$ are $q^q \tT$-bounded.
\end{lemma}

\nib{Proof.}
For each signed element $e$ of $J$
we add to $\Phi$ a uniformly random $Q^e$
with $e \in Q^e \in K^r_q(K^r_n)$,
where the sign of $Q^e$ in $\Phi$
is the same as that of $e$ in $J$.

For any $f \in \tbinom{[n]}{r-1}$ and $k \in [r]$
there are at most $\tbinom{r-1}{k-1} \tT n^k$ signed elements $e$ 
of $J$ with $|e \sm f|=k$. For each such $e$
there are $\tbinom{n-r}{q-r}$ choices of $Q^e$,
of which at most $n^{q-r-(k-1)}$ contain $f$,
so $\mb{P}(f \sub Q^e) < (1+O(n^{-1})) (q-r)! n^{-k+1}$. 
Then $\pl^\pm \Phi(f)$ are pseudobinomial
with mean at most $2^r (q-r)! \tT n$, so w.h.p.\ 
$J'$ and $\pl^\pm \Phi$ are $q^q \tT$-bounded.

It remains to bound $|J'_{e^*}|$ for any $e^* \in K^r_n$.
Each signed element $e'$ counted by $J_{e^*}$
is cancelled by $Q^{e'}$ in $\Phi$,
so any other nonzero contribution to $J'_{e^*}$
comes from $Q^e$ with $e \ne e^*$.
Similarly to above, for any $k \in [r]$
there are at most $\tbinom{r}{k} \tT n^k$ 
signed elements $e$ with $|e \sm e^*|=k$. 
For each such $e$ at most $n^{q-r-k}$ 
choices of $Q^e$ contain $e^*$,
so $\mb{P}(e^* \sub Q^e) < (1+O(n^{-1})) (q-r)! n^{-k}$. 
Thus $J'_{e^*}$ is pseudobinomial
with mean at most $2^r (q-r)! \tT$, 
so w.h.p.\ $|J'_{e^*}| < n^{0.1}$. 
\qed

\medskip

The next lemma will allow us to focus 
within a sparse random subgraph $L$.
The cost in boundedness is only a constant factor;
it is crucial that this is independent of the density $d(L)$.

\begin{lemma} \label{bddint:random} 
Suppose $J \in \mb{Z}^{K^r_n}$ is $\tT$-bounded,
where $n>n_0(q)$ is large.
Let $L \sub K^r_n$ be $(c,Q)$-typical 
and such that $J$ is $(\tT,Q)$-bounded w.r.t.\ $L$.
Then there is some $J' \in \mb{Z}^L \sub \mb{Z}^{K^r_n}$
and $\Phi \in \mb{Z}^{K^r_q(K^r_n)}$ 
such that $\pl \Phi = J-J'$
and $J'$ and $\pl^\pm \Phi$ are $q^q \tT$-bounded.
\end{lemma}

\nib{Proof.}
We define $\Phi$ by including 
for each signed element $e$ of $J$
a uniformly random $Q^e$ with 
$e \in Q^e \in K^r_q(K^r_n)$
and $Q^e\sm\{e\} \sub L$,
where the sign of $Q^e$ in $\Phi$
is the same as that of $e$ in $J$.
Then $J' := J-\pl \Phi \in \mb{Z}^L$.

We claim for any $e' \in L$ that
$E_{e'} := \sum_e \mb{P}(e' \sub Q^e) 
< 1.3 (q-r)! 2^r \tT d(L)^{-1}$.
To see this, first note that for any $k \in [r]$,
as $J$ is $(\tT,Q)$-bounded w.r.t.\ $L$
there are at most
$\tbinom{r}{k} d(L)^{\tbinom{k+r}{r}-2} \tT n^k$
signed elements $e$ of $J$ with $|e \sm e'|=k$
and $\tbinom{e' \cup e}{r} \sm \{e\} \sub L$.
For each such $e$, as $L$ is $(c,Q)$-typical,
there are at least $0.9 d(L)^{Q-1} \tbinom{n}{q-r}$ 
choices of $Q^e$, of which at most 
$1.1 d(L)^{Q-\tbinom{k+r}{r}} n^{q-r-k}$ contain $e'$,
so $\mb{P}(e' \sub Q^e) < 
1.3 (q-r)! d(L)^{1-\tbinom{k+r}{r}} n^{-k}$.
Summing over $k$ gives the claim.

Now for any $f \in \tbinom{[n]}{r-1}$,
by typicality $|L(f)| < 1.1 d(L) n$,
so $\pl^\pm \Phi(f)$ are pseudobinomial
with mean at most $1.5 (q-r)! 2^r \tT n$ by the claim,
so w.h.p.\ $J'$ and $\pl^\pm \Phi$ are $q^q \tT$-bounded. \qed

\medskip

We conclude by proving the
main result of this section.

\medskip

\nib{Proof of Lemma \ref{bddint}.}
Suppose $J \in \mb{Z}^{K^r_n}$ is
$K^r_q$-divisible and $\tT$-bounded.
By Lemma \ref{bddint:flat}
there is some $J^0 \in \mb{Z}^{K^r_n}$
and $\Phi^0 \in \mb{Z}^{K^r_q(K^r_n)}$ 
such that $\pl \Phi^0 = J-J^0$,
all $|J^0_e| < n^{0.1}$ and $J^0$ and
$\pl^\pm \Phi^0$ are $q^q \tT$-bounded.

Let $L \sim K^r_n(\nu)$, where $\nu = n^{-1/3Q}$.
By Lemma \ref{extrandom} w.h.p.\ $L$ is $(n^{-1/9},Q)$-typical
and by Lemma \ref{Jrandom}
w.h.p.\ $J^0$ is $(1.1 \cdot q^q \tT, Q)$-bounded w.r.t.\ $L$.
As w.h.p.\ $L$ is $1.1\nu$-bounded,
by Lemma \ref{bddgen} there is $S \sub K^r_q(K^r_n)$
such that $\pl S$ is $3q! (8r)^r \nu^{1/b}$-bounded and
$\bgen{K^r_q(K^r_n)} \cap \mb{Z}^L \sub \bgen{S}$.

By Lemma \ref{bddint:random}
there is some $J^1 \in \mb{Z}^L$
and $\Phi^1 \in \mb{Z}^{K^r_q(K^r_n)}$ 
such that $\pl \Phi^1 = J^0-J^1$, 
and $J^1$ and $\pl^\pm \Phi^1$ are $2q^{2q} \tT$-bounded.
As $J^1 \in \bgen{K^r_q(K^r_n)} \cap \mb{Z}^L \sub \bgen{S}$
there is $\Psi \in \mb{Z}^S$ with $\pl \Psi = J^1$.

Let $\Phi^2 \in [N]^S$ be such that 
$\Psi - \Phi^2 \in N\mb{Z}^S$.
Then $\pl \Phi^2$ is $3q! (8r)^r \nu^{1/b} N$-bounded 
(as $\pl S$ is $3q! (8r)^r \nu^{1/b}$-bounded)
and $J^2 := J^1 - \pl \Phi^2
= \pl (\Psi-\Phi^2) \in N\mb{Z}^{K^r_n}$
is $4q^{2q} \tT$-bounded, as $\tT \gg \nu^{1/b}$.
By Lemma \ref{bddint:N}
there is $\Phi^3 \in \mb{Z}^{K^r_q(K^r_n)}$ such that 
$\pl \Phi^3 = J^2$ and $\pl^\pm \Phi^3$ are 
$(2q)^{3q} N \tT$-bounded.

Let $\Phi=\sum_{i=0}^3 \Phi^i$. Then $\pl \Phi = J$
and $\pl^\pm \Phi$ are $N^2 \tT$-bounded. \qed

\section{Absorption} \label{sec:ab}

\begin{figure}
\begin{center}
\includegraphics{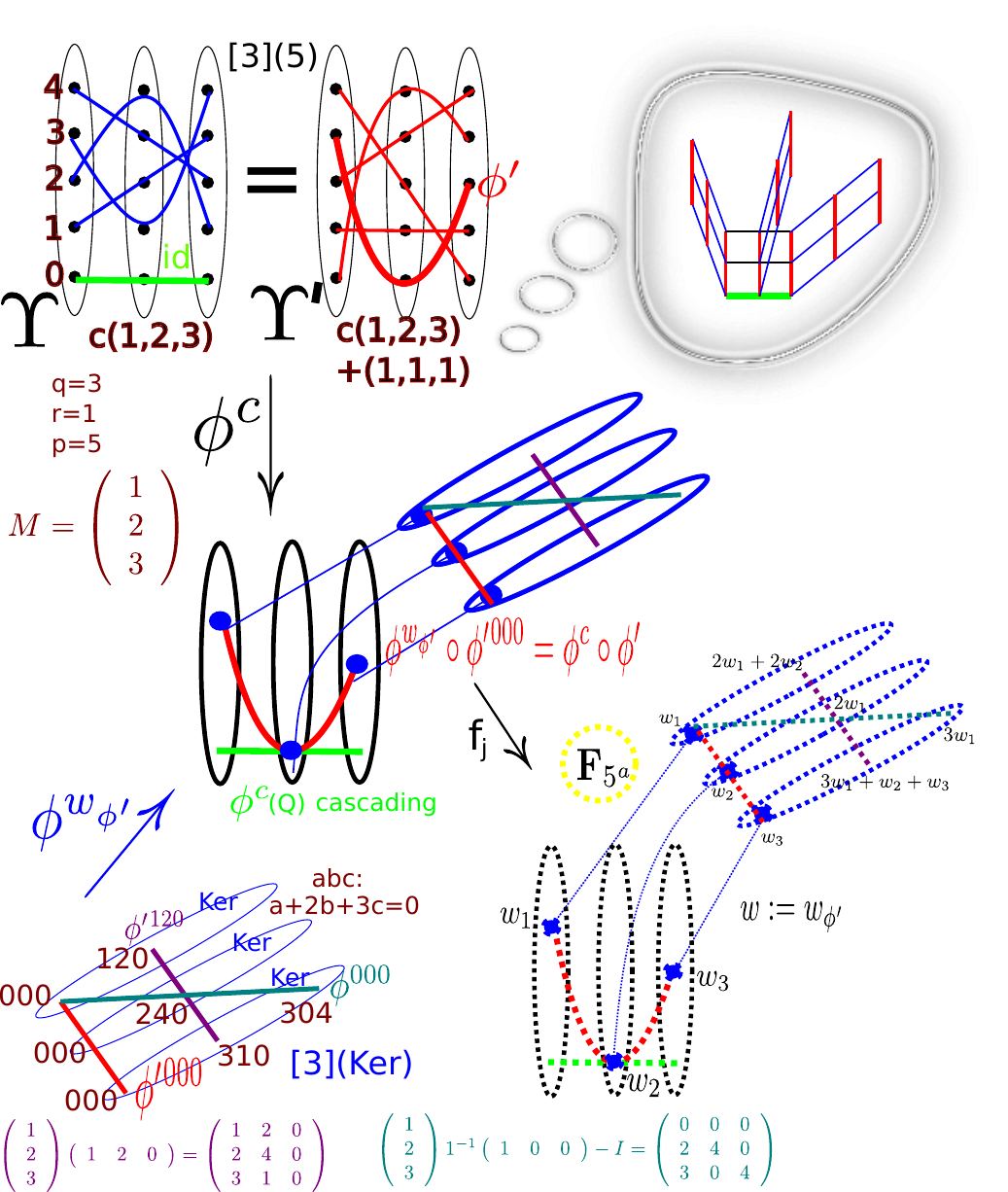}
\caption{A cascade}
\end{center}
\label{fig:ab}
\end{figure}

Now we describe the structure of
absorbable cliques in the template;
it is here that the algebraic properties 
of the template construction will come into play.
As this section is rather technical,
we start by illustrating the constructions
in the first subsection,
with reference to Figure \ref{fig:ab},
in the case $q=3$ and $r=1$,
that is, $3$-graph matchings
(it would be hard to make
a figure for $r>1$).
In the second subsection we construct absorbers.
The third subsection combines absorbers
to create cascades. The last subsection
obtains lower bounds on extensions involving
cascading cliques that are required for the
analysis of the Clique Exchange Algorithm 
in Section \ref{sec:cea}.

\subsection{Illustrations}

We start with the `thought bubble' in the top
right of the picture, which contains 
a `cartoon cascade'.
The blue diagonal triples represent
some triples of the template.
The green horizontal triple at the bottom
represents the `target':
we want to modify the template 
so that it contains the green triple,
without changing the set of vertices that it covers.
To achieve this, we first replace
the blue triples by the vertical red triples,
which is valid as they are both matchings
covering the same set of vertices.
Then the three vertical red triples in the square
can be replaced by three horizontal triples
that cover the same vertices, and include
the green triple, as desired.

The cartoon cascade was obtained by gluing
together four copies of a simpler structure,
namely a set of nine vertices with two
decompositions into three triples.
Three of these copies use template edges,
and correspond to what we later call `absorbers':
these are subsets (in general subgraphs)
of the template with two decompositions,
one of which only uses template triples
(in general $q$-cliques).
The red triples in the picture correspond
to cliques that we will call `absorbable':
these can be included in the template
by `flipping' the relevant absorber,
with no need for a cascade.

The reader may wonder why we do not also 
describe the green triple as `absorbable',
given that it is obtained by the net result
of the above replacements, which take the
nine blue template triples and replace
them by nine other triples 
that include the green one.
The reason is that the algebraic structure
naturally associates to any clique
a simple configuration that acts as
an absorber if it is present in $G$
(e.g. for triangle decompositions in \cite{K2}
we associate octahedra to triangles).
Thus we have a naturally defined subfamily
of cliques with simple absorbers,
which we combine into more complicated
structures (cascades) that absorb
a larger family of `cascading' cliques.

In our illustration we glue three
absorbers onto a `base', which we chose to be 
isomorphic to an absorber. However, this is not necessary,
and in general it will be convenient
to use a different structure for the base,
which is simpler than that of the absorbers.
 
Now we turn to the details of an actual cascade,
in the case $q=3$ and $r=1$. We will use the prime $p=5$
(which is not as large as advertised elsewhere,
but the construction still works).
We fix the generic $3 \times 1$ matrix $M = (1\ 2\ 3)^t$.

The top left of the picture illustrates the `blueprint'
for the base of the cascade, which consists
of two perfect $3$-graphs matchings 
$\Ups$ and $\Ups'$ on a set of $15$ points
(the same set, drawn twice for clarity), 
divided into $3$ parts of size $5$,
where each triple is transverse to the partition.
Reading each triple of $\Ups$ or $\Ups'$ as a vector, 
$\Ups$ consists of all $c(1,2,3)$  
and $\Ups'$ of all $c(1,2,3)+(1,1,1)$,  
where $c \in \mb{F}_5$.

The base of any cascade is defined
by some embedding $\phi^c$ of this blueprint
of the base in the template.
Note that here `embedding' only constrains
the vertices (in general $r$-edges);
the triples (in general $q$-cliques)
are contained in the underlying graph of the template
but may not belong to the template decomposition.

Similarly to the cartoon cascade, the cascade will flip
in two stages. The first stage will provide absorbers
for the cliques in $\Ups'$, which can be
flipped so that all cliques of $\Ups'$ 
are present in the decomposition.
The second stage is to flip the base;
we replace $\Ups'$ by $\Ups$.

The green triple $(0,0,0)$ of $\Ups$ is mapped by $\phi^c$ 
to the target of the cascade. 
It is notationally convenient to identify
$\mb{F}_5$ with $[5]$ so that $0$ is identified
with $1$, and identify the vertices of the green triple
with $[3]$. Recalling our notation to identify vectors
with functions, the green triple is thus identified
with $id_{[3]}$ (the identity map on $[3]$),
so the target clique is $\phi^c(Q)=\phi^c([3])$.

We require each clique of $\Ups'$ to be absorbable,
so the remainder of the cascade will be defined by
gluing absorbers onto these cliques.
We illustrate this for the clique labelled $\phi'$,
where $\phi'(1)=3$, $\phi'(2)=1$, $\phi'(3)=2$.
In the centre of the figure this is the red clique,
which has been drawn twice for clarity,
once in the base of the cascade, and once in an absorber, 
where three vertices of the absorber are identified
with the corresponding vertices of the base,
and the absorber is otherwise vertex-disjoint
from all other parts of the construction.

The blueprint for absorbers is illustrated
in the bottom left of the figure. 
Similarly to the base of the cascade,
it consists of two perfect $3$-graph matchings
of the same set of points, divided into $3$ parts,
so that each triple is transverse to the partition.
However, now the parts have size $25$, and
each is identified with the left kernel of $M$:
the set of vectors $(a,b,c) \in \mb{F}_5^3$
(also written as $abc$) with $a+2b+3c=0$.

The absorbers in the cascade are defined by
various embeddings of the blueprint absorber.
These embeddings are specified with reference
to one of the template embeddings 
$f_j:[n] \to \mb{F}_{p^a}$, where each
cascade fixes some $j \in [z]$ 
for all of its absorbers.
Each clique $\phi'$ of $\Ups'$ corresponds to some 
$w_{\phi'} = f_j \phi^c \phi' \in (\mb{F}_{5^a})^3$,
a vector that can be identified with
a function $w_{\phi'}: [3] \to \mb{F}_{5^a}$
where each $w_{\phi'}(i) = f_j(\phi^c(\phi'(i)))$.

For convenient notation in the remainder of
the illustration we fix $\phi'$ and write $w:=w_{\phi'}$.
The actual absorber for $\phi^c \phi'$
(the red clique in the middle of the figure)
is obtained by embedding the blueprint absorber.
This embedding is specified by a map $\phi^w=\phi^{w_{\phi'}}$
satisfying $f_j\phi^w((i,a)) = w_i + a \cdot w$,
where $(i,a)$ denotes the copy of $a$ in the $i$th part,
for any $i \in [3]$ and $a \in Ker$.
We require the base embedding $\phi^c$
to be such that $\{w_1,w_2,w_3\}$ 
has full dimension (viewing $\mb{F}_{5^a}$
as a vector space over $\mb{F}_5$);
it then follows that $\phi^w$ is injective.

To relate the red clique $\phi^c \phi'$
to the embedding $\phi^w$, we note that each
$f_j \phi^c \phi'(i) = w_i + 0 \cdot w$, where 
$0=000$ is the zero vector in $\mb{F}_{5^a}$.
We view triples in the blueprint absorber
as $3 \times 3$ matrices in which the 
$i$th row is the vector corresponding 
to the vertex chosen from the $i$th part.
Thus $\phi^c \phi'=\phi^w \phi'{}^0$,
where $\phi'{}^0=\phi'{}^{000}$ maps
each $i \in [3]$ to $000$, and so can
be viewed as the $3 \times 3$ zero matrix.

The essential feature of absorbers is that
they have two decompositions, one of which
uses the target absorbable clique
$\phi^c \phi'=\phi^w \phi'{}^0$,
and the other of which is contained
within the template decomposition.
We can specify these decompositions in the
blueprint absorber and then transfer
them to the absorber via $\phi^w$.
One decomposition consists of all
triples $\phi'{}^a$ with $a \in Ker$,
specified in $3 \times 3$ matrix form
as the outer product $Ma$.
Concretely, for each $a \in Ker$ the triple 
$\phi'{}^a$ uses vertex $ia = (ia_1,ia_2,ia_3)$ 
in part $i$ for $i \in [3]$.
(This agrees with our above
notation for $\phi'{}^0$.)
The purple clique illustrates
this for $a = 120$.

The other decomposition consists of all
triples $\phi^a$ with $a \in Ker$,
specified in $3 \times 3$ matrix form
as $\phi^a = M(a+e_1)-I$, where $e_1=(1\ 0 \ 0)$;
the teal clique illustrates this for $a=000$.
Note that all such triples are contained
in the blueprint absorber, 
as $\phi^a M = M(a+e_1)M-M = Me_1 M - M = 0$,
so each row of $\phi^a$ is in $Ker$.
Furthermore,
as $f_j\phi^w((i,b)) = w_i + b \cdot w$,
we have $f_j \phi^w \phi^a = w + \phi^a w
= w + (M(a+e_1)-I)w = M(a+e_1)w$:
 each $f_j \phi^w \phi^a$ is in the image of $M$,
so $\phi^w \phi^a$ can be a template clique
(if the activation and compatibility conditions
of the template construction also hold).

\subsection{Absorbers}

Now we will implement the previous illustration
in our general setting.
The construction of absorbers will use 
the left kernel of $M$: let\footnote{
Our notation uses `a' in many ways:
here it is a vector, 
later it will be a matrix,
and throughout it appears
in the notation for the field $\mb{F}_{p^a}$;
we hope that the intended uses
will be clear from their contexts.} 
\[Ker := \{ a \in \mb{F}_p^q: aM=0 \}.\]
The following properties of $Ker$ are
immediate from the construction of $M$,
so we omit their proofs.

\begin{lemma} \label{Kprops} $ $
\begin{enumerate}
\item $\dim(Ker)=q-r$, so $|Ker|=p^{q-r}$,
\item if $a \in Ker$ with $a \ne 0$
then $|\{i: a_i \ne 0\}|>r$,
\item for any $I \in Q$ and $i \in [q] \sm I$,
the unique $a \in Ker$ with $a_i=-1$
and $a_j = 0$ for all $j \in [q] \sm (I \cup \{i\})$
is $a = e_i M_I^{-1} e_I - e_i$.
\end{enumerate}
\end{lemma}

Given $a = (a^i: i \in I)$ with each $a^i \in \mb{F}_p^q$,
we identify $a$ with a matrix $a \in \mb{F}_p^{I \times [q]}$ 
having entries $(a^i_{i'}: i \in I, i' \in [q])$.
For $a \in Ker^r \sub \mb{F}_p^{[r] \times [q]}$ 
and $w \in \mb{F}_{p^a}^q$ we write
\[ v^w_a = M M_{[r]}^{-1} (e_{[r]} + a) w
\quad \text{ and } \quad
v'{}^w_a = w+Maw .\]
For example, if $a=0 \in Ker^r$ then $v'{}^w_0=w$
and $v^w_0 = M M_{[r]}^{-1} w_{[r]}$ is a vector
in the image of $M$ that might correspond 
to a template clique containing
an edge that corresponds to $w_{[r]}$.

Now we come to the key definition of this subsection,
which implements the general form of the construction
of absorbers illustrated in the previous subsection.
An absorber is defined by an embedding $\phi^w$ 
of a `blueprint absorber',
which consists of a complete $q$-partite
$r$-graph with vertex set $[q] \times Ker$.
For any fixed clique $\phi(Q) \sub G^*_j$,
the vector $w:=f_j\phi \in \mb{F}_{p^a}^q$
determines the embedding $\phi^w$ of the blueprint
to the absorber $A^{\phi(Q)}$ for $\phi(Q)$
via point (i) of the definition.
Point (ii) of the definition requires
certain cliques in the absorber 
to belong to the template $M^*_j$;
we will see in Lemma \ref{shuffle}
that these form a $K^r_q$-decomposition
$\Psi(\phi^w)$ of the absorber $A^{\phi(Q)}$.
We write $[q](Ker)$ for the set of partite maps
$f:[q] \to [q] \times Ker$,
meaning that each $f(i)$ is some $(i,a)$ with $a \in Ker$.

\begin{defn} \label{def:absorb} (absorbers)
Suppose $\phi(Q) \sub G^*_j$ with $j \in [z]$ and
$w:=f_j\phi \in \mb{F}_{p^a}^q$ has $\dim(w)=q$.

Suppose $\phi^w: [q] \times Ker \to [n]$ such that
\begin{enumerate}
\item $f_j\phi^w((i,a)) = w_i + a \cdot w$ 
for each $i \in [q]$, $a \in Ker$,
\item if $\phi' \in [q](Ker)$ with
$f_j \phi^w \phi' = v^w_a$ for some $a \in Ker^r$
then $\phi^w \phi'(Q) \in M^*_j$. 
\end{enumerate}

We say that $\phi$ is absorbable
and call $\phi^w$ the absorber for $\phi$.

We also call $A^{\phi(Q)}=A^w=\phi^w(K^r_q(Ker))$
the absorber for $\phi(Q)$.%
\footnote{Recall that $K^r_q(Ker)$ is the complete
$q$-partite $r$-graph with each part
identified with $Ker$.} 
\end{defn}

The essential property of absorbers 
(see Lemma \ref{shuffle} below)
is that they can be decomposed in two ways, 
one of which uses cliques that all belong to template, 
and the other of which uses any 
absorbable `target' clique.

First we make some comments on the definition.
The notation $A^w$ is ambiguous, but will be clear from
the context, as $j$ in Definition \ref{def:absorb}
is uniquely determined by $\phi(Q) \sub G^*_j$.
The notation $A^{\phi(Q)}$ is also ambiguous in that
we could reorder $\phi$ without changing $\phi(Q)$,
but the order will be clear from the context
(we will only consider $\pi$-compatible $\phi$).

Next we introduce some notation for edges in absorbers.
The edges of the complete $q$-partite $r$-graph
$A^w=\phi^w(K^r_q(Ker))$ correspond 
to choices $I \in Q = \tbinom{[q]}{r}$ of $r$ parts
and any choices of vertices $\phi^w((i,a^i))$
in these parts for each $i \in I$.
We identify $(a^i: i \in I)$ with $a \in Ker^I$
and denote the corresponding edge by $e^w_a$.
By Definition \ref{def:absorb}.i we have $f_j(e^w_a) = (e_I+a)w$. 

An easy but important property of absorbers 
is established by the following lemma,
which shows that all edges
have full dimension in their relevant embedding
(so, in particular, all $w_i + a \cdot w$ 
are distinct, so $\phi^w$ is injective).

\begin{lemma} \label{absdim}
Suppose $A^w$ is the absorber of $Q' \in K^r_q(G^*_j)$.
Then each $\dim(f_j(e^w_a))=r$.
\end{lemma}

\nib{Proof.}
Suppose $a \in Ker^I$ and $c \in \mb{F}_p^r$ with $c (e_I + a)w = 0$.
As $\dim(w)=q$ we must have $c(e_I+a)=0$.
Then $0 \ne -ce_I = ca \in Ker$, has at most $r$
nonzero coordinates, contradicting Lemma \ref{Kprops}. \qed

\medskip

We require some more notation to specify the
clique decompositions of absorbers.
We can write Definition \ref{def:absorb}.i
as $f_j \phi^w \phi' = w + \phi' w$ 
for all $\phi' \in [q](Ker)$,
viewing $\phi'$ as a matrix in $\mb{F}_p^{q \times q}$.
For $a \in Ker^r$ we define $\phi^a$ and $\phi'{}^a$
in matrix form by
\begin{equation} \label{eq:phia}
\phi^a:=M M_{[r]}^{-1} (e_{[r]} + a)-I
\quad \text{ and } \quad
\phi'{}^a = Ma .
\end{equation}
Then for $\phi^w$ as in Definition \ref{def:absorb}.i
we have
\[ f_j \phi^w \phi^a = w + \phi^a w = v^w_a 
\quad \text { and } \quad
f_j \phi^w \phi'{}^a = w + \phi'{}^a w = v'{}^w_a. \]
We write
\[ \Psi(\phi^w) = \{ \phi^w\phi^a(Q): a \in Ker^r\}
\quad \text{ and } \quad
\Psi'(\phi^w) = \{ \phi^w\phi'{}^a(Q): a \in Ker^r\}. \]
Note that $\Psi'(\phi^w)$ contains
the clique $\phi^w\phi'{}^0(Q)$ with
vertex set $f_j^{-1}(v'{}^w_0)=f_j^{-1}(w)=Im(\phi)$,
and $\Psi(\phi^w) \sub M^*_j$ by Definition \ref{def:absorb}.ii.
Thus the following lemma shows that the absorber for $\phi(Q)$
can be used to modify the template,
replacing $\Psi(\phi^w)$ by $\Psi'(\phi^w)$,
so that it contains $\phi(Q)$
(we say that we `flip' $A^{\phi(Q)}$).

\begin{lemma} \label{shuffle}
$\Psi(\phi^w)$ and $\Psi'(\phi^w)$
are both $K^r_q$-decompositions of $A^{\phi(Q)}$.
\end{lemma}

\nib{Proof.}
First we claim that each clique 
in $\Psi(\phi^w)$ and $\Psi'(\phi^w)$ intersects each 
part $A^w_i = f_j^{-1}\{  w_i + a \cdot w: a \in Ker \}$.
To see this, note that
a $q$-set intersects each $A^w_i$ if and only if
it can be written as $f_j^{-1}((I+B)w)$ 
for some $q \times q$ matrix $B \in Ker^q$.
As $f_j \phi^w\phi'{}^a  = v'{}^w_a = (I+Ma) w$
and $f_j \phi^w\phi^a  = v^w_a 
= M M_{[r]}^{-1} (e_{[r]} + a) w$,
the claim follows from $(I+Ma)M=IM=M$
and $M M_{[r]}^{-1} (e_{[r]} + a)M 
= M M_{[r]}^{-1} M_{[r]} = M$.

Now consider any $e^w_a \in A^w$,
where $a \in Ker^I$ for some $I \in Q$. Then 
$e^w_a \in \phi^w\phi'{}^{M_I^{-1}a}(Q) \in \Psi'(\phi^w)$,
and $e^w_a \in \phi^w\phi^{a'}(Q) \in \Psi(\phi^w)$ where
$e_I + a = M_I M_{[r]}^{-1} (e_{[r]} + a')$,
that is, $a' = M_{[r]} M_I^{-1} (e_I + a) -  e_{[r]}$
(note that $a'M = M_{[r]} M_I^{-1} M_I - M_{[r]}=0$).
As $|\Psi(\phi^w)|=|\Psi'(\phi^w)|=|Ker|^r=Q^{-1}|A^w|$,
the lemma follows. \qed

\subsection{Cascades}

Absorbable cliques are plentiful but not ubiquitous.
Here we will describe a much wider class of cliques
that can be included in the template via a series
of modifications using absorbable cliques.

First we describe our clique exchange tool,
which will also be used in Section \ref{sec:cea}
for the Clique Exchange Algorithm.
It consists of two suitable
decompositions of a small fixed $r$-graph:
we use the complete $q$-partite $r$-graph
$K^r_q(p)$ with $p$ vertices in each part.

\begin{lemma} \label{exchange}
There are $K^r_q$-decompositions  $\Ups$ and $\Ups'$ 
of $\OO=K^r_q(p)$ such that
\begin{enumerate}
\item $|V(f) \cap V(f')| \le r$
for all $f \in \Ups$ and $f' \in \Ups'$,
\item if $f \in \Ups$ and $\{f',f''\} \sub \Ups'$
with $|V(f) \cap V(f')|=|V(f) \cap V(f'')|=r$

then $(V(f') \sm V(f)) \cap (V(f'') \sm V(f))=\es$.
\end{enumerate}
\end{lemma}  

The construction requires a matrix of the same type 
as that used in constructing the template, 
with an additional technical property.

\begin{defn} \label{defgen'}
Let $M' \in \mb{F}_p^{q \times r}$ be 
such that every square submatrix of $M'$ is nonsingular
and for any $r \times r$ submatrix $A$ of $M'$
and row $v$ of $M'$ not in $A$ each
entry of $vA^{-1}$ is not $0$ or $1$.
\end{defn}

To see that such $M'$ exists 
we again consider a uniform random $M'$,
and recall from the construction of $M$
that the probability of having
any singular square submatrix is at most $2^{q+r+1}p^{-1}$.
Now fix any $r \times r$ submatrix $A$ of $M'$
and row $v$ of $M'$ not in $A$:
there are fewer than $q2^q$ choices.
There are fewer than $2r p^{r-1}$ row vectors 
$w \in \mb{F}_p^r$ such that some entry is $0$ or $1$.
We fix any such $w$ and bound $\mb{P}(Aw=v)$.
We can assume $w \ne 0$, as otherwise $v=0$,
so $M'$ has zero entries, which are
singular $1$ by $1$ submatrices.
Without loss of generality $w_1 \ne 0$.
We condition on any value of $v$ and
all but the first column of $A$.
Then $Aw$ is uniformly random,
so $\mb{P}(Aw=v) = p^{-r}$.
Thus the required properties of $M'$
fail with probability at most
$2^{q+r+1}p^{-1} + 2r q2^q p^{-1} < 1$,
so $M'$ exists.

\medskip

\nib{Proof of Lemma \ref{exchange}.}
We identify each part $\OO_i$ of $\OO$ with $\mb{F}_p$.
We let $\Ups$ consist of all $q$-cliques $Q'$ 
of the form $M'x$, that is, for some $x \in \mb{F}_p^r$ 
we have $V(Q') \cap \OO_i = (M'x)_i$ for all $i \in [q]$.
We choose $c \in (\mb{F}_p \sm \{0\})^q$ 
uniformly at random and let $\Ups'$ consist of all 
$q$-cliques $Q'$ of the form $M'x+c$.

Now for any $I \in Q$ and $z_I \in \mb{F}_p^I$,
there is a unique $q$-clique $v = M'(M'_I)^{-1}z_I$ 
in $\Ups$ containing $z_I$, and 
a unique $q$-clique $v' = M'(M'_I)^{-1}(z_I-c_I)+c$ 
in $\Ups'$ containing $z_I$.
Thus $\Ups$ and $\Ups'$ are $K^r_q$-decompositions of $\OO$.

To show properties (i) and (ii)
we show that the failure of each corresponds 
to a nontrivial linear equation in $c$,
so with positive probability 
there is some $c$ such that (i) and (ii) hold.

If property (i) failed we would have
$i \notin I \in Q$ with
$c_i = M'_i (M'_I)^{-1} c_I$,
which is an equation for $c$
with a nonzero coefficient of $c_i$
by construction of $M'$
(no entry of $M'_i (M'_I)^{-1}$ is equal to $1$).

If property (ii) failed we would have
$\{I,I'\} \sub Q$ and $i \notin I \cup I'$ 
with $M'_i (M'_I)^{-1} c_I = M'_i (M'_{I'})^{-1} c_{I'}$.
We can choose $i' \in I \sm I'$,
and then $c_{i'}$ appears with 
nonzero coefficient in the equation
(no coefficient of $M'_i (M'_I)^{-1}$ is equal to $0$).

This gives at most $qQ+Q^2$ equations for $c$,
each holding with probability at most $(p-1)^{-1}$,
so we can choose $c$ such that (i) and (ii) hold. \qed

\medskip

We identify $\Ups$ and $\Ups'$ with subsets of 
the set $[q](p)$ of partite maps from $[q]$ to $[q] \times [p]$.
We identify $[q]$ with $\{ (i,1): i \in [q] \} \sub V(\OO)$
and with the corresponding map $id_{[q]}$;
by relabelling we can assume $[q] \in \Ups$.
Next we require some more terminology.

\begin{defn}
Suppose $U' \sub U \sub [n]$. 
We say that $U$ is $j$-generic for $U'$ 
if $\dim(f_j(U))=\dim(f_j(U'))+|U|-|U'|$.
We say that $U$ is generic for $U'$ 
if $U$ is $j$-generic for $U'$ for all $j \in [z]$.
\end{defn}

Note that given $U' \in \tbinom{[n]}{u'}$, 
all but $O(n^{u-u'-1})$ sets $U \in \tbinom{[n]}{u}$
containing $U'$ are generic for $U'$.
Now we can give the key definition of this subsection.

\begin{defn} \label{def:cascade} (cascades)
Suppose $Q' = \phi(Q) \sub G^*_j$ 
and $\phi^c$ is an embedding of $K^r_q(p)$ in $G^*_j$
where $\phi^c id_{[q]} = \phi$ and $Im(\phi^c)$ 
is $j$-generic for $Im(\phi)$, such that each 
$\phi^c\phi'$ with $\phi' \in \Ups'$ is absorbable,
with absorber $A^{\phi^c\phi'(Q)}=\phi^{w_{\phi'}}(K^r_q(Ker))$, 
and $C_{\phi^c} = \sum \{A^{\phi^c\phi'(Q)}: \phi' \in \Ups'\}$
is a set (without multiple elements). 
We call $C_{\phi^c}$ a cascade for $Q'$.
\end{defn}

A cascade for $Q'$ provides a two-step process
for modifying the template so as to include $Q'$:
we flip all of the absorbers in the cascade,
and then flip the $K^r_q$-decomposition
of the base embedding of $\OO$.
Formally, to flip a cascade $C_{\phi^c}$ we replace 
\[\Psi(C_{\phi^c}) := \bigcup \{ \Psi(w_{\phi'}): \phi' \in \Ups' \}
\qquad \text{ by } \]
\[ \Psi'(C_{\phi^c}) := \{ \phi^c\phi'(Q): \phi' \in \Ups \} \cup \bigcup 
\{  \Psi'(w_{\phi'}) \sm \{\phi^c\phi'(Q)\}: \phi' \in \Ups' \}.\]
This has the desired property as
$\phi(Q) = \phi^c (Q) \in \Psi'(C_{\phi^c})$.
Next we define the set of cliques 
for which we will show (Lemma \ref{cascade})
that we have many cascades.

\begin{defn} \label{def:cascading} (cascading cliques)
Let $\mc{Q}^* = \cup_{j \in [z]} \mc{Q}_j$,
where each $\mc{Q}_j$ is the set of all $\phi(Q) \sub G^*_j$
where $\phi$ is $\pi$-compatible and $\dim(f_j\phi)=q$.
\end{defn}

In analysing the choice of cascades, we will often
need to know how fixing the image of one edge
constrains the possible images of some other edge.
We will define bipartite graphs $B^{aI'}_{j\phi'}$ 
describing which pairs of edges $(e,e')$ in some $G^*_j$ 
can satisfy $e=e^{w_{\phi'}}_a$ and $e'=\phi^c(I')$.
In the accompanying lemma we show that
$f_j$ embeds $B^{aI'}_{j\phi'}$
in an algebraically defined regular 
bipartite graph $AB^{aI'}_J$,
for which we can describe neighbourhoods 
as certain affine linear spaces, 
all of the same dimension.

\begin{defn} \label{def:ee'} 
(cascade edge compatibility graphs)

Let $\phi' \in \Ups'$ and $J \sub I' \in Q$
where $J = \{i \in [q]: \phi'(i)=(i,1)=i\}$.
Let $a \in Ker^I$ with $I \in Q$.

Let $B^{aI'}_{j\phi'} \sub G^*_j \times G^*_j$
be the bipartite graph where $(e,e')$ is an edge
if there is a cascade $C_{\phi^c}$
for some $Q'=\phi^c([q])$ with $\phi^c(I')=e'$ where 
$A^{\phi^c\phi'(Q)}=\phi^{w_{\phi'}}(K^r_q(Ker))$
with $e=e^{w_{\phi'}}_a$.

Let\footnote{Note that the use of `$a$' in $\mb{F}_{p^a}$
is unrelated to its use in this definition.} 
 $AB^{aI'}_J \sub \mb{F}_{p^a}^I \times \mb{F}_{p^a}^{I'}$
be the bipartite graph where $(x,x')$ is an edge
if there is a solution $w \in \mb{F}_p^q$ 
to the simultaneous equations
$w_i + a_i \cdot w = x_i$ for $i \in I$
and $w_i = x'_i$ for $i \in J$. 

We let $r_{a\phi'} = \rk Y$ where
$Y = (e_I + a)e_{[q] \sm J}^T$.
\end{defn}

\begin{lemma} \label{lem:ee'} 
With notation as in Definition \ref{def:ee'},
\begin{enumerate}
\item if $(e,e') \in B^{aI'}_{j\phi'}$
then $(f_j(e),f_j(e')) \in AB^{aI'}_J$,
regarded in $\mb{F}_{p^a}^I \times \mb{F}_{p^a}^{I'}$
via $(\pi_e,\pi_{e'})$,
\item every vertex neighbourhood in $AB^{aI'}_J$
is an affine linear space of dimension $r_{a\phi'}$.
\end{enumerate}
\end{lemma}

Here (and later) we require the following easy fact
from linear algebra (we omit the proof).

\begin{lemma} \label{lem:linalg}
Let $\mb{F}$ be a field, $A \in \mb{F}^{a \times q}$,
$B \in \mb{F}^{b \times q}$ and 
$c \in \{ Bw: w \in \mb{F}^q \}$.
Then $\{Aw: Bw=c\}$ is an affine linear space
of dimension $\rk \tbinom{A}{B} - \rk B $.
\end{lemma}

\begin{proof}[Proof of Lemma \ref{lem:ee'}]
To see (i), note that
if $e=e^w_a$ with $w=w_{\phi'}$ then
$f_j(e) = (e_I + a) w$, which has $i$-coordinate 
$w_i + a_i \cdot w$ for each $i \in I$,
and if $\phi^c(I')=e'$ then
$f_j \phi^c (i) = w_i$ for each $i \in J$.

For (ii), we write $e_I + a = (X \ Y)$
where $X = (e_I + a)e_J^T$
and $Y = (e_I + a)e_{[q] \sm J}^T$,
that is, $X$ contains the columns indexed by $J$
and $Y$ those indexed by $[q] \sm J$.

Consider any $x' \in \mb{F}_{p^a}^{I'}$.
If $(x,x')$ is an edge with some solution $w$
then we can write 
$x = (e_I + a)w = Xw_J + Yw_{[q] \sm J}
= Xx'_J + Yw_{[q] \sm J}$,
so $x$ must lie in the affine space
$Xx'_J + \text{Im}(Y)$.

Now consider any $x \in \mb{F}_{p^a}^I$.
In choosing $x' \in \mb{F}_{p^a}^{I'}$,
the coordinates $x'_i$ with $i \in I' \sm J$
are unconstrained, whereas those in $J$ satisfy
$x'_J \in \{ e_J w: (e_I + a)w = x \}$.
By Lemma \ref{absdim}, $x'$ lies in an 
affine space $W$ with
$\dim(W) = r-|J|+\rk Z - \rk (e_I+a)$,
where $Z = \bigl( \begin{smallmatrix}
e_J&0\\ X&Y
\end{smallmatrix} \bigr)$.
Note that $Z$ is row-equivalent to
$\bigl( \begin{smallmatrix}
e_J&0\\ 0&Y
\end{smallmatrix} \bigr)$,
so has rank $|J|+r_{a\phi'}$,
and $e_I+a$ has rank $r$ by Lemma \ref{absdim}.
Thus $\dim(W) = r-|J|+|J|+r_{a\phi'}-r=r_{a\phi'}$, as required.
\end{proof}

\begin{rem} \label{rem:nondegen}
We will only ever consider $\phi'$ and $a$
as in Definition \ref{def:ee'} 
such that setting $e=e^{w_{\phi'}}_a$
and $e'=\phi^c(I')$ does not imply
that $e$ and $e'$ belong to the same 
clique of the template $M^*$.
This condition is equivalent to $r_{a\phi'} > 0$,
that is, $Y = (e_I + a)e_{[q] \sm J}^T \ne 0$.
To see this, we note for any $i \in I$
that if row $i$ of $Y$ is $0$ 
then $1_{i=j} + a_{ij}=0$
for all $j \in [q] \sm J$,
which by Lemma \ref{Kprops}.iii
implies $i \notin J=I' \in Q$ and
$a_i = e_i M M_{I'}^{-1} - e_i$.
However, if this holds for all $i \in I$
then $e'=\phi^c(I')$ and
$e=e^{w_{\phi'}}_a$ imply that
$e$ is the $I$-edge of $M^*(e')$. 
\end{rem}

The proof of concentration in Lemma \ref{cascade} 
uses the following upper bound on the number of cascades
for a given clique using a given edge;
this bound will also be used in the analysis of the
cascade algorithm in the proof of Theorem \ref{main+}.
  
\begin{lemma} \label{cascade:role}
Suppose $Q' = \phi(Q) \in \mc{Q}_j$ and $e \in G^*$.
Let $\phi' \in \Ups'$, 
$I \in Q$, $a \in Ker^I$.
Then there are at most $(p^a)^{q(p-1)-r_{a\phi'}}$ 
cascades $C_{\phi^c}$ for $Q'$ such that the absorber
$A^{\phi^c\phi'(Q)}=\phi^{w_{\phi'}}(K^r_q(Ker))$
for $\phi^c\phi'$ satisfies $e^{w_{\phi'}}_a=e$.
\end{lemma}

\nib{Proof.}
Any cascade $C_{\phi^c}$ for $Q'$ determines some
$w^C = f^j\phi^c \in \mb{F}_{p^a}^{[q] \times [p]}$
with $w' := w^C_{[q]} = f_j \phi$, and any such $w^C$
corresponds to at most one cascade for $Q'$.
The condition $e^{w_{\phi'}}_a=e$ imposes the
additional constraint $(e_I + a)w = b$,
where $w := f_j \phi^c \phi'$ and $b := f_j(e)$
(regarded in $\mb{F}_{p^a}^I$ via $\pi_e$).
Let $J = \{i \in [q]: \phi'(i)=i\}$.
Noting that $w'_J=w_J$, we write 
$w^C = (w^{C'},w'_{[q] \sm J},w)$,
and apply Lemma \ref{lem:linalg}
to see that the choices for $w^C$
lie in an affine space $W$ with
$\dim W = \rk \tbinom{A}{B} - \rk B $,
where $A = I_{pq}$ (an identity matrix) 
and $B = \Big( \begin{smallmatrix}
0&0&0 \\ 0&I_{q-|J|}&0\\ 0&0&Z
\end{smallmatrix} \Big)$,
with $Z$ as in the proof of Lemma \ref{lem:ee'}.
Recalling that $\rk Z = |J| + r_{a\phi'}$, we have
$\dim W = pq - (q-|J|+|J|+r_{a\phi'}) = q(p-1)-r_{a\phi'}$.
\qed

\medskip

Now we give a lower bound on the number of cascades
on any cascading clique (to see that it is effective recall 
$p^{q^2} < 2^{9q^3} < h^{1/5}$ and $\oO>n^{-b^{-1}h^{-2}}$).

\begin{lemma} \label{cascade}
For any $Q' = \phi(Q) \in \mc{Q}_j$ there are
at least $\oO^{p^{q^2}} n^{q(p-1)}$ cascades for $Q'$.
\end{lemma}

\nib{Proof.}
We condition on local events $\mc{E} = \cap_{e \in Q'} \mc{E}^e$
such that $Q' \in \mc{Q}_j$. Then $\dim(f_j\phi)=q$ 
and $\phi$ is $\pi$-compatible, so for each $e \in \phi(Q)$ 
we can write $M^*(e)=\phi^e(Q)$ 
with $\pi_e \phi^e = \pi_e \phi = id$.
Let $U$ be the set of vertices touched by $\mc{E}$.

Now we consider any fixed combinatorial 
structure that could be a cascade for $Q'$ 
if it satisfies the necessary algebraic constraints.
We fix any embedding $\phi^c$ of $K^r_q(p)$ in $G$
with $\phi^c id_{[q]} = \phi$ and
$Im(\phi^c) \sm Im(\phi)$ disjoint from $U$;
recalling the illustration above, 
this specifies the base of the cascade,
and $\phi(Q)=\phi^c(Q)$ is represented by 
the green clique in Figure \ref{fig:ab}.

We also need to specify the combinatorial
structure of the absorbers.
For each $\phi' \in \Ups'$ we fix
any embedding $\phi^{\phi'}$ of $K^r_q(Ker)$ in $G$
with $\phi^{\phi'} \phi'{}^0 = \phi^c \phi'$,
recalling that $\phi'{}^0=(0,\dots,0)$
is identified with $[q]$; this is illustrated by 
the red clique $\phi^c \phi'(Q)$ in Figure \ref{fig:ab}
(we will add algebraic constraints below
so that $\phi^{w_{\phi'}}=\phi^{\phi'}$).

There is an additional constraint on $\phi^{\phi'}$
for each $\phi' \in \Ups'$ with $|Im(\phi') \cap [q]|=r$.
Indeed, then the red clique $\phi^c \phi'(Q)$
shares an edge with the green clique $\phi(Q)=\phi^c(Q)$;
such a $\phi'$ is illustrated in Figure \ref{fig:ab}
(where $r=1$, so an `edge' is a vertex).
If this edge is $e \in \phi(Q)$ 
we denote $\phi'$ by $\phi'_e$. The absorber
for $\phi^c \phi'_e(Q)$ must contain the
template clique $M^*(e)$ which contains $e$.
Accordingly, for each $e \in \phi(Q)$
we let $a_e \in Ker^r$ be such that 
$Im(\pi_e) \sub Im(\phi^{a_e})$,
where we identify $[q]$ with 
$[q] \times \{0\} \sub [q] \times Ker$
and recall $\phi^{a_e}$ from (\ref{eq:phia}).
Then $\phi^{a_e}$ must correspond
to the template clique containing $e$,
so we require $\phi^{\phi'_e} \phi^{a_e} = \phi^e$.

The final combinatorial condition on the cascade
is that the base and absorbers should be 
`as disjoint as possible' subject to the
gluing of the absorbers onto the base.
For each $\phi' \in \Ups'$, 
the set of `private' vertices
of the absorber for $\phi^c \phi'$ is
$I_{\phi'} = Im(\phi^{\phi'}) \sm Im(\phi^c\phi')$
if $\phi'$ is not some $\phi'_e$, 
or $I_{\phi'_e} = Im(\phi^{\phi'_e}) \sm 
( Im(\phi^c\phi'_e) \cup Im(\phi^e) )$.
We choose the $\phi^{\phi'}$ so that the $I_{\phi'}$
are pairwise disjoint and disjoint from $U \cup Im(\phi^c)$.
This is possible for $I_{\phi'_e}$ by Lemma \ref{exchange}.ii
as $\dim(f_j\phi)=q$ so $Im(\phi^e) \sm e$ 
are pairwise disjoint for all $e \in \phi(Q)$. 

As $G$ is $(\oO,h)$-extendable, the number of such choices
for $\phi^c$ and $\phi^{\phi'}$ given $\phi$ and $\mc{E}$
is at least $0.9\oO n^{q(p-1)+v_+}$, where 
$v_+ = \sum_{\phi'} |I_{\phi'}| = p^r q(p^{q-r}-1)-Q(q-r)$.

Next we specify the algebraic constraints.
We condition on $f_j\phi^c$
such that $Im(\phi^c)$ is $j$-generic for $Im(\phi)$,
which occurs with probability $1-O(n^{-1})$.
We define $w_{\phi'}=f_j\phi^c\phi'$ for $\phi' \in \Ups'$
and note that each $\dim(w_{\phi'})=q$, 
as $Im(\phi^c)$ is $j$-generic for $Im(\phi)$ and
$\dim(f_j\phi)=q$ as $Q' = \phi(Q) \in \mc{Q}_j$.

Then $\phi^c$ will define a cascade
$C_{\phi^c} = \sum \{A^{\phi^c\phi'(Q)}: \phi' \in \Ups'\}$
with each $A^{\phi^c\phi'(Q)}=\phi^{w_{\phi'}}(K^r_q(Ker))
= \phi^{\phi'}(K^r_q(Ker))$ as in Definitions 
\ref{def:absorb} and \ref{def:cascade} if 
\begin{enumerate}
\item
$f_j\phi^{\phi'}((i,a)) = (w_{\phi'})_i + a \cdot w_{\phi'}$ 
for each $\phi' \in \Ups'$, $i \in [q]$, $a \in Ker$, and
\item
$\phi^{\phi'} \phi^a(Q)$ is activated and
$T_e=j$ and $\pi_e \phi^{\phi'} \phi^a=id$ 
for all $\phi' \in \Ups'$, $a \in Ker^r$,
$e \in \phi^{\phi'} \phi^a(Q)$.
\end{enumerate}
Given $f_j\phi^c$, as all cliques are activated 
independently with probability at least $\oO^2$,
these events occur with probability
$(1+O(n^{-1})) (p^{-a})^{v_+} ((z(q)_r)^{-Q} \oO^2 )^{p^{r(q-r+1)}}$,
provided that (i) does not contradict injectivity of $f_j$:
we need to show that $(w_{\phi'})_i+a \cdot w_{\phi'}$ are distinct 
for distinct choices of $(\phi',i,a)$ with
$\phi' \in \Ups'$, $i \in [q]$, $a \in Ker \sm \{0\}$.

Suppose for contradiction that we have some identity
$(w_{\phi^1})_{i^1}+a^1 \cdot (w_{\phi^1}) 
= (w_{\phi^2})_{i^2}+a^2 \cdot (w_{\phi^2})$.
If $\phi^1=\phi^2$ then as $\dim(w_{\phi^1})=q$ we have
$e_{i^1} - e_{i^2} = a^1 - a^2 \in Ker$,
so $i^1=i^2$ by Lemma \ref{Kprops}, so $a^1=a^2$.
If $\phi^1 \ne \phi^2$ then
as $|\{i: a^2_i \ne 0\}|>r$ by Lemma \ref{Kprops}
we can find $i \in [q]$ with $a^2_i \ne 0$
and $\phi^2(i) \notin Im(\phi^1)$.
Then $(w_{\phi^2})_i=f_j\phi^c\phi^2(i)$
appears with a nonzero coefficient in the identity
but is not in the span of $w_{\phi^1}$ and
the other coordinates of $w_{\phi^2}$,
which gives the required contradiction.

We deduce (using $\oO<\oO_0$)
that the number $X$ of cascades for $Q'$ satisfies 
\[\mb{E}[X \mid \mc{E} ] > 0.9\oO n^{q(p-1)+v_+} \cdot
(1+O(n^{-1})) (p^{-a})^{v_+} ((z(q)_r)^{-Q} \oO^2)^{p^{r(q-r+1)}}
> 3 \oO^{p^{q^2}} n^{q(p-1)}.\]

For concentration of $X \mid \mc{E}$
we will apply Lemma \ref{lip3} 
similarly to Subsection \ref{sec:nibble}.
We start by showing concentration of 
$\mb{E}[X \mid \mc{E},f_j]$, which is the conditional
expectation where we reveal the embedding $f_j$
but not the other random choices 
in the construction of the template.
We claim that changing any $f_j(x)$
with $x \notin U$ from $\aA$ to $\aA'$ 
affects $X$ by $O(n^{q(p-1)-1})$.
To see this, it suffices to show for any 
$\phi' \in \Ups'$, $i \in [q]$, $a \in Ker$,
writing $w=w_{\phi'}$, that if the constraint
$w_i + a \cdot w = \aA$ can affect $X$
then it defines a strict subspace
of the cascade variables 
$f_j \phi^c \sm f_j \phi
\in \mb{F}_{p^a}^{[q] \times [2,p]}$.
We need to show that
the constraint is non-trivial,
in that it has a non-zero
coefficient of some variable 
in $f_j \phi^c \sm f_j \phi$.
If  $a=0$ then $w_i$ is such a variable, 
as the constraint can only affect $X$
when $\phi'(i) \ne i$. 
If $a \ne 0$ then Lemma \ref{Kprops}.iii implies
that if the constraint only depends on $f_j\phi$
then it has the form $f_j(x)=\aA$ for some $x \in U$, 
but we have already conditioned on this.
This proves the claim. 
Thus $\mb{E}[X \mid \mc{E},f_j]$
is $O(n^{2q(p-1)-1})$-varying (in the choice of $f_j$),
so w.h.p.\ $\mb{E}[X \mid \mc{E},f_j]
> 2\oO^{p^{q^2}} n^{q(p-1)}$
by Lemma \ref{lip3}.
 
Now we fix $f_j$ with $\mb{E}[X \mid \mc{E},f_j]
> 2\oO^{p^{q^2}} n^{q(p-1)}$ and show concentration
of $X$ under the remaining random choices.
Consider any clique $\phi^*(Q) \in K^r_q(G)$ such that
there is some $y \in \mb{F}_{p^a}^r$ 
with $f_j(\phi^*(i)) = (My)_i$ for all $i \in [q]$.
The random choices of any $T_e$ or $\pi_e$ 
for $e \in \phi^*(Q)$ or whether $\phi^*(Q)$ is activated
can affect whether $\phi^*(Q)$ is in the template,
and so all cascades containing $\phi^*(Q)$.
We consider the separately the effect according 
to the role of $\phi^*(Q)$ in the cascade:\
we fix $\phi' \in \Ups'$, $I \in Q$, $a \in Ker^I$,
where we can assume $r_{a\phi'}>0$
by Remark \ref{rem:nondegen},
and consider cascades where $\phi^*(Q)$ is the template
clique containing $e^w_a$, where $w=w_{\phi'}$, 
that is, $f_j \phi^* = M M_I^{-1} (e_I + a) w$.
By Lemma \ref{lem:ee'}, there are at most $(p^a)^{r_{a\phi'}}$
choices of $e^w_a$, which determines $\phi^*(Q)$,
and by Lemma \ref{cascade:role} there are at most
$(p^a)^{q(p-1)-r_{a\phi'}}$ cascades 
with $\phi^*(Q)$ in this role.
Summing over all $\phi'$ and $a$ we see that
$X \mid \mc{E}, f_j$ is $O(n^{2q(p-1)-1})$-varying,
so by Lemma \ref{lip3} w.h.p.\ 
$X > \oO^{p^{q^2}} n^{q(p-1)}$. \qed

\subsection{Cascading extensions}

The Clique Exchange Algorithm in Section \ref{sec:cea}
will take an integral decomposition $\Phi$ of the spill $S$ 
and modify it into a signed decomposition 
in which all positive cliques are cascading.
Here we establish some lower bounds on extensions 
required for the analysis of the algorithm.
Throughout we use the same notation as in Lemma \ref{exchange}.

To motivate the bounds in this subsection 
we start with an informal description of the algorithm.
We will repeatedly use the clique exchange tool from Lemma \ref{exchange}
to replace some signed cliques by another signed combination
of cliques while preserving the property $\pl \Phi = S$.
The goal is to ensure that all positive cliques are cascading.
In the final Elimination Phase of the algorithm 
we will eliminate certain `cancelling pairs',
which consist of two cliques of opposite sign in $\Phi$
sharing one common edge, without introducing
any other clique containing that edge.
This allows us to eliminate high multiplicity 
uses of any edge, and also uses of non-edges $K^r_n - G$.
We will also ensure that all new positive cliques are cascading.

\begin{figure}
\begin{center}
\includegraphics[width=0.8\textwidth]{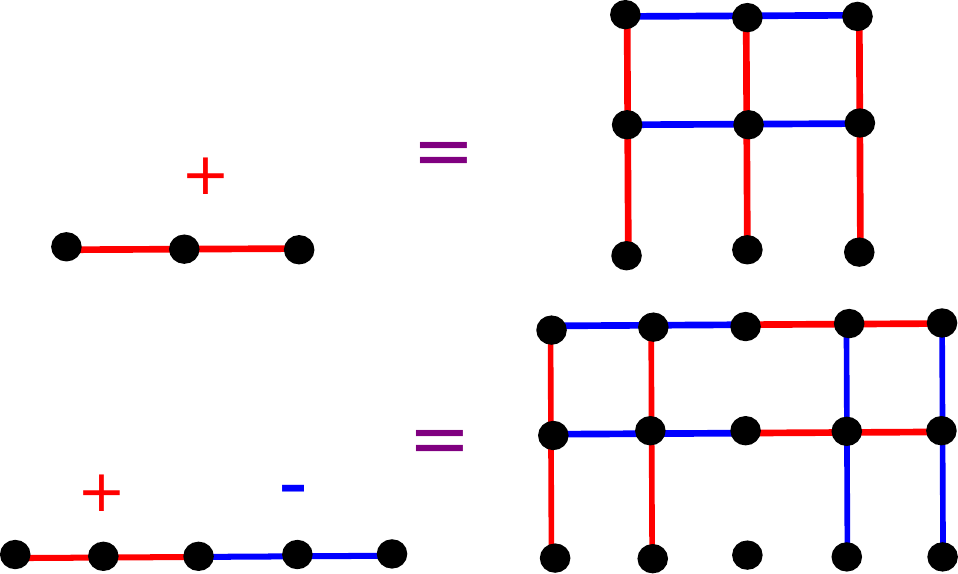}
\caption{Splitting and Elimination `cartoons' for $3$-graph matchings.}
\end{center}
\label{fig:elim}
\end{figure}

Two preparatory phases are required before the Elimination Phase.
The first Splitting Phase will address the issue that a given signed
element of $\Phi$ may be required for more than one cancelling pair.
In this phase we replace $\Phi$ by $\Phi'$, 
preserving $\pl \Phi' = \pl \Phi = S$, 
so that we can choose the cancelling pairs of cliques 
required for the Elimination Phase,
and any signed element of $\Phi'$ is in at most one such pair.
We will also ensure that all new cliques are rainbow,
which facilitates the later phases.
The second Solo Phase will replace the set of cliques
not in cancelling pairs by an equivalent set such that
all positive cliques are cascading. Then the Elimination Phase
will remove all other cliques and replace them by new cliques
that achieve the goal of making all positive cliques cascading.
 
We start with some notation pertaining to extensions $E(\phi)$
that correspond to choosing a copy of $\OO=K^r_q(p)$
containing some base clique $\phi(Q)$,
which will be a positive clique in $\Phi$.
Recall that $\OO$ is equipped with two clique decompositions
$\Ups,\Ups' \sub [q](p)$, where $id_{[q]} \in \Ups$.
Given any extension $\phi^+$ of $\phi$ to an embedding of $\OO$,
we can modify the decomposition $\Phi$, while preserving $\pl \Phi = S$,
by adding new cliques $\phi^+ \phi'(Q)$ that are positive
for all $\phi' \in \Ups'$ and negative for all $\phi' \in \Ups$.
Note that one of the negative cliques is $\phi' id_{[q]}(Q)=\phi(Q)$,
which we can remove together with its original positive copy in $\Phi$. 
Thus we have removed $\phi(Q)$ and added new cliques (positive and negative),
all of which are edge-disjoint from $\phi(Q)$, except that for each $e \in Q$
there is a new positive clique that intersects $\phi(Q)$ in $\phi(e)$.
This is illustrated by the top row of Figure \ref{fig:elim}.

We introduce two variants of the above extension
in which we specify desirable properties 
of the new positive cliques $\phi^+ \phi'(Q)$ with $\phi' \in \Ups$.
In the first variant, we use rainbow extensions 
(each edge is compatible with its 
own private embedding $f_j$ with $j \in [z]$),
and also require compatibility with edge orders
(indicated in the notation by the superscript $\pi$).
In the second variant, we require that all positive cliques are cascading
(indicated by the superscript $c$), in a rainbow manner,
meaning that each of these cascading cliques
is compatible with its own private embedding $f_j$.

\begin{defn} \label{def:Xpi}
Consider any extension $E(\phi)=(\phi,[q],\OO)$ 
where $\OO=K^r_q(p)$ and $\phi(Q) \in K^r_q(K^r_n)$.

We let $X^\pi_{E(\phi)}(G^*)$ be the set or number of rainbow
extensions $\phi^+ \in X'_{E(\phi)}(G^*)$ 
(recall Remark \ref{rem:X'}) such that each
$\phi^+ \phi'(Q)$ with $\phi' \in \Ups'$ 
is $\pi$-compatible, meaning that we can order it
as $\phi^+ \phi'(Q) = \psi(Q)$ for some $\psi$
such that $\pi_e \psi = id$ for all $e \in \psi(Q)$
(recall $\pi_e$ was defined for all $e \in K^r_n$).

We let $X^c_{E(\phi)}(G^*)$ be the set or number 
of $\phi^+ \in X_{E(\phi)}(G^*)$ 
that are `rainbow $\Ups'$ cascading', 
meaning that $\phi^+ \psi(Q) \in \mc{Q}^*$ 
(recall Definition \ref{def:cascading}) for all $\psi \in \Ups'$,
and $j \ne j'$ whenever $\{\psi,\psi'\} \sub \Ups'$ with
$\phi^+ \psi(Q) \in \mc{Q}_j$, $\phi^+ \psi'(Q) \in \mc{Q}_{j'}$.
\end{defn}

By Lemma \ref{ext*} and Remark \ref{rem:X'} w.h.p.\ 
\begin{equation} \label{Xpi}
X^\pi_{E(\phi)}(G^*) > \oO n^{pq-q} (z\rho/2(q)_r)^{Qp^r};
\end{equation}
this estimate will be used in the Splitting Phase of the algorithm,
as it shows that there are many ways to replace any positive clique 
by an equivalent set of rainbow positive and negative cliques.
The next estimate will be used in the Solo Phase of the algorithm,
as it shows that there are many ways to replace any positive clique 
(these are now all rainbow and $\pi$-compatible)
by an equivalent set of rainbow positive and negative cliques
in which all positive cliques are cascading.

\begin{lemma} \label{split:cascade}
If $\phi(Q) \sub G^*$ is rainbow and $\phi$ is $\pi$-compatible
then w.h.p.\ \[X^c_{E(\phi)}(G^*) > 0.9 \oO^{3Q^2 p^r} n^{pq-q}.\]
\end{lemma}

The proof of Lemma \ref{split:cascade}
requires the following analogue of Lemma \ref{e|E}.

\begin{lemma} \label{Q*|eE}
Let $S \sub G$ with $|S|<h$
and $\mc{E}^S = \cap_{f \in S} \mc{E}^f$ be an atom of $\mc{F}^S$.
Suppose $\phi(Q) \sub G$
has at most one edge in $S$
and the edges of $\phi(Q) \sm S$
are not touched by $\mc{E}^S$.
Let $j \in [z]$ be such that the set $R$ 
of $v \in \phi([q])$ such that $f_j(v)$ is revealed 
by $\mc{E}^S$ has $|R| \le r$ and $\dim f_j(R)=|R|$.
If $|R|=r$ suppose also that $R=e \in G^*_j$
with $\pi_e\phi=id$ and $\phi(Q) \cap S = \{e\}$.
Then $\mb{P}( \phi(Q) \in \mc{Q}_j \mid \mc{E}^S ) > \oO^{3Q^2}$.
\end{lemma}

\nib{Proof.}
Let $1_e$ be $1$ if $\phi(Q) \cap S = \{e\}$ or $0$ otherwise.
For $e' \in \phi(Q) \sm S$ let $\pi^0_{e'}:e' \to [q]$
be such that $\pi^0_{e'} \phi = id$.
For each $e' \in \phi(Q) \sm S$ we fix 
$\phi^{e'}(Q) \in K^r_q(G)$ with $\pi^0_{e'}\phi^{e'}=id$
and estimate the probability that all such
$e' \in G^*_j$ with $M^*(e') = \phi^{e'}(Q)$.  
As $G$ is $(\oO,h)$-extendable, there are at least
$(1-O(n^{-1})) \oO n^{(Q-1_e)(q-r)}$ choices for all $\phi^{e'}$
such that the sets $Im(\phi^{e'}) \sm e'$ are pairwise
disjoint and disjoint from $Im(\phi)$, and no edge
of any $\phi^{e'}(Q)$ is touched by $\mc{E}^S$.
The probability that $\phi^{e'}(Q)$ is activated,
$T_f=j$ and $\pi_f\phi^{e'}=id$ for all such $e'$
and $f \in \phi^{e'}(Q)$ is at least
$((z(q)_r)^{-Q} \oO^2 )^{Q-1_e}$.

As $M$ is generic, when we reveal $f_j\!\mid_{Im(\phi) \sm R}$
then each $e' \in \phi(Q) \sm S$ will have
a unique $y^{e'} \in \mb{F}_{p^a}^r$ such that 
$(My^{e'})_i = f_j\pi_{e'}^{-1}(i)$ for all $i \in Im(\pi_{e'})$.
We claim that with probability $1-O(n^{-1})$
the choice of $f_j\!\mid_{Im(\phi) \sm R}$ 
is such that its entries and all $(My^{e'})_i$ 
with $e' \in \phi(Q) \sm S$ 
and $i \in [q] \sm \pi_{e'}(e')$
are distinct and not equal to 
any $f_j(x)$ revealed by $\mc{E}^S$.
To see this, first note that there are at most $hr$ 
such revealed $f_j(x)$, and each $(My^{e'})_i$
is a non-trivial linear function of the entries of
$f_j\!\mid_{Im(\phi) \sm R}$, so is a revealed
value with probability $O(1/n)$. Also, with probability
$1-O(n^{-1})$ we have $\dim(f_j\phi)=q$,  
which implies the distinctness property,
so the claim holds.

With probability $(1+O(n^{-1}))(p^{-a})^{(q-r)(Q-1_e)}$
we have $f_j(\phi^{e'}(i))=(My^{e'})_i$ 
for all such $e'$ and $i \in [q] \sm Im(\pi_{e'})$.
Therefore $\mb{P}(\cap_{e'} \{ M^*(e') = \phi^{e'}(Q) \} \mid \mc{E}^S)
> (1+O(n^{-1})) [ (z(q)_r \oO^{-2} )^Q p^{a(q-r)} ]^{1_e-Q}$.
Summing over all choices for $\phi^{e'}$ gives 
$\mb{P}( \phi(Q) \in \mc{Q}_j \mid \mc{E}^S ) 
> \oO^{3Q^2}$. \qed

\begin{rem} \label{rem:Q*|eE}
Note that the condition on $j$ 
in Lemma \ref{Q*|eE} holds trivially if 
$T_{e'} \ne j$ for all $e' \in S \sm \phi(Q)$,
which will be the case for our application 
to the proof of Lemma \ref{split:cascade};
we allow the more general condition on $j$
in the statement for use 
in the proof of Lemma \ref{elim:cascade}.
\end{rem}

\nib{Proof of Lemma \ref{split:cascade}.}
As $G$ is $(\oO,h)$-extendable, 
there are at least $\oO n^{pq-q}$ choices 
of $\phi^+ \in X_{E(\phi)}(G)$. We fix any such $\phi^+$
and estimate $\mb{P}(\phi^+ \in X^c_{E(\phi)}(G^*))$
by repeated application of Lemma \ref{Q*|eE}.
We consider sequentially each $\psi \in \Ups'$, and fix 
$j \in [z]$ distinct from all previous choices 
such that if $Im(\psi) \cap Im(\phi) = e$ then $e \in G^*_j$
(this is possible by the hypotheses of the lemma).

We let $\mc{E}$ be the intersection of all local events $\mc{E}^e$ 
where $e \sub Im(\phi)$ or $e \sub Im(\phi^+\psi')$ 
for some previously considered $\psi' \in \Ups'$.
If any edge of $\phi^+\psi(Q)$ is touched by $\mc{E}$ 
we discard $\phi^+$; thus we discard $O(n^{pq-q-1})$ choices.
Then $\mb{P}( \phi^+\psi(Q) \in \mc{Q}_j \mid \mc{E} ) 
> \oO^{3Q^2}$ by Lemma \ref{Q*|eE}.
Multiplying all conditional probabilities and summing over $\phi^+$
gives $\mb{E}X^c_{E(\phi)}(G^*) 
> (1-O(n^{-1}))(\oO^{3Q^2})^{p^r} n^{pq-q}$.

Similarly to several earlier proofs of concentration,
$X^c_{E(\phi)}(G^*) \mid \mc{E}$ is $O(n^{2(pq-q)-1})$-varying,
so by Lemma \ref{lip3} w.h.p.\ 
$X^c_{E(\phi)}(G^*) > 0.9 \oO^{3Q^2 p^r} n^{pq-q}$. \qed

\medskip

Next we describe the construction used in
the Elimination Phase of the algorithm.
We combine two copies of our $\OO = K^r_q(p)$,
each equipped as above with clique decompositions $\Ups$, $\Ups'$,
forming a new $r$-graph $\OO^*$ with clique decompositions $\Ups^\pm$,
containing two specified cliques of opposing sign overlapping 
in a single edge and inducing a subgraph that contains no other edges.

\begin{defn} \label{def:elim}
Let $\OO_1$ and $\OO_2$ be two copies of $\OO$.
Fix $f \in \Ups$ and $f' \in \Ups'$ 
with $|V(f) \cap V(f')|=r$.
For $j=1,2$ we denote the copies of
$\Ups$, $\Ups'$, $f$, $f'$ in $\OO_j$
by $\Ups_j$, $\Ups'_j$, $f_j$, $f'_j$.
Let $\OO^*$ be obtained by identifying
$\OO_1$ and $\OO_2$ so that $f'_1=f'_2$.
Let $\Ups^+ = \Ups_1 \cup (\Ups'_2\sm\{f'_2\})$
and $\Ups^- = \Ups_2 \cup (\Ups'_1\sm\{f'_1\})$.
Then $\Ups^+$ is a $K^r_q$-decomposition of $\OO^*$ containing $f_1$
and $\Ups^-$ is a $K^r_q$-decomposition of $\OO^*$ containing $f_2$.
\end{defn}

Note that $\OO^*[Im(f_1) \cup Im(f_2)] = f_1(Q) \cup f_2(Q)$
as stated before the definition: every edge of $\OO^*$ contained in 
$Im(f_1) \cup Im(f_2)$ is contained in $f_1(Q)$ or $f_2(Q)$.

Next we introduce some notation, 
analogous to that in Definition \ref{def:Xpi},
but here pertaining to extensions $E(\phi^\pm)$
that correspond to choosing a copy of $\OO^*$
containing some base pair of cliques $\phi^\pm(Q)$,
which will be a cancelling pair $\phi^\pm(Q)$ in $\Phi$.
This cancelling pair will correspond to the two specified cliques 
$f_1,f_2$ of opposite sign in the decompositions $\Ups^\pm$ of $\OO^*$.

We can use any extension $\phi^+$ to an embedding of $\OO^*$
to modify the decomposition $\Phi$ while preserving $\pl \Phi = S$:
we add new cliques $\phi^+ \phi'(Q)$ that are positive
for all $\phi' \in \Ups^-$ and negative for all $\phi' \in \Ups^+$.
The new cliques corresponding to $f_1,f_2$ will have opposite signs
to those in $\phi^\pm(Q)$, so we can remove these cliques
together with $\phi^\pm(Q)$, with the net effect
of removing the cancelling pair $\phi^\pm(Q)$ and replacing it by new cliques 
(positive and negative), all of which are edge-disjoint from $\phi(Q)$, 
except that for each $e$ in exactly one of $\phi^\pm(Q)$ there is
a new clique of the same sign that intersects 
$Im(\phi^+) \cup Im(\phi^-)$ in $\phi(e)$.
This is illustrated by the bottom row of Figure \ref{fig:elim}.

We require the variant of this extension with the same desirable
property as in the second variant of Definition \ref{def:Xpi},
namely that all new positive cliques are cascading,
again in a rainbow manner,
meaning that each of these cascading cliques
is compatible with its own private embedding $f_j$.

\begin{defn} \label{def:elim2}
Let $Q^\pm = \phi^\pm(Q)$ be cliques with 
$Q^+ \cap Q^- = \{e\}$ and $\phi^\pm$ are $\pi$-compatible.
We label $\OO^*$ so that $f_1=f_2=[q]$ 
and $f_1 \cap f_2 = Im(\pi_e)$ 
consistently with both copies of $[q]$.
Consider the extension
$E(\phi^\pm) := (\phi^\cup,f_1 \cup f_2,\OO^*)$ where
$\phi^\cup f_1 = \phi^+$ and $\phi^\cup f_2 = \phi^-$.

We let $X^c_{E(\phi^\pm)}(G^*)$ be the set or number 
of $\phi^* \in X_{E(\phi^\pm)}(G^*)$ 
that are `rainbow $\Ups^- \sm \{f_2\}$ cascading', meaning that
$\phi^* \psi(Q) \in \mc{Q}^*$ for all $\psi \in \Ups^- \sm \{f_2\}$
and $j \ne j'$ whenever $\{\psi,\psi'\} \sub \Ups^- \sm \{f_2\}$ with
$\phi^* \psi(Q) \in \mc{Q}_j$, $\phi^* \psi'(Q) \in \mc{Q}_{j'}$.
\end{defn}

\begin{lemma} \label{elim:cascade}
Suppose $Q^\pm = \phi^\pm(Q)$ with $Q^+ \cap Q^- = \{e\}$,
where $\phi^\pm$ are $\pi$-compatible 
and $Q^\pm \sm \{e\}$ are both rainbow in $G^*$.
Then w.h.p.\ $X^c_{E(\phi^\pm)}(G^*) 
> 0.9 \oO^{6Q^2 p^r} n^{2pq-3q+r}$.
\end{lemma}

\nib{Proof.}
The proof is very similar to that of Lemma \ref{split:cascade}.
As $G$ is $(\oO,h)$-extendable, 
there are at least $\oO n^{2pq-3q+r}$ choices 
of $\phi^* \in X_{E(\phi^\pm)}(G)$. We fix any such $\phi^*$
and estimate $\mb{P}(\phi^* \in X^c_{E(\phi^\pm)}(G^*))$
by repeated application of Lemma \ref{Q*|eE}.
We consider sequentially each $\psi \in \Ups^- \sm \{f_2\}$;
by the hypotheses of the lemma we can fix
$j \in [z]$ distinct from all previous choices,
except that if $Im(\psi) \cap Im(\phi') = e'$
with $\phi' \in \{\phi^+,\phi'\}$ 
then we choose $j$ so that $e \in G^*_j$ and we allow
this $j$ to occur once for both of $\phi^\pm$.

We let $\mc{E}$ be the intersection of all local events $\mc{E}^{e'}$ 
where $e' \sub Im(\phi^\pm)$ or $e' \sub Im(\phi^*\psi')$ 
for some previously considered $\psi' \in \Ups^-\sm \{f_2\}$.
If any edge of $\phi^*\psi(Q)$ is touched by $\mc{E}$ 
we discard $\phi^*$; thus we discard $O(n^{2pq-3q+r-1})$ choices.
Then $\mb{P}( \phi^*\psi(Q) \in \mc{Q}_j \mid \mc{E} ) 
> \oO^{3Q^2}$ by Lemma \ref{Q*|eE}.
Multiplying all conditional probabilities and summing over $\phi^+$
gives $\mb{E}X^c_{E(\phi)}(G^*) 
> (1-O(n^{-1}))(\oO^{3Q^2})^{2p^r} n^{2pq-3q+r}$.
The lemma now follows by a concentration argument
very similar to those given elsewhere. \qed

\section{Clique Exchange Algorithm} \label{sec:cea}

This section contains the proof of the following lemma,
which takes an integral decomposition $\Phi$ of the spill $S$ 
obtained via the results of Section \ref{sec:int} 
and modifies it into a signed decomposition 
in which all positive cliques are cascading.
This is the main remaining step in the proof
of Theorem \ref{main+}, which will then follow
quite easily in the next section. 
Indeed, after Lemma \ref{hole} we will have
$G = \sum M^* + \sum M' - \sum M^+$,
where $M'= M^n \cup M^c \cup M^-$.
It will then suffice to find edge-disjoint
cascades for each clique in $M^+$,
as then flipping these and removing $M^+$
will give a decomposition of $G$.

\begin{lemma} \label{hole}
Suppose $S \sub G^*$ is $K^r_q$-divisible
and $M^*(S)$ is a set that is 
$c_2$-bounded and linearly $c_2$-bounded.
Then there are $M^\pm \sub K^r_q(G^*)$
such that every clique in $M^+$ is cascading,
$M^*(\sum M^+)$ is a set 
(with no multiple edges)
and linearly $3c_4$-bounded,
and $\sum M^+$ is the disjoint union
of $\sum M^-$ and $S$.
\end{lemma}

We will now start the proof of Lemma \ref{hole},
which will occupy the remainder of this section.

We first apply Lemma \ref{bddint} (using $c_2>n^{-1/4Qb}$)
to obtain $\Phi \in \mb{Z}^{K^r_q(K^r_n)}$ 
such that $\pl \Phi = S$ and $\pl^\pm \Phi$ 
are $N^2 c_2$-bounded, where $N=(2q)^q !$.
We will modify $\Phi$, maintaining $\pl \Phi = S$,
with the goal of making all positive cliques cascading,
using the Clique Exchange Algorithm,
which was informally described in the previous section.
Below we will formally describe and analyse the three phases
of the algorithm, namely the Splitting Phase,
the Solo Phase and the Elimination Phase.
Each uses random greedy algorithms
(similarly to the proof of Lemma \ref{cover})
that maintain edge-disjointness,
and moreover disjointness of all $M^*(Q^+)$
where $Q^+$ is a positive clique,
so that it will be possible to choose
edge-disjoint cascades for all positive cliques.

Throughout this section we fix $\OO$,
$\Ups$ and $\Ups'$ as in Lemma \ref{exchange},
identify $\Ups$ and $\Ups'$ with subsets of $[q](p)$,
identify $[q]$ with $\{ (i,1): i \in [q] \} \sub V(\OO)$,
and assume $[q] \in \Ups$.
We write $\OO' = \OO \sm Q = K^r_q(p) \sm \tbinom{[q]}{r}$.
Now we define the first phase
of the Clique Exchange Algorithm.
(Recall $E(\phi)=(\phi,[q],\OO)$ 
and $X^\pi_{E(\phi)}(G^*)$
from Definition \ref{def:Xpi}.)

\begin{alg} (Splitting Phase)
Let $(Q^i=\phi_i(Q): i \in [|\Phi|])$ be any ordering 
of the signed elements of $\Phi$.
We apply a random greedy algorithm
to choose $\phi^*_i \in X^\pi_{E(\phi_i)}(G^*)$.
Write $A_i = \cup_{i'<i} M^*(\phi^*_{i'}(\OO'))$.
We choose $\phi^*_i \in X^\pi_{E(\phi_i)}(G^*)$
uniformly at random subject to
$\phi^*_i(\OO') \cap (M^*(S) \cup A_i) = \es$.
(If there is no such choice of $\phi^*_i$
then the algorithm aborts.)
\end{alg}

Note that if $\phi^*_i \in X^\pi_{E(\phi_i)}(G^*)$
then $\phi^*_i(\OO') \sub G^*$ is rainbow,
so $M^*(\phi^*_i(\OO'))$ is a set.
The proof of the following lemma 
is very similar to that of Lemma \ref{cover}.

\begin{lemma} \label{split}
W.h.p.\ the Splitting Phase does not abort
and $A_{|\Phi|}$ is $c_3$-bounded
and linearly $c_3$-bounded
(recall $c_3 = \oO^{-h/20Q} c_2$). 
\end{lemma}

\nib{Proof.} 
For $i \in [|\Phi|]$ we let $\mc{B}_i$ be the bad event 
that $A_i$ is not $c_3$-bounded
or not  linearly $c_3$-bounded. Let $\tau$ be the 
smallest $i$ for which $\mc{B}_i$ holds or the algorithm 
aborts, or $\infty$ if there is no such $i$. 
It suffices to show w.h.p.\ $\tau=\infty$. 
We fix $i_0 \in [|\Phi|]$ and bound 
$\mb{P}(\tau=i_0)$ as follows. 

We claim that for any $i<i_0$
the conditions on $\phi^*_i$ forbid at most 
half of the possible choices of $\phi^*_i$.
To see this, recall from (\ref{Xpi}) that
$X^\pi_{E(\phi_i)}(G^*) > \oO n^{pq-q} (z\rho/2(q)_r)^{Qp^r}$. 
As $A_i$ is $c_3$-bounded and $M^*(S)$ is $c_2$-bounded, 
at most $2|\OO'| c_3 n^{pq-q}$ choices of $\phi^*_i$ are forbidden.
The claim follows.

For each $e \in G^*$ let 
$r_e = \sum_{i<i_0} \mb{P}'(e \in M^*(\phi^*_i(\OO')))$,
where $\mb{P}'$ denotes conditional probability
given the choices made before step $i$.
Note that $e \in M^*(\phi^*_i(\OO'))$ if and only if
$e' \in \phi^*_i(\OO')$ for some $e' \in M^*(e)$.
For fixed $i$ and $e' \in M^*(e)$, 
writing $r'=|e' \sm V(Q^i)|$, there are 
at most $r!|\OO'|n^{pq-q-r'}$ choices of $\phi^*_i$
such that $e' \in \phi^*_i(\OO')$, so by the claim 
$\mb{P}'(e' \in \phi^*_i(\OO')) 
< 2r!|\OO'|\oO^{-1} (z\rho/2(q)_r)^{-Qp^r} n^{-r'}$.
Also, given $r' \in [r]$, as $\pl^+ \Phi$ 
is $N^2 c_2$-bounded there are at most
$N^2 c_2 \tbinom{r}{r'} n^{r'}$ choices of $i$
such that $|e' \sm V(Q^i)|=r'$. Therefore 
$r_e < q^r|\OO'|\oO^{-1} (z\rho/2(q)_r)^{-Qp^r} 2^{r+1} N^2 c_2$.

Now fix any $f \in \tbinom{[n]}{r-1}$ and let $X =  \sum_{i<i_0} X_i$,
where $X_i = |\{e: f \sub e \in M^*(\phi^*_i(\OO'))\}|$.
Then each $|X_i|<Q|\OO'|$ and $\sum_{i<i_0} \mb{E}'X_i
= \sum \{r_e: f \sub e \in G^*\}$,
so by Lemma \ref{dom} w.h.p.\ $X < c_3 n$.

Similarly, consider any $j \in [z]$, $I \in Q$ and $I$-line $L$,
and let $X' = \sum_{i<i_0} X'_i$, where 
$X'_i = |\{e \in M^*(\phi^*_i(\OO')): f_j(e) \in L\}|$.
Then each $|X'_i|<Q|\OO'|$ and $\sum_{i<i_0} \mb{E}'X_i
= \sum \{r_e: f_j(e) \in L \}$,
so by Lemma \ref{dom} w.h.p.\ $X < c_3 p^a$.
 Thus w.h.p.\ $A_i$ is $c_3$-bounded 
and linearly $c_3$-bounded for all $i<i_0$, 
so w.h.p.\ $\tau=\infty$, as required. \qed

\medskip

We let $\Phi' = \Phi + \sum_{i \in [|\Phi|]} s(Q^i)
( \phi^*_i(\Ups')-\phi^*_i(\Ups) )$.
Then $\pl \Phi' = \pl \Phi = S$, and all signed elements
of $\Phi$ are cancelled, so $\Phi'$ is supported 
on cliques $Q'$ added during the Splitting Phase,
all of which are $\pi$-compatible,
have all but at most one edge in $G^*$,
and are rainbow.

\medskip

We classify cliques added during the Splitting Phase as near or far, 
where near cliques are those of the form $\phi^*_i \phi(Q)$ 
for $\phi \in \Ups'$ with $|Im(\phi) \cap [q]|=r$.
Also, for each pair $(e,Q')$ 
where $Q'$ is added during the Splitting Phase and $e \in Q'$,
we call $(e,Q')$ near if $Q'=\phi^*_i \phi(Q)$ is near 
and $e = \phi^*_i(Im(\phi) \cap [q])$, 
otherwise we call $(e,Q')$ far.
We also classify cliques and near pairs as positive 
or negative according to their sign in $\Phi'$.

Note that for each edge $e$ such that there is some
far pair $(e,Q')$ there are exactly two such far pairs
and they have opposite sign in $\Phi'$.
The $r$-graph $\GG \sub G^*$ of all such $e$ satisfies 
$\GG = \cup_{i \in [|\Phi|]} \phi^*_i(\OO')$
and $M^*(\GG)=A_{|\Phi|}$, which is disjoint from $M^*(S)$.

For each edge $e$, there are 
$\pl^+ \Phi_e$ positive near pairs $(e,Q')$ 
and $\pl^- \Phi_e$ negative near pairs $(e,Q')$.
We group the near pairs on $e$ into `cancelling' pairs,
each consisting of one positive and one negative near pair,
and one additional positive near pair $(e,Q^e)$ if $e \in S$,
which we call `solo',
where $Q^e \in K^r_q(G^*)$ as $S \sub G^*$.
Note that each cancelling pair on $e$
intersects only in $e$ by Lemma \ref{exchange}.i.
In a cancelling pair $\{(e,Q^+),(e,Q^-)\}$
the common edge may be any $e \in K^r_n$, 
but every other edge of $Q^+ \cup Q^-$ is in $G^*$.
 
\medskip

Let $(\{(e_i,Q_i^+),(e_i,Q_i^-): i \in [P''])$ 
be any ordering of the cancelling pairs,
where $Q_i^+$ is positive and $Q_i^-$ is negative.
Let $(Q^i=\phi'_i(Q): i \in [P'])$ be the 
positive far cliques and
cliques in solo near pairs $(e,Q^e)$.
By definition of $X^\pi_{E(\phi_i)}(G^*)$,
we can choose $\pi$-compatible orderings $\phi'_i$.
Furthermore, each $Q^i \in K^r_q(G^*)$ is rainbow,
so we can apply Lemma \ref{split:cascade}.
We will now process the solo pairs and positive far cliques
in the second phase of the algorithm,
which is the same as the Splitting Phase,
but now we ensure that all positive cliques are cascading
(recall $X^c_{E(\phi)}(G^*)$ from Definition \ref{def:Xpi}).

\begin{alg} (Solo Phase)
We apply a random greedy algorithm to choose 
$\phi^*_i \in X^c_{E(\phi'_i)}(G^*)$ for each $i \in [P']$.
Write $A'_i = \cup_{i'<i} M^*(\phi^*_{i'}(\OO'))$. 
We choose $\phi^*_i \in X^c_{E(\phi'_i)}(G^*)$
uniformly at random such that 
$\phi^*_i(\OO')$ is edge-disjoint
from $A'_i \cup A_{|\Phi|} \cup M^*(S)$.
\end{alg}

We note in each step of the Solo Phase
that $M^*(\phi^*_i(\OO'))$ is edge-disjoint 
from $A'_i \cup A_{|\Phi|} \cup M^*(S)$,
as the latter is a union of cliques in $M^*$.

\begin{lemma} \label{absorbI}
W.h.p.\ the Solo Phase does not abort,
and $A'_{P'}$ is $c_4$-bounded
and linearly $c_4$-bounded
(recall $c_4 = \oO^{-h/20Q} c_3$). 
\end{lemma}

\nib{Proof.}
For $i \in [P']$ we let $\mc{B}_i$ be the bad event 
that $A'_i$ is not $c_4$-bounded
or not linearly $c_4$-bounded. Let $\tau$ be the 
smallest $i$ for which $\mc{B}_i$ holds or the algorithm 
aborts, or $\infty$ if there is no such $i$. 
It suffices to show w.h.p.\ $\tau=\infty$. 
We fix $i_0 \in [P']$ and bound 
$\mb{P}(\tau=i_0)$ as follows. 

By Lemma \ref{split:cascade} we have
$X^c_{E(\phi'_i)}(G^*) > 0.9 \oO^{3Q^2 p^r} n^{pq-q}$.
At most half of the choices of 
$\phi^*_i \in X^c_{E(\phi'_i)}(G^*)$
are forbidden, as $M^*(S)$ is $c_2$-bounded,
$A'_i$ is $c_4$-bounded
and $A_{|\Phi|}$ is $c_3$-bounded.

For each $e \in G^*$ let 
$r_e = \sum_{i<i_0} \mb{P}'(e \in M^*(\phi^*_i(\OO')))
= \sum_{i<i_0} \sum_{e' \in M^*(e)}
\mb{P}'(e' \in \phi^*_i(\OO'))$.
Given $r' \in [r]$,
as $A_{|\Phi|} \cup M^*(S)$ is $2c_3$-bounded
there are at most
$2\tbinom{r}{r'}c_3 n^{r'}$ choices of $i$
such that $|e' \sm V(Q^i)|=r'$.
For each such $i$ we have
$\mb{P}'(e' \in \phi^*_i(\OO')) 
< 2r! |\OO'| \oO^{-3Q^2 p^r} n^{-r'}$, so
$r_e < q^r |\OO'| \oO^{-3Q^2 p^r} 2^{r+2} c_3$.
As in the proof of Lemma \ref{split},
by Lemma \ref{dom} w.h.p.\ $A'_i$ is 
$c_4$-bounded and linearly $c_4$-bounded
for all $i<i_0$, so w.h.p.\ $\tau=\infty$, as required. \qed

\medskip

We let $\Phi'' = \Phi' + \sum_{i \in [|P'|]} s(Q^i)
( \phi^*_i(\Ups')-\phi^*_i(\Ups) )$.
Then $\pl \Phi'' = S$, all solo near pairs are cancelled,
and for each positive clique $Q'$ added during Solo Phase
we have $Q'$ cascading, $M^*(Q')$ is a set
and all such $M^*(Q')$ are disjoint.
These positive cliques include all $M^*(e')$ with 
$e' \in Q^e \sm \{e\}$ for some solo near pair $(e,Q^e)$;
we note that $M^*(e') \sub A_{|\Phi|}$ and $M^*(e) \sub M^*(S)$.
All other positive cliques added during Solo Phase
are contained in $A'_{P'}$, which is
$c_4$-bounded, linearly $c_4$-bounded
and disjoint from $M^*(S) \cup A_{|\Phi|}$.

\medskip
 
Recall that the cancelling pairs are
$(\{(e_i,Q_i^+),(e_i,Q_i^-): i \in [P''])$
with each $(Q_i^+ \cup Q_i^-) \sm \{e_i\} \sub G^*$.
By definition of the Splitting Phase we can write each
$Q_i^\pm = \phi_i^\pm(Q)$ where $\phi^\pm$ are $\pi$-compatible,
so Lemma \ref{elim:cascade} can be applied.
We adopt the notation of Definitions 
\ref{def:elim} and \ref{def:elim2}
and write $\OO^*{}' = \OO^* \sm (Q_i^+ \cup Q_i^-)$.

\begin{alg} (Elimination Phase)
We choose $\phi^*_i \in X^c_{E(\phi_i^\pm)}(G^*)$
by a random greedy algorithm.
Write $A''_i = \cup_{i'<i} M^*(\phi^*_{i'}(\OO^*{}'))$.
We choose $\phi^*_i \in X^c_{E(\phi_i^\pm)}(G^*)$
uniformly at random subject to $\phi^*_i(\OO^*{}') \cap 
(A_{|\Phi|} \cup A'_{P'} \cup M^*(S) \cup A''_i) = \es$.
\end{alg} 

Similarly to a remark above on the Solo Phase,
we note in each step of the Elimination Phase
that $M^*(\phi^*_i(\OO^*{}'))$ 
is edge-disjoint from 
$A_{|\Phi|} \cup A'_{P'} \cup M^*(S) \cup A''_i$,
as the latter is a union of cliques in $M^*$.

\begin{lemma} \label{absorbII}
W.h.p.\ the Elimination Phase does not abort
and $A''_{P''}$ is $c_4$-bounded
and linearly $c_4$-bounded.
\end{lemma}

\nib{Proof.} 
For $i \in [P'']$ we let $\mc{B}_i$ be the bad event 
that $A''_i$ is not $c_4$-bounded
or not linearly $c_4$-bounded. Let $\tau$ be the 
smallest $i$ for which $\mc{B}_i$ holds or the algorithm 
aborts, or $\infty$ if there is no such $i$. 
It suffices to show w.h.p.\ $\tau=\infty$. 
We fix $i_0 \in [P'']$ and bound 
$\mb{P}(\tau=i_0)$ as follows. 

By Lemma \ref{elim:cascade} we have
$X^c_{E(\phi_i^\pm)}(G^*) > 0.9 \oO^{6Q^2 p^r} n^{2pq-3q+r}$.
At most half of the choices of 
$\phi^*_i \in X^c_{E(\phi^\pm)}(G^*)$
are forbidden, as $M^*(S)$ is $c_2$-bounded,
$A_{|\Phi|}$ is $c_3$-bounded,
$A'_{P'}$ is $c_4$-bounded
and $A''_i$ is $c_4$-bounded.

For each $e \in G^*$ let 
$r_e = \sum_{i<i_0} \mb{P}'(e \in M^*(\phi^*_i(\OO^*{}')))
= \sum_{i<i_0} \sum_{e' \in M^*(e)}
\mb{P}'(e' \in \phi^*_i(\OO^*{}'))$.
Given $r' \in [r]$, as $A_{|\Phi|} \cup M^*(S)$ is $2c_3$-bounded
there are at most $4\tbinom{r}{r'}c_3 n^{r'}$ choices of $i$
such that $|e' \sm V(Q^+_i)|=r'$ or $|e' \sm V(Q^-_i)|=r'$.
For each such $i$ we have
$\mb{P}'(e' \in \phi^*_i(\OO')) 
< 1.2 |\OO^*{}'| \oO^{-6Q^2 p^r} n^{-r'}$, so
$r_e < Q|\OO^*{}'| \oO^{-6Q^2 p^r} 2^{r+3} c_3$.
As in the proof of Lemma \ref{split},
by Lemma \ref{dom} w.h.p.\ $A''_i$ is $c_4$-bounded
and linearly $c_4$-bounded
for all $i<i_0$, as required. \qed

\medskip

We let $\Phi^* = \Phi'' + \sum_{i \in [P'']} 
(\phi^*_i(\Ups^-)-\phi^*_i(\Ups^+))$.
Then $\pl \Phi^* = S$, and all cancelling pairs are cancelled,
as by Definitions \ref{def:elim} and \ref{def:elim2}
each $\phi^+_i = \phi^*_i f_1 \in \phi^*_i(\Ups^+)$
and $\phi^-_i = \phi^*_i f_2 \in \phi^*_i(\Ups^-)$.
For each positive clique $Q'$ added during Elimination Phase
we have $Q'$ cascading, $M^*(Q')$ is a set
and all such $M^*(Q')$ are disjoint.
These positive cliques include all $M^*(e')$ 
with $e' \in Q_i^+ \sm \{e_i\}$ 
for some cancelling pair $\{(e_i,Q_i^+),(e_i,Q_i^-)\}$;
we note that $M^*(e') \sub A_{|\Phi|}$.
The others are contained in $A''_{P''}$,
which is $c_4$-bounded, linearly $c_4$-bounded 
and disjoint from $M^*(S) \cup A_{|\Phi|} \cup A'_{P'}$.
This concludes the proof of Lemma \ref{hole}. 

\section{Conclusion} \label{sec:conclude}

We conclude with the proof of our main theorem and some remarks 
on possible future directions for research. For convenient reference,
we start by summarising the overall structure of the proof and identifying its key steps, 
all but one of which have been completed in the preceding sections.
We consider $G$ as in the statement of Theorem \ref{main+}, 
so $G$ is a $K^r_q$-divisible $(K^r_q,c,\oO)$-regular 
$(\oO,h)$-extendable $r$-multigraph on $n$ vertices.
\begin{enumerate}
\item In Section \ref{sec:template} we constructed the template $G^* \subset G$,
which is an edge-disjoint union of cliques $M^*$, most of which will appear in the final 
decomposition of $G$, but some of which will be replaced during later steps.
\item In Lemma \ref{nibble} we applied the nibble to choose a set of edge-disjoint 
cliques $M^n$ in $G \sm G^*$ that cover almost all of $G \sm G^*$; this was possible 
because of the regularity assumption on $G$ and the random construction of $G^*$.
We then chose a set of edge-disjoint cliques $M^c$ in Lemma \ref{cover} that covered the `leave' 
$L = G \sm \bigcup (M^* \cup M^n)$ and also some `spill' $S := \bigcup M^c \cap \bigcup M^*$;
this was possible by the combinatorial extendability properties of $G^*$ established in Lemma \ref{ext*}.
\item The remainder of the proof consists of creating a `hole' in the template which exactly matches the spill $S$,
so that these edges, which are currently covered twice (by the template and the cover) will be covered exactly once,
as required for a decomposition. The construction of the hole has three steps, two of which have been implemented:
Lemma \ref{bddint} provides an integral decomposition of $S$, which is modified in Lemma \ref{hole} via the
Clique Exchange Algorithm to a signed decomposition, namely two sets $M^+$ (the `positive' cliques) and $M^-$ (the `negative' cliques) 
of edge-disjoint cliques such that $\bigcup M^+$ is the disjoint union of $\bigcup M^-$ and $S$.
\item The final step of the proof, to be implemented in Lemma \ref{hole+} below, is to complete the construction of the hole by modifying
the signed decomposition $M^\pm$ into another signed decomposition, consisting of $M^o$ (the `outer' cliques)
and $M^i$ (the `inner' cliques), where similarly to before $\bigcup M^o$ is the disjoint union of $\bigcup M^i$ and $S$,
but now we have the additional property $M^o \sub M^*$, that is, each outer clique is a template clique. This is the key
step of the proof, as we can remove $M^o$ from $M^*$ and replace it by $M^i$ to exactly correct for the spill.
\end{enumerate}

We now implement the final piece of the argument,
via a cascade algorithm for absorption,
which we formulate as a separate lemma as it 
may be useful for other problems.
The proof applies a random greedy algorithm to choose 
edge-disjoint cascades for each positive clique,
which is possible as we ensured in Lemma \ref{hole}
that each positive clique is cascading 
and every template clique shares an edge
with at most one positive clique.

\begin{lemma} \label{hole+}
Suppose $S \sub G^*$ is $K^r_q$-divisible
and $M^*(S)$ is a set that is $c_2$-bounded
and linearly $c_2$-bounded. 
Then there is $M^o \sub M^*$ and $M^i \sub K^r_q(G^*)$
such that $\sum M^o$ is the disjoint union
of $\sum M^i$ and $S$.
\end{lemma}

\nib{Proof.}
By Lemma \ref{hole} we can choose $M^\pm$
such that every clique in $M^+$ is cascading,
$M^*(\sum M^+)$ is a set and linearly $3c_4$-bounded,
and $\sum M^+$ is the disjoint union
of $\sum M^-$ and $S$.

We apply a random greedy algorithm
to choose cascades for each clique in $M^+$.
Write $M^+ = (Q^i: i \in [P])$.
At step $i$ we choose a cascade $C^i=C^{Q^i}$ 
for $Q^i$ and write $C'{}^i = C^i \sm M^*(Q^i)$.
We choose $C^i$ uniformly at random such that
$C'{}^i$ is disjoint from $M^*(\sum M^+)$
and $A^i := \cup_{i'<i} C'{}^{i'}$.

For $i \in [P]$ we let $\mc{B}_i$ be the bad event 
that $A^i$ is not linearly $c_5$-bounded. Let $\tau$ be the 
smallest $i$ for which $\mc{B}_i$ holds or the algorithm 
aborts, or $\infty$ if there is no such $i$. 
It suffices to show w.h.p.\ $\tau=\infty$. 
We fix $i_0 \in [P]$ and bound $\mb{P}(\tau=i_0)$.

For any $i<i_0$, we claim that at most half of the
choices of $C^i$ are forbidden by the disjointness condition.
To see this, fix any $\phi' \in \Ups'$, $I \in Q$, $a \in Ker^I$,
where we can assume $r_{a\phi'}>0$
by Remark \ref{rem:nondegen},
and consider which cascades
are forbidden due to $e^w_a \in M^*(\sum M^+) \cup A^i$,
where $w = w_{\phi'}$.
Given $Q^i \sub G^*_j$, by Lemma \ref{lem:ee'}, the possible 
$f_j(e^w_a)$ lie an affine linear space of dimension $r_{a\phi'}$.
As $M^*(\sum M^+)$ is linearly $3c_4$-bounded
and $A^i$ is linearly $c_5$-bounded 
($\mc{B}_i$ does not hold), we have at most 
$2c_5 (p^a)^{r_{a\phi'}}$ choices 
of $e^w_a \in M^*(\sum M^+) \cup A^i$,
and by Lemma \ref{cascade:role} each forbids at most 
$(p^a)^{qp-q-r_{a\phi'}}$ choices of $C^i$.
Summing over $\phi'$ and $a$, at most
$2c_5 p^{q^2} (p^a)^{qp-q}$ choices of $C^i$ are forbidden,
which by Lemma \ref{cascade} is at most half the total, as claimed.

For each $e \in G^* \sm M^*(\sum M^+)$ let 
$r_e = \sum_{i<i_0} \mb{P}'(e \in C'{}^i)$,
where $\mb{P}'$ denotes conditional probability
given the choices made before step $i$.
We consider the contributions to $r_e$ from cascades 
choosing $e=e^w_a$ where $w = w_{\phi'}$ 
for any fixed $\phi' \in \Ups'$, $I \in Q$, $a \in Ker^I$,
where as $e \notin M^*(\sum M^+)$
we can assume $r_{a\phi'}>0$
by Remark \ref{rem:nondegen}.
Let $J = \{i \in [q]: \phi'(i)=(i,1)=i\}$
and fix $J \sub I' \in Q$. Suppose $e \in G^*_j$.
By Lemma \ref{lem:ee'}, if a cascade $C_{\phi^c}$
for some $Q^i=\phi^c([q])$ such that
$A^{\phi^c\phi'(Q)}=\phi^{w_{\phi'}}(K^r_q(Ker))$
has $e=e^{w_{\phi'}}_a$ then we have $\phi^c(I')=e'$
for some edge $e'$ such that $f_j(e')$ lies in an
affine space of dimension $r_{a\phi'}$.
As $M^*(\sum M^+)$ is linearly $3c_4$-bounded,
there are at most $3c_4 (p^a)^{r_{a\phi'}}$
choices for $e'$, which determines $Q^i$.
For each such $Q^i$, by Lemma \ref{cascade:role} 
there are at most $(p^a)^{q(p-1)-r_{a\phi'}}$ choices of $C^i$ 
such that $e^{w_{\phi'}}_a=e$. We choose $C^i$ with probability
at most $2 (p^a)^{q(p-1)-r_{a\phi'}}
/ \oO^{p^{q^2}} n^{q(p-1)}$, by Lemma \ref{cascade} 
and the bound on excluded choices.
Summing over all $\phi'$, $a$ and $i$, we get
$r_e < p^{q^2} 3c_4 (p^a)^{r_{a\phi'}}
2 (p^a)^{q(p-1)-r_{a\phi'}}
/ \oO^{p^{q^2}} n^{q(p-1)}
< c_4 p^{q^2} \gG^{-pq} \oO^{-p^{q^2}} < c_5/2$.

As in the proof of Lemma \ref{split},
by Lemma \ref{dom} w.h.p.\ $\tau=\infty$,
so the algorithm does not abort,
and we can choose cascades $C^i=C^{Q^i}$ for all $i \in [P]$.
Then $M^o = \bigcup_{Q' \in M^+} \Psi(C^{Q'}) \sub M^*$
and $M^i = M^- \cup \bigcup_{Q' \in M^+} (\Psi'(C^{Q'})\sm\{Q'\})$
are as required (in words, $M^o$ contains the template
decomposition of the cascade of each positive clique;
$M^i$ is obtained by flipping these cascades 
and then replacing $M^+$ by $M^-$). \qed

\medskip

To complete the proof of our main theorem,
we take the matchings $M^*$, $M^n$ and $M^c$
obtained from the template, nibble and cover,
then use the previous lemma to correct for the spill $S$.

\medskip

\nib{Proof of Theorem \ref{main+}.}
We choose a template $M^*$ and $G^*$ 
as in Definition \ref{def:template} that 
satisfies all of the w.h.p.\ statements in the paper.
We let $M^n$ be obtained from Lemma \ref{nibble}
and $M^c$ and $S$ from Lemma \ref{cover}.
Note that $S$ is $K^r_q$-divisible, 
as $S = \sum M^* + \sum M^n + \sum M^c - G$,
and $M^*(S)$ is a set that is $c_2$-bounded
and linearly $c_2$-bounded, so we can apply
Lemma \ref{hole+} to obtain $M^o \sub M^*$ and $M^i \sub K^r_q(G^*)$
such that $\sum M^o$ is the disjoint union
of $\sum M^i$ and $S$.
Our final $K^r_q$-decomposition of $G$ 
is $M = M^n \cup M^c \cup (M^* \sm M^o) \cup M^i$. \qed

\medskip

\nib{Remarks.}
As discussed in Section \ref{sec:subseq}, there has been an explosion of progress 
in the new probabilistic approach to Design Theory 
in the decade following the first arXiv version \cite{K} of this paper.
These new results suggest many further questions
(including several that can be found in the relevant papers),
so here we just include our original remarks from 2014 
on potential directions for future research,
which are still instructive, although (perhaps unsurprisingly)
they fail to capture the richness and potential further applications of these ideas.

These original remarks discussed two main directions for further research.
One direction was structural characterisations 
of the perfect matching problem in hypergraphs, 
noting the analogies between the sparse setting of auxiliary hypergraphs
for designs and the dense setting that we characterised in earlier work 
with Mycroft \cite{KM} and Knox and Mycroft \cite{KKM} 
(see also the survey by R\"odl and Ruci\'nski \cite{RR} for many further references).
In this context we identified Ryser's Conjecture on transversals in Latin Squares 
as a target for further progress, which turned out to be a good prediction!
It will be interesting to investigate the scope of Montgomery's approach for such problems,
as discussed in his survey \cite{Msurvey}.

The other direction was to pursue potential connections with the study of
random matchings in Probability and Statistical Physics along the lines of results 
obtained by Kahn \cite{KaICM,Ka2}, Kahn and Kayll \cite{KaKa} 
and Barvinok and Samorodnitsky \cite{BS}. This still seems very hard,
and a recent disproof by Lee \cite{Lee} of one of Kahn's conjectures
indicates that we do not yet even know what we should be trying to prove in general.
As discussed above, there has been progress for specific models, 
particularly random Latin Squares / Steiner Systems, 
although even here we still cannot answer basic questions
such as asymptotic counting of these structures, so we still 
seem quite far from a detailed probabilistic description of these models.

\medskip

\renewcommand{\theenumi}{\alph{enumi}}

\nib{Acknowledgements.} 
I would like to thank: 
\begin{enumerate}
\item Rick Wilson for encouraging me to work on the Existence Conjecture,
\item Stefan Glock, Daniela K\"uhn, Allan Lo and Deryk Osthus
for pointing out an error in the first version of this paper,
\item Eoin Long and Lisa Sauermann for carefully reading the second version
and many comments that have improved the presentation,
\item Lisa Sauermann for pointing out an incorrect `simplification'
of the proof in going from the first to the second version
(see Remark \ref{rem:linbdd}).
\item The anonymous referees for helpful comments.
\end{enumerate}

\end{document}